\documentclass[11pt]{smfart}

\usepackage{amssymb, amsxtra, amsmath, amsfonts}

\usepackage{smfthm}
\theoremstyle{plain}

\usepackage[french,english]{babel}

\makeatletter
\def\amsbb{\use@mathgroup \M@U \symAMSb}
\makeatother

\usepackage[all]{xy}
\SelectTips{cm}{} 

\usepackage{tikz}
\usepackage{tikz-cd}
\usetikzlibrary{matrix,arrows,backgrounds}

\usepackage[a4paper, marginparwidth=60pt, margin=1.0in]{geometry}

\usepackage{mathrsfs}

\usepackage{marginnote}

\usepackage{stmaryrd}

\usepackage[outline]{contour}

\usepackage{enumitem}

\usepackage{tensor}

\usepackage{mathtools}

\usepackage{hyperref}

\numberwithin{equation}{section}

\setitemize[0]{leftmargin=*, itemsep=\the\smallskipamount}
\setenumerate[0]{leftmargin=*, itemsep=\the\smallskipamount}

\providecommand{\amsbb}{\mathbb}
\newcommand{\BQ}{\amsbb{Q}}
\newcommand{\BZ}{\amsbb{Z}}
\newcommand{\BL}{\amsbb{L}}
\newcommand{\BN}{\amsbb{N}}
\newcommand{\BP}{\amsbb{P}}
\newcommand{\BV}{\amsbb{V}}
\newcommand{\BG}{\amsbb{G}}
\newcommand{\BF}{\amsbb{F}}
\newcommand{\BD}{\amsbb{D}}
\newcommand{\BH}{\amsbb{H}}

\newcommand{\CB}{\mathcal{B}}
\newcommand{\CL}{\mathcal{L}}

\newcommand{\CO}{\mathcal{O}}

\newcommand{\CU}{\mathcal{U}}
\newcommand{\CV}{\mathcal{V}}

\newcommand{\CI}{\mathcal{I}}
\newcommand{\CJ}{\mathcal{J}}
\newcommand{\CT}{\mathcal{T}}

\newcommand{\sC}{\mathscr{C}}
\newcommand{\sD}{\mathscr{D}}
\newcommand{\sG}{\mathscr{G}}
\newcommand{\sL}{\mathscr{L}}

\newcommand{\sU}{\mathscr{U}}
\newcommand{\sW}{\mathscr{W}}

\newcommand{\bPi}{\mathbf{\Pi}}
\newcommand{\bV}{\mathbf{V}}

\newcommand{\bD}{\mathbf{D}}
\newcommand{\bm}{\mathbf{m}}
\newcommand{\bn}{\mathbf{n}}

\newcommand{\fka}{\mathfrak{a}}
\newcommand{\fkm}{\mathfrak{m}}
\newcommand{\fkg}{\mathfrak{g}}
\newcommand{\fkn}{\mathfrak{n}}
\newcommand{\fkr}{\mathfrak{r}}

\newcommand{\Qp}{\BQ_p}
\newcommand{\Qpt}{\Qp^{\times}}
\newcommand{\Qpbar}{\overline{\BQ}_p}
\newcommand{\Qpd}{\BQ_{p^2}}
\newcommand{\Zp}{\BZ_p}

\newcommand{\czG}{\cz{G}}
\newcommand{\GL}{\mathrm{GL}}
\newcommand{\Gal}{\mathrm{Gal}}
\newcommand{\rmM}{\mathrm{M}}
\newcommand{\rmO}{\mathrm{O}}
\newcommand{\rmD}{\mathrm{D}}
\newcommand{\rmOD}{\rmO_{\rmD}}
\newcommand{\rmODt}{\rmO_{\rmD}^{\times}}
\newcommand{\rmU}{\mathrm{U}}

\newcommand{\rmH}{\mathrm{H}}

\newcommand{\Harg}[2]{\!\big(#1\,;\,#2\big)}

\newcommand{\Bdr}{\mathrm{B}_{\dR}}
\newcommand{\Bdrp}{\mathrm{B}_{\dR}^{+}}
\newcommand{\Bst}{\mathrm{B}_{\st}}
\newcommand{\Bstp}{\mathrm{B}_{\st}^{+}}
\newcommand{\Bcrisp}{\mathrm{B}_{\cris}^{+}}

\newcommand{\pet}{\mathrm{p}\acute{\mathrm{e}}\mathrm{t}}
\newcommand{\et}{\acute{\mathrm{e}}\mathrm{t}}
\newcommand{\dR}{\mathrm{dR}}
\newcommand{\st}{\mathrm{st}}
\newcommand{\cris}{\mathrm{cr}}
\newcommand{\HT}{\mathrm{HT}}
\newcommand{\pst}{\mathrm{pst}}
\newcommand{\lalg}{\mathrm{lalg}}
\newcommand{\Fil}{\mathrm{Fil}}

\newcommand{\JL}{\mathrm{JL}}
\newcommand{\LL}{\mathrm{LL}}
\newcommand{\WD}{\mathrm{WD}}
\newcommand{\sWD}{\sW\!\sD}
\newcommand{\Sp}{\mathrm{Sp}}
\newcommand{\St}{\mathrm{St}}
\newcommand{\Tors}{\mathrm{Tors}}
\newcommand{\Irr}{\mathrm{Irr}}
\newcommand{\Ext}{\mathrm{Ext}}
\newcommand{\ps}{\mathrm{ps}}
\newcommand{\tr}{\mathrm{tr}}
\newcommand{\Stab}{\mathrm{Stab}}
\newcommand{\Ban}{\mathrm{Ban}}
\newcommand{\Gr}{\mathrm{Gr}}
\newcommand{\niv}{\mathrm{niv}}
\newcommand{\ind}{\mathrm{ind}}
\newcommand{\lan}{\mathrm{lan}}
\newcommand{\finie}{\mathrm{finie}}

\newcommand{\Hom}{\mathrm{Hom}}
\newcommand{\End}{\mathrm{End}}
\newcommand{\Ind}{\mathrm{Ind}}
\newcommand{\diag}{\mathrm{diag}}
\newcommand{\Sym}{\mathrm{Sym}}
\newcommand{\id}{\mathrm{id}}
\newcommand{\dif}{\mathrm{dif}}
\newcommand{\nrd}{\mathrm{nrd}}
\newcommand{\nr}{\mathrm{nr}}

\DeclareMathOperator{\Coker}{Coker}

\newcommand{\intnn}[2]{\!\big(#1,\,#2\big)}
\newcommand{\intn}[1]{\big(#1\big)}
\newcommand{\restr}[2]{{\left.\kern-\nulldelimiterspace#1\vphantom{\big|}\right|_{#2}}}

\newcommand{\wh}{\widehat}
\newcommand{\wt}{\widetilde}
\newcommand{\cz}[1]{\check{#1}}
\newcommand{\brv}[1]{\breve{#1}}
\newcommand{\wotimes}{\widehat{\otimes}}

\newcommand{\Spec}{\mathrm{Spec}}
\newcommand{\Spm}{\mathrm{Spm}}
\newcommand{\an}{\mathrm{an}}
\newcommand{\presp}[1]{\prescript{p}{}{#1}}

\newcommand{\itemb}{\item[$\bullet$]}
\newcommand{\cf}{{\it cf.}\ }
\newcommand{\ie}{{\it i.e.}\ }
\newcommand{\resp}{{\it resp.}\ }

\newcommand{\cfc}[1]{(\cf \cite{#1})}
\newcommand{\num}{n\textsuperscript{o}\,}

\newcommand{\rec}{\mathrm{rec}}
\newcommand{\phinl}{\Phi N_{\mkern-1mu{\lambda}}}
\newcommand{\wtphinl}{\wt{\Phi N}_{\mkern-1mu{\lambda}}}
\newcommand{\czW}{\check{W}}
\newcommand{\rmb}{\mathrm{b}}
\newcommand{\rmB}{\mathrm{B}}
\newcommand{\sfs}{\sf{s}}
\newcommand{\sft}{\sf{t}}
\newcommand{\incl}{\hookrightarrow}
\newcommand{\surj}{\twoheadrightarrow}
\newcommand{\sesi}{\mathrm{ss}}

\title[Factorisation de la cohomologie de systèmes locaux $p$-adiques]{Factorisation de la cohomologie de systèmes locaux $p$-adiques sur le demi-plan de Drinfeld}
\author{Arnaud Vanhaecke}
\address{{\tiny Morningside Center of Mathematics, No. 55, Zhongguancun East Road, Beijing, 100190, China.} }
\email{arnaud@amss.ac.cn}
\date{\today}
\thanks{Ce travail a été financé par le Morningside Center of Mathematics, CAS}

\begin{document}

\maketitle

\begin{abstract}
On calcule le premier groupe de cohomologie de l'algèbre symétrique du système local étale $p$-adique universel sur la tour de revêtement du demi-plan $p$-adique de Drinfeld. Le résultat s'énonce sous une forme factorisée en utilisant la correspondance de Langlands $p$-adique en famille sur les anneaux de Kisin. Ce travail étend les résultats correspondants de Colmez, Dospinescu et Nizio{\l} pour les coefficients triviaux. Il repose sur le calcul des multiplicités automorphes dans le groupe de cohomologie étale du système local, obtenu dans un article antérieur, et la détermination des anneaux de Kisin pour le type spécial comme fonctions sur un ouvert analytique de la droite projective.
\end{abstract}
\begin{altabstract}
We compute the first cohomology group of the symmetric algebra of the universal étale $p$-adic local system on the tower of coverings of Drinfeld’s $p$-adic half-plane. The result takes a factorized form, using the $p$-adic Langlands correspondence in families over Kisin rings. This work extends the corresponding results of Colmez, Dospinescu, and Nizio{\l} for trivial coefficients. It relies on the computation of automorphic multiplicities in the étale cohomology group of the local system, done in a previous paper, as well as on the determination of the Kisin rings for the special type as functions on an analytic open subset of the projective line.
\end{altabstract}

\setcounter{tocdepth}{2}
\tableofcontents

\section{Introduction}

Soit $p$ un nombre premier, $\sG_{\Qp}\coloneqq \Gal(\Qpbar/\Qp)$ le groupe de Galois absolu de $\Qp$, $G=\GL_2(\Qp)$ et $\czG\coloneqq \rmD^{\times}$ le groupe des unités de $\rmD$, l'algèbre des quaternions non scindée sur $\Qp$. Soit $C$ le complété d'une clôture algébrique $\Qpbar$ de $\Qp$, et soit $L$ une extension finie, assez grande, de $\Qp$ qui est notre corps des coefficients. Dans \cite{van}, nous avons calculé la multiplicité des représentations galoisiennes dans la cohomologie (c'est-à-dire, dans le premier groupe de cohomologie) de l'algèbre symétrique du système local universel, notée $\Sym \BV(1)$, sur la tour de revêtements de Drinfeld $\{\brv{\rmM}_C^{n}\}_{n\geqslant 0}$, en termes de la correspondance de Langlands $p$-adique. 
 
Ce travail étendait des résultats de Colmez Dospinescu et Nizio{\l} (\cf \cite{codoni}) qui traitent le cas des coefficients triviaux $L(1)$. Les deux différences principales sont les suivantes :
\begin{itemize}
 \itemb Alors que la cohomologie étale de $\{\brv{\rmM}_C^{n}\}_{n\geqslant 0}$ à coefficients dans $L(1)$ fait apparaître les représentations potentiellement cristallines à poids de Hodge-Tate $(0,1)$, la cohomologie à coefficients dans $\Sym \BV(1)$ fait apparaître toutes les représentations potentiellement cristallines à poids de Hodge-Tate réguliers positifs.
 \itemb Une seconde différence, plus importante, est que la cohomologie étale de $\brv{\rmM}_C^{0}$ à coefficients dans $L(1)$ est simplement le dual d'une steinberg continue, alors qu'à coefficients dans $\Sym \BV(1)$, elle fait apparaître toutes les représentations semi-stables non cristallines (et même une représentation cristalline dite \emph{exceptionnelle}), avec la multiplicité attendue (\ie donnée par la correspondance de Langlands $p$-adique).
\end{itemize}
À la suite de \cite{codoni} et pour $p>3$, Colmez Dospinescu et Nizio{\l} déterminent la structure complète de la cohomologie  étale de\footnote{Deux petites différences avec le cas précédent sont de remplacer $C$ par $\Qpbar$ dans les coefficients de l'espace de Drinfeld et de considérer $\presp{\rmM}_{\Qpbar}^{n}=\brv{\rmM}_{\Qpbar}^{n}/p$, où $p$ est vu comme un élément de $G$.} $\{\presp{\rmM}_{\Qpbar}^{n}\}_{n\geqslant 0}$ à coefficients dans $L(1)$ comme $L[G\times \sG_{\Qp}\times \czG]$-module (\cf \cite{codonifac}). Le résultat s'énonce sous une forme factorisée, similaire à la description d'Emerton (\cf \cite{eme}) dans le cas de la cohomologie complétée des courbes modulaires, et fait intervenir les anneaux de déformation potentiellement cristallins introduits par Kisin (\cf \cite{kis}), dont certaines propriétés essentielles (en particulier que se sont des produits finis d'anneaux principaux) ont été établies dans \cite{codonikis}.

\subsection{Enoncés des résultats}
Le but de cet article est de poursuivre cette entreprise, à partir des résultats de \cite{van}, pour obtenir une description similaire de la cohomologie de $\presp{\rmM}_{\Qpbar}^{n}$ à coefficients dans $\Sym \BV(1)$, comme $L[G\times \sG_{\Qp}\times \czG]$-module, sous une forme factorisée :

\medskip
\begin{theo}\label{thm:princint}
Pour $p>3$, on a un isomorphisme de $L[G\times \sG_{\Qp}\times \czG]$-modules topologiques :
$$
\rmH^1_{\et}\Harg{\presp{\rmM}_{\Qpbar}^{n}}{\Sym \BV(1)}\cong \bigoplus_{\lambda}\bigoplus_{M}\bigl (\wh{\bigoplus_{\CB}}\ \bPi(\rho^{\lambda}_{\CB,M})'\otimes \rho_{\CB,M}^{\lambda}\otimes\cz{R}_{\CB,M}^{\lambda} \bigr ) \otimes_L \JL_M^{\lambda}.
$$
\end{theo}
\medskip
Expliquons brièvement les termes de cet énoncé, qui seront détaillés dans la suite.
\begin{itemize}
\itemb La première somme porte sur les \emph{poids de Hodge-Tate} $\lambda=(\lambda_1,\lambda_2)\in \BN^2$ tels que $\lambda_2>\lambda_1\geqslant 0$. La seconde somme porte sur les $L$-$(\varphi,N,\sG_{\Qp})$-modules spéciaux et cuspidaux $M$ compatibles à $\lambda$ et de niveau $\leqslant n$. La troisième somme, complétée pour la topologie $p$-adique, porte sur les blocs $\CB$ des représentations de $G$ modulo $p$ que l'on identifie (\cf \cite{pas2}) aux représentations $\bar\rho_{\CB}$ modulo $p$ de $\sG_{\Qp}$ semi-simples et de dimension $2$.
\itemb Les produits tensoriels sont au-dessus de l'anneau de Kisin $R_{\CB,M}^{\lambda}$ qui est l'anneau de déformation de $\bar\rho_{\CB}$ des représentations de type $M$ et à poids de Hodge-Tate $\lambda$ ; $\cz{R}_{\CB,M}^{\lambda}$ est le $L$-dual continu de $R_{\CB,M}^{\lambda}$. Le $R_{\CB,M}^{\lambda}[\sG_{\Qp}]$-module $\rho_{\CB,M}^{\lambda}$ est la déformation universelle et $\bPi(\rho_{\CB,M}^{\lambda})'$ est le dual stéréotypique du $R_{\CB,M}^{\lambda}[G]$-module associé à $\rho_{\CB,M}^{\lambda}$ par la correspondance de Langlands $p$-adique.
\itemb Le dernier terme $\JL_M^{\lambda}$ est la représentation localement algébrique de $\czG$ de poids $(\lambda_1,\lambda_2-1)$ associée à $M$ par la correspondance de Jacquet-Langlands.
\end{itemize}

Pour $M$ cuspidal (apparaissant dès que $n\geqslant 1$), le résultat est exactement le même que celui de \cite{codonifac}, et fait intervenir les poids supérieurs. Il suffit alors, en utilisant les résultats de \cite{codonikis} et les calculs de \cite{van}, de reprendre pas à pas la démonstration de \cite{codonifac}. L'un des points les plus délicat de de \cite{codonifac} est la démonstration de propriétés de finitude la cohomologie étale modulo $p$ ; à coefficients, on peut directement utiliser ces résultats puisque le système local en question modulo $p$ devient trivial sur un revêtement fini, puis utiliser la suite spectrale d'Hochschild-Serre pour obtenir ces propriétés pour tous les étages.

Pour $M$ spécial (c'est-à-dire, $n=0$ à torsion par un caractère près), alors que le résultat à coefficients dans $L(1)$ est simplement que la cohomologie étale du demi-plan de Drinfeld est le dual d'une steinberg continue, il est bien plus riche à coefficients non triviaux. La preuve de ce cas constitue le noyau dur du présent travail. En particulier, pour démontrer le théorème \ref{thm:princint} dans le cas spécial, on doit adapter les résultats de \cite{codonikis} aux anneaux de Kisin de type spécial ; le théorème utilisé, obtenu dans le chapitre $2$, est le suivant :
\medskip
\begin{theo}\label{thm:princint2}
Supposons $p>3$, si $\lambda=(\lambda_1,\lambda_2)\in \BN^2$ sont tels que $\lambda_2-\lambda_1>1$ et $M$ un $L$-$(\varphi,N,\sG_{\Qp})$-module de rang $2$ tel que $N\neq 0$, alors l'anneau de Kisin $R_{\CB,M}^{\lambda}$ est l'anneau des fonctions analytiques bornées sur un ouvert analytique de $\BP^{1,\an}_L$ ; de plus, c'est un produit fini d'anneaux principaux.
\end{theo}
\medskip


\subsection{Notations}
On suppose que $p>3$. Soient $\sW_{\Qp}\subset \sG_{\Qp}$ le groupe de Weil et $\sWD_{\Qp}$ le groupe de Weil-Deligne. On note $\rmOD\subset \rmD$ l'unique ordre maximal et $\varpi_{\rmD}\in \rmOD$ une uniformisante. On suppose que $L$ est assez grande, en particulier qu'elle contient toutes les extensions quadratiques de $\Qp$. Soient $\rmO_L\subset L$ l'anneau des entiers, $\fkm_L\subset \rmO_L$ son idéal maximal et $k_L\coloneqq \rmO_L/\fkm_L$ son corps résiduel. On fixe $\varpi_{L}\in \fkm_L$ une uniformisante. Pour $A$ un anneau (topologique), une \emph{$A$-représentation (continue)} désigne une représentation sur un $A$-module (topologique). De plus, tous les $\Hom$ entre $A$-représentations topologiques seront continus et munis de la topologie de la convergence uniforme sur les parties fortement bornées.

Soit $\rec\colon \sW_{\Qp}\rightarrow \Qpt$ le morphisme de réciprocité, normalisé de sorte que le frobenius arithmétique soit envoyé sur $p$, $\det \colon G\rightarrow \Qpt$ le déterminant et $\nrd \colon \czG\rightarrow \Qpt$ la norme réduite. Au travers de ces caractères tout $L$-caractère (\ie $L$-représentation de dimension $1$) de $\Qpt$ définit un $L$-caractère de $\sW_{\Qp}$, $G$ et $\czG$ et s'il est unitaire (\ie à coefficients dans $\rmO_L^{\times}$), de $\sG_{\Qp}$ ; on identifie ainsi les caractères de ces groupes. On note $\omega\colon \Qpt\rightarrow \rmO_L^{\times}$ le \emph{caractère cyclotomique} défini pour $x\in \Qpt$ par $\omega(x)=x\lvert x\rvert_p$.
\subsubsection{Types spéciaux et cuspidaux}
On fixe une injection $\rmD\incl M_2(L)$ qui définit une injection $\czG\incl \GL_2(L)$, ce qui est possible puisque $L$ contient une extension quadratique de $\Qp$. On note 
$$
P\coloneqq \{(\lambda_1,\lambda_2)\in \BN^2\mid \lambda_2>\lambda_1\},
$$
l'ensemble des \emph{poids de Hodge-Tate réguliers} et $P_+\subset P$ l'ensemble des poids positifs, \ie tels que $\lambda_1\geqslant 0$. À $\lambda\in P$ on associe $w(\lambda)\coloneqq \lambda_2-\lambda_1\in \BN$ et $\lvert \lambda \rvert\coloneqq \lambda_2+\lambda_1\in \BZ$ et on dit que $\lambda\in P$ est \emph{très régulier} si $w(\lambda)>1$. Si $V$ est une $L$-représentation de Rham de $\sG_{\Qp}$ on dit qu'elle est à poids de Hodge-Tate $\lambda=(\lambda_1,\lambda_2)$ si ses poids de Hodge-Tate sont $\lambda_1$ et $\lambda_2$.

Un $L$-$(\varphi,N,\sG_{\Qp})$-module est un $L\otimes_{\Qp}\Qp^{\nr}$-module $M$ muni d'un frobenius $\varphi$, d'un opérateur de monodromie $N$ tel que $N\varphi=p\varphi N$ et d'une action lisse de $\sG_{\Qp}$ commutant à $\varphi$ et $N$. On associe à $M$ une $L$-représentation de Weil-Deligne $\WD(M)$, \ie une représentation du groupe $\sWD_{\Qp}$. Si $M$ est de rang $2$, on dit que $M$ est
\begin{itemize}
\itemb \emph{cuspidal} si $\WD(M)$ est absolument irréductible ; dans ce cas $N=0$,
\itemb \emph{exceptionnel} si $\WD(M)\cong \chi_1\oplus \chi_2$ où $\chi_1,\chi_2\colon \sW_{\Qp}\rightarrow L^{\times}$ sont deux caractères lisses tels que $\chi_1\chi_2^{-1}=\lvert \rec \rvert_p^{\pm 1}$, où $\rec\colon \sW_{\Qp}\rightarrow \Qpt$ est le morphisme de réciprocité,
\itemb \emph{spécial} si $M$ est exceptionnel ou bien si $N\neq 0$ (notons que $N\neq 0$ implique que $\WD(M)^{\sesi}$ est exceptionnel).
\end{itemize}
Pour $\lambda\in P$, on note (\cf \cite[Définition 11.2]{van}) $\wtphinl$ (\resp $\phinl$) l'ensemble des $L$-$(\varphi,N,\sG_{\Qp})$-modules spéciaux\footnote{Pour $M\in \phinl$ non cuspidal on dira simplement que $M$ est spécial, sous-entendu non exceptionnel. La distinction entre $\wtphinl$ et $\phinl$ apparait surtout pour la fluidité de l'introduction. Dans la suite on utilisera surtout $\phinl$ puisque les représentations exceptionnelles apparaissent par déformation des représentations spéciales non exceptionnelles.}  (\resp non exceptionnel) ou cuspidaux dont la moyenne des pentes est $1-\lvert \lambda \rvert/2$ ($\phinl$ ne dépend donc que de $\lvert \lambda \rvert$). À $\lambda\in P$ et $M\in \wtphinl$ on associe des représentations localement algébriques unitaires de $G$ et $\czG$ :
\begin{itemize}
\itemb Soit $\LL_M$ la $L$-représentation lisse de $G$ obtenue en appliquant la correspondance de Langlands locale à $\WD(M)$. On définit une $L$-représentation unitaire localement algébrique de $G$ de poids $(\lambda_1,\lambda_2-1)$ :
$$
\LL_M^{\lambda}\coloneqq \LL_M\otimes_L \Sym_L^{w(\lambda)-1}\otimes_L{\det}^{\lambda_1},
$$
où $\Sym_L^{w(\lambda)-1}$ est la puissance $w(\lambda)-1=(\lambda_2-\lambda_1)-1$ symétrique de la représentation standard $L^2$ de $G$. On note $\zeta_M^{\lambda}$ le $L$-caractère central (unitaire) de $\LL_M^{\lambda}$.

Si $M$ est spécial non exceptionnel, alors $\LL_M^{\lambda}$ est la torsion par $\zeta_M^{\lambda}\circ \det$ de la \emph{steinberg localement algébrique de poids $(\lambda_1,\lambda_2-1)$}, définie par $\St_{\lambda}^{\lalg}=\St_L^{\infty}\otimes_L W_{\lambda}^*$, où
$$
W_{\lambda}^*\coloneqq \Sym_L^{w(\lambda)-1}\otimes_L{\det}^{\lambda_1}\otimes_L \lvert \det \rvert_p^{\tfrac{\lvert \lambda \rvert -1}{2}},
$$
et $\St_L^{\infty}$ est la steinberg lisse de $G$ \ie l'espace des fonctions localement constantes sur $\BP^1(\Qp)$ modulo les fonctions constantes. Notons que le dernier facteur de la définition de $W_{\lambda}^*$, où l'on rappelle que $\lvert \lambda \rvert=\lambda_2+\lambda_1$, est présent pour s'assurer que $\St_{\lambda}^{\lalg}$ est unitaire. Si $M$ est exceptionnel, $\LL_M^{\lambda}$ est\footnote{Ce n'est pas  la définition \og usuelle \fg\,  de la correspondance de Langlands locale (\cf \cite[VI 6. 11]{colp}).}  la torsion par $(\zeta_M^{\lambda}\circ\det)$ de l'unique extension localement algébrique de $W_{\lambda}^*$ par $\St_{\lambda}^{\lalg}$ que l'on note $\wt{\St}_{\lambda}^{\lalg}$.

\itemb Soit\footnote{Il faut ici supposer que $L$ est assez grand pour que $\JL_M$ soit défini sur $L$.} $\JL_M$ la $L$-représentation lisse de $\czG$ associée à $\LL_M$ par la correspondance de Jacquet-Langlands à $\JL_M$. On définit une $L$-représentation unitaire localement algébrique de $\czG$ de poids $(\lambda_1,\lambda_2-1)$ :
$$
\JL_M^{\lambda}=\JL_M\otimes_L\Sym_L^{w(\lambda)-1}\otimes_L{\nrd}^{\lambda_1},
$$
où $\Sym_L^{w(\lambda)-1}$ est vu comme représentation de $\czG$ au travers de l'injection $\czG\incl \GL_2(L)$. Le caractère central de $\JL_M^{\lambda}$ est le même que celui de $\LL_M^{\lambda}$, c'est-à-dire $\zeta_M^{\lambda}$.

Si $M$ est spécial, $\JL_M^{\lambda}\cong \czW_{\lambda}\otimes_L(\zeta_M^{\lambda}\circ \nrd)$ où 
\begin{equation}\label{eq:defczw}
\czW_{\lambda}\coloneqq \Sym_L^{w(\lambda)-1}\otimes_L{\nrd}^{\lambda_1}\otimes_L \lvert \nrd \rvert_p^{\tfrac{\lvert \lambda \rvert -1}{2}}.
\end{equation}
\end{itemize}

Côté galoisien, si $\lambda\in P$ et $M\in \wtphinl$ alors on dit que $V$, une $L$-représentation de Rham  de $\sG_{\Qp}$ de dimension $2$, est \emph{de type\footnote{\label{ftn:convtype}On n'a pas recours ici à l'abus de notation dans \cite{van} où le type d'une représentation exceptionnelle est donné par $M\in \phinl$.} $(M,\lambda)$} si elle est à poids de Hodge-Tate $\lambda$ et \footnote{Pour $D$ un $\varphi$-module et $k\in \BZ$ on note $D[k]$ le module $D$ muni du frobenius $p^k\varphi$.} $\bD_{\pst}(V)=M[-1]$. En particulier, on dit que $V$ est \emph{cuspidale} (\resp \emph{spéciale}, \emph{exceptionnelle}) si $M$ est cuspidal (\resp spécial, exceptionnel). 

Si $V$ est de type $(M,\lambda)$, elle est caractérisée par un invariant $\sL\in \BP^1(M_{\dR})$ qui est le paramètre de la filtration à deux crans sur $M_{\dR}\coloneqq (M\otimes_{\Qp^{\nr}}\Qpbar)^{\sG_{\Qp}}=\bD_{\dR}(V)$. De plus, $\lambda\in P$, $M\in \phinl$ et $\sL\in \BP^1(M_{\dR})$ définissent un $(\varphi,N,\sG_{\Qp})$-module filtré $M_{\sL}^{\lambda}$ qui, s'il est (faiblement) admissible\footnote{Si $M$ est cuspidal alors $M_{\sL}^{\lambda}$ est admissible pour tout $\sL\in \BP^1(M_{\dR})$, si $M$ est spécial il existe un seul $\sL\in \BP^1(M_{\dR})$ tel que $M_{\sL}^{\lambda}$ n'est pas admissible et si $M$ est exceptionnel tous les $M_{\sL}^{\lambda}$ admissibles sont isomorphes.}, définit $V_{M,\sL}^{\lambda}$, une $L$-représentation de $\sG_{\Qp}$, par \hbox{\cfc{cofo}} : 
\begin{equation}\label{eq:defvst}
V_{M,\CL}^{\lambda}\coloneqq \bV_{\st}(M_{\sL}^{\lambda}[1])\coloneqq(M\otimes_{\Qp^{\nr}}\Bst)^{\varphi=p}\cap \Fil^0(\Bdr\otimes_{\Qp}{M_{\dR}}),
\end{equation}
et $V=V_{M,\CL}^{\lambda}$, ce qui permet de reconstruire $V$ à partir de $(M,\lambda, \sL)$. De plus, $V_{M,\CL}^{\lambda}$ est de déterminant $\zeta_M^{\lambda}\omega$. 

La correspondance de Langlands $p$-adique (\cf \cite{colp}), $V\mapsto \bPi(V)$ associe une $L$-représentation de Banach unitaire admissible de $G$ à $V$ une $L$-représentation de $\sG_{\Qp}$ de dimension $2$. Si $\lambda\in P$ et $M\in \phinl$ ce foncteur établit une bijection entre les $L$-représentations absolument irréductibles $V$ de $\sG_{\Qp}$ de type $(M,\lambda)$ et les complétés admissibles unitaires irréductibles de $\LL_M^{\lambda}$ (\cf \cite{codopa}) qui sont définis par $\bPi_{M,\CL}^{\lambda}\coloneqq \bPi(V_{M,\CL}^{\lambda})$. Si $M$ est cuspidal (\resp spécial, exceptionnel), on dit que $\Pi_{M,\CL}^{\lambda}$ est \emph{cuspidale} (\resp \emph{spéciale}, \emph{exceptionelle}).

Réciproquement, Colmez construit (\cf \cite{colp}) un foncteur $\Pi\mapsto \bV(\Pi)$ qui associe à une $L$-représentation unitaire admissible de $G$ une $L$-représentation de $\sG_{\Qp}$ de dimension $2$. Ces foncteurs satisfont $\bV\circ \bPi=\id$, ce qui signifie que $\bPi(V)$ caractérise entièrement $V$.
\subsubsection{Blocs et anneaux de déformations}
Soit $\zeta\colon \Qpt\rightarrow \rmO_L^{\times}$ un caractère lisse unitaire. Notons $\Tors\, G$ (\resp $\Tors^{\zeta} G$) la catégorie abélienne des $\rmO_L$-représentations lisses de longueur finie (et de caractère central $\zeta$). On définit sur $\Irr(\Tors^{(\zeta)} G)$, l'ensemble des $k_L$-représentations irréductibles de $\Tors^{(\zeta)} G$, la relation d'équivalence $\sim$ engendrée par 
$$
\pi,\pi'\in\Irr(\Tors^{(\zeta)} G) \quad \Ext^1(\pi,\pi')=0 \implies\pi\sim\pi',
$$
et on appelle \emph{bloc (absolu)} de $\Tors\, G$ (\resp $\Tors^{\zeta} G$) une classe d'équivalence de $\Irr(\Tors^{(\zeta)} G)$ pour cette relation. La théorie de Gabriel \cfc{gab} fournit une décomposition suivant les blocs
$$
\Tors^{(\zeta)} G=\prod_{\CB} \Tors_{\CB}^{(\zeta)} G,
$$
où $\Tors_{\CB}^{(\zeta)} G\subset \Tors^{(\zeta)} G$ est la sous-catégorie pleine des $\rmO_L$-représentations de $\Tors^{(\zeta)} G$ dont toutes les composantes de Jordan-Hölder sont dans $\CB$. Dans \cite{pas}, Paškūnas décrit les blocs $\CB$ de cette catégorie (\cf \num \ref{subsubsec:azu}). La correspondance de Langlands $p$-adique semi-simple modulo $p$ s'énonce alors comme une bijection (\cf \cite{pas2}) :
\begin{equation}
\begin{array}{ccc}
\left\{
\begin{array}{c}
\text{Blocs absolues de}\\[2pt]
\text{$\Tors\, G$}
\end{array}
\right\}
&
\longleftrightarrow
&
\left\{
\begin{array}{c}
\text{$k_L$-représentations }\\[2pt]
\text{continues semi-simples de $\sG_{\Qp}$ }\\[2pt]
\text{de dimension $2$}
\end{array}
\right\}
\\[12pt]
\CB & \mapsto & \bar{\rho}_{\CB}
\end{array}
\end{equation}
On fixe un bloc $\CB$ et soit $\bar \rho_{\CB}$ la $k_L$-représentation de $\sG_{\Qp}$ associée. On note $R_{\CB}^{\ps,\zeta}$ l'anneau de déformation du pseudo-caractère $(\tr\bar\rho_{\CB},\det \bar \rho_{\CB})$ et de déterminant $\zeta\omega$. Pour $\lambda\in P$ et $M\in \phinl$, Kisin (\cf \cite{kis}) construit un quotient $R_{\CB}^{\ps,\zeta}\surj R_{\CB,M}^{\lambda,+}$ caractérisé par la propriété suivante : pour tout idéal maximal $\fkm_x\subset R_{\CB}^{\ps,\zeta}[1/p]$, dont on note $L_x$ le corps résiduel, le morphisme $R_{\CB}^{\ps,\zeta}\rightarrow L_x$ se factorise par $R_{\CB,M}^{\lambda,+}\rightarrow L_x$ si et seulement si la $L_x$-représentation définie par $R_{\CB}^{\ps,\zeta}\rightarrow L_x$ est, soit de type $(M,\lambda)$, ou bien de type $(M',\lambda)$ où $M'$ est le $L$-$(\varphi,\sG_{\Qp})$-module $M$ muni de la monodromie $N=0$ (\cf la remarque \ref{rem:degenexp}). On note $\rho_{\CB,M}^{\lambda,+}\colon \sG_{\Qp}\rightarrow \GL_2(R_{\CB,M}^{\lambda,+})$ la $R_{\CB,M}^{\lambda,+}$-représentation universelle de $\sG_{\Qp}$.

On pose 
$$
R_{\CB,M}^{\lambda}\coloneqq R_{\CB,M}^{\lambda,+}\bigl [\tfrac 1 p\bigr],\quad \rho_{\CB,M}^{\lambda}\coloneqq \rho_{\CB,M}^{\lambda,+} [\tfrac 1 p\bigr],\quad X_{\CB,M}^{\lambda}\coloneqq \Spec\, R_{\CB,M}^{\lambda}.
$$

Soit $\lambda\in P$,  $M\in \phinl$ et $\CB$ un bloc de $\Tors^{\zeta_M^{\lambda}}G$ où l'on rappelle que $\zeta_M^{\lambda}$ est le caractère central de $\LL_M^{\lambda}$ (et de $\JL_M^{\lambda}$).
\begin{itemize}
\itemb La correspondance de Langlands $p$-adique en famille (\cf la définition \ref{def:plangfam}) permet d'associer à $ \rho_{\CB,M}^{\lambda}$ une $R_{\CB,M}^{\lambda}$-représentation topologique $\bPi(\rho_{\CB,M}^{\lambda})$. On note $\bPi(\rho_{\CB,M}^{\lambda})'$ son $L$-dual stéréotypique\footnote{C'est le $L$-dual continu pour la topologie de la convergence uniforme sur les parties fortement bornées, \ie la dualité pour laquelle les espaces de Banach et de Fréchets sont réflexifs.}. Pour $\fkm_x\subset R_{\CB,M}^{\lambda}$ un idéal maximal dont on note $L_x$ le corps résiduel, on note $\rho_x$ la $L_x$-représentation de $\sG_{\Qp}$ donnée par $\fkm_x$. La représentation $\bPi(\rho_{\CB,M}^{\lambda})$ interpole les $\bPi(\rho_x)$ au sens où $ \bPi(\rho_{\CB,M}^{\lambda})'\otimes_{R_{\CB,M}^{\lambda}}L_x\cong \bPi(\rho_x)'$, c'est-à-dire que les $\bPi(\rho_x)'$ sont des quotients de $\bPi(\rho_{\CB,M}^{\lambda})'$.
\itemb On définit $\wh{\LL}_M^{\lambda}$ le \emph{complété unitaire universel} de $\LL_M^{\lambda}$ (\cf \cite{codo}) et on choisit un $\rmO_L$-réseau $\wh{\LL}_M^{\lambda,+}\subset \wh{\LL}_M^{\lambda}$. Le \emph{complété $\CB$-adique} (\cf \cite{codonikis} et \num \ref{subsubsec:defibcomp}) $\LL_{M,\CB}^{\lambda,+}$ de $\LL_M^{\lambda,+}$ est défini comme la limite projective de ses quotients de longueur finie contenus dans $\Tors_{\CB}^{\zeta_{M}^{\lambda}} G$. En particulier,
$$
\wh{\LL}_M^{\lambda,+}=\prod_{\CB}\LL_{M,\CB}^{\lambda,+},
$$
et on pose $\LL_{M,\CB}^{\lambda}\coloneqq \LL_{M,\CB}^{\lambda,+}\otimes_{\rmO_L}L$. 

\end{itemize}

\subsubsection{Cohomologie du système local $p$-adique}
Le \emph{demi-plan $p$-adique de Drinfeld} est défini comme
$$
\BH_{\Qp}\coloneqq\BP^{1,\an}_{\Qp}\setminus \BP^1(\Qp),
$$
en tant que $\Qp$-espace rigide analytique. Le groupe $G$ agit par homographie sur $\BH_{\Qp}$ et $\czG$ agit trivialement. Posons
$$
\BG\coloneqq G\times \sG_{\Qp}\times \czG.
$$
Soit $K$ une extension finie de $\Qp$ contenant une extension quadratique de $\Qp$, on définit $\presp{\rmM}_K^0\coloneqq \BH_K\times \BZ/2$ qui est muni d'une action de $\BG$ où l'action sur $\BZ/2$ est donnée par $+v_p(\det(g))$ pour $g\in G$, par le quotient non ramifié $\sG_{\Qp}\rightarrow \BZ/2$ pour $\sG_{\Qp}$ et par $+v_p(\nrd(\check{g}))$ pour $\check g \in \czG$. Par convention, pour toute extension finie $K$ de $\Qp$, on définit $\presp{\rmM}_{K}^0$ par descente galoisienne à partir d'une extension quadratique de $K$. Le théorème de Drinfeld \cfc{drin} fournit une application surjective
\begin{equation}\label{eq:drinsurj}
\pi_1^{\et}(\BH_{\Qp})\surj \czG^+\subset \czG,
\end{equation}
où $\czG^+\coloneqq \rmODt$. On peut alors associer à $\presp{\rmM}_K^0$ des revêtements étales et des systèmes locaux.
\begin{itemize}
\itemb Pour $n\geqslant 1$, le sous-groupe $\czG_n^+\coloneqq 1+\varpi_{\rmD}^n\rmOD\subset \czG^+$ définit un revêtement étale galoisien $\presp{\rmM}_K^n\rightarrow \presp{\rmM}_K^0$. Pour $K=\Qp$ c'est un revêtement de groupe de Galois\footnote{De plus si $K$ est assez grand c'est un revêtement de groupe de Galois $\czG^1/\czG_n^1$ où $\czG^1\coloneqq\ker \nrd$ est le noyau de la norme réduite et $\czG_n^1\coloneqq\czG^1\cap \czG^1_n$.} $\czG^+/\czG_n^+$. On obtient ainsi une action de $\BG$ sur $\presp{\rmM}_K^n$.
\itemb Pour $\lambda\in P_+$, la $L$-représentation $\czW_{\lambda}$ (\cf (\ref{eq:defczw})) de $\czG$ contient un réseau stable $\czW_{\lambda}^+$ et la surjection (\ref{eq:drinsurj}) permet de définir une $\rmO_L$-représentation de $\pi_1^{\et}(\presp{\rmM}_K^n)$, c'est-à-dire \cfc{dej}, un $\rmO_L$-système local étale $\BV_{\lambda}^+$ et son $L$-système local associé $\BV_{\lambda}\coloneqq \BV_{\lambda}^+\otimes_{\rmO_L}L$ sur $\presp{\rmM}_K^n$.
\end{itemize}
On veut calculer la cohomologie étale du système local 
$$
\Sym \BV(1)\coloneqq \bigoplus_{\lambda\in P_+}\BV_{\lambda}(1)\otimes_{L}\czW_{\lambda},
$$
où $(1)$ désigne la torsion par le caractère cyclotomique. Dans \cite{van}, on donne une interprétation de ce système local comme l'algèbre symétrique du module de Tate du $\rmOD$-module formel spécial universel sur $\presp{\rmM}_K^n$. On n'utilisera pas cette interprétation ici. On pose 
$$
\begin{gathered}
\rmH^1_{\et}\Harg{\presp{\rmM}_{\Qpbar}^n}{\Sym \BV(1)}\coloneqq \bigoplus_{\lambda\in P_+}\rmH^1_{\et}\Harg{\presp{\rmM}_{\Qpbar}^n}{\BV_{\lambda}(1)}\otimes_{L}\czW_{\lambda},\\
\rmH^1_{\et}\Harg{\presp{\rmM}_{\Qpbar}^n}{\BV_{\lambda}(1)}\coloneqq \left ( \varprojlim_n\varinjlim_{[K\colon \Qp]<\infty}\rmH^1_{\et}\Harg{\presp{\rmM}_{K}^n}{\BV^+_{\lambda}(1)/p^n}\right ) \otimes_{\rmO_L}L.
\end{gathered}
$$
Rappelons que dans \cite{van}, nous avons montré comment les $L$-représentations $V_{M,\CL}^{\lambda}$ et $\Pi_{M,\CL}^{\lambda}$ apparaissent par entrelacement dans la cohomologie de $\presp{\rmM}_{C}^n$ à coefficients dans $\Sym \BV(1)$. Il s'agit ici d'utiliser ce résultat pour obtenir une description complète de la cohomologie de $\presp{\rmM}_{\Qpbar}^n$ à coefficients dans $\Sym \BV(1)$.

\subsection{Plan de l'article}
Après cette introduction, l'article est constitué de deux sections. La seconde section est consacrée à la preuve du théorème \ref{thm:princint2} qui est ensuite utilisée pour démontrer le théorème \ref{thm:princint}. 
\begin{itemize}
\itemb Dans la section $2$, après avoir énoncé le résultat principal, on commence, (\cf \num \ref{subsec:seriesperap}), par quelques rappels sur les représentations spéciales et on construit un réseau dans le complété unitaire universel de la steinberg localement algébrique (\cf \num \ref{subsubsec:latt}). On calcule alors les déformations de de Rham infinitésimales universelles des représentations spéciales (\cf \num\ref{subsec:definfgal} et \num\ref{subsec:definfaut}) avant d'en déduire les complétés de longueur finie de la steinberg localement algébrique (\cf la proposition \ref{prop:complst}). On utilise ensuite ces résultats pour déterminer les anneaux locaux complétés de l'anneau de Kisin (\cf \num \ref{subsubsec:lafamun}) après avoir rappellé quelques propriétés des $\CB$-complétés (\cf \num\ref{subsec:Bcomkis}). Finalement, ceci permet de définir une application analytique de l'analytifié du spectre de l'anneau de Kisin vers la droite projective et de montrer que cet application est une immersion ouverte (\cf  \num \ref{subsubsec:lopen}), ce qui achève la preuve du théorème (\cf \num\ref{subsec:finsecd}).
\itemb Dans la section $3$, on commence par définir les derniers termes de l'énoncé du théorème \ref{thm:princint}. On complète ensuite certains résultats de \cite{van}, notamment les théorèmes de finitude (\cf \num\ref{subsubsec:finit}) et le passage de $\brv{\rmM}_C^n$ à $\presp{\rmM}^n_{\Qpbar}$ (\cf \num\ref{subsubsec:vpl}). On démontre ensuite la décomposition en blocs au niveau entier (\cf \num\ref{subsubsec:decentier}). Enfin, après quelques lemmes techniques (\cf \num\ref{subsubsec:azu} et le lemme \ref{lemm:priptto}), on termine la preuve du théorème \ref{thm:princint} au \num \ref{subsubsec:final}.
\end{itemize}
\subsection{Remerciements}
Je tiens à remercier Colmez, Dospinescu et Nizio\l ; il est évident que ce travail leur doit énormément. 
Je remercie le Morningside Cente of Mathematics de la Chinese Academy of Sciences, où ce travail a été réalisé, pour ses excellentes conditions de travail. Je remercie particulièrement Zicheng Qian pour de nombreuses discussions éclairantes.

\section{Anneaux de Kisin}
Soit $\lambda \in P$ et $M\in \phinl$. Si $M$ est spécial supposons de plus que $w(\lambda)>1$.
\begin{equation}\label{eq:defr}
T_{M,\CB}^{\lambda}\coloneqq \End(\LL_{M,\CB}^{\lambda}),\quad Y_{M,\CB}^{\lambda}\coloneqq\Spec(T_{M,\CB}^{\lambda}),\quad \sigma_{M,\CB}^{\lambda}\coloneqq\bV(\LL_{M,\CB}^{\lambda}).
\end{equation}
On veut démontrer que $T_{M,\CB}^{\lambda}$ est un anneau de déformation et que $\sigma_{M,\CB}^{\lambda}$ est la famille universelle sur ce dernier. Le théorème s'énonce de la façon suivante :

\medskip
\begin{theo}\label{conj:kisin}
\
\begin{enumerate}
\item L'anneau $T_{M,\CB}^{\lambda}$ s'identifie aux fonctions bornées sur un ouvert analytique de $\BP^{1,\an}_L$ : algébriquement, c'est le produit d'anneaux principaux.
\item On a un isomorphisme d'anneaux $R_{\CB,M}^{\lambda}\cong T_{M,\CB}^{\lambda}$. De plus $\sigma_{\CB}^{\lambda}$ est canoniquement isomorphe à la représentation universelle sur $R_{\CB,M}^{\lambda}$, \ie $\sigma_{M,\CB}^{\lambda}\cong \rho_{\CB,M}^{\lambda}$.
\end{enumerate}
\end{theo}
\medskip
\medskip
\begin{rema}\label{rem:kisthm}
\ 
\begin{itemize}
\itemb Quitte prendre $L$ plus grand, on peut supposer que tous blocs sont absolus (\cf \cite[Remarque 1.9]{codonikis}), ce que l'on fera dans toute la suite.
\itemb Pour $M$ cuspidal, ce théorème est démontré dans \cite{codonikis}. Dans la section $2$, on traite le cas où $M$ est spécial par les mêmes méthodes.
\itemb La difficulté principale provient de la \emph{représentation exceptionnelle} qui correspond à l'idéal maximal de $R_{\CB,M}^{\lambda}$, pour $M$ spécial et $\CB$ un bloc précis, pour lequel la représentation associée est cristalline. La raison est que cette représentation admet une déformation cristalline et on doit montrer qu'elle n'apparaît pas dans $R_{\CB,M}^{\lambda}$. 
\itemb Pour $M$ spécial et $w(\lambda)$ satisfaisant $1\leqslant w(\lambda) \leqslant p-1$, en utilisant les résultats de Chitrao et Ghate (\cf \cite{chgh}) on peut déterminer explicitement l'ouvert $(X_{\CB,M}^{\lambda})^{\an}$.
\itemb Dans le cas où $w(\lambda)=1$, la convention sera que $R_{\CB,M}^{\lambda}= L$ pour le bloc de la steinberg, auquel cas $\rho_{M,\CB}^{\lambda}=\zeta_M^{\lambda}\omega$, et $R_{\CB,M}^{\lambda}=0$ sinon.
\end{itemize}
\end{rema}
\medskip
Le but de cette section est donc de démontrer la proposition dans le cas où $M$ est spécial, \ie $M=\Sp(\lvert \lambda \rvert-2 )$ (\cf \num \ref{subsubsec:defrepspe}). Dans ce cas, $\LL_M^{\lambda}=\St_{\lambda}^{\lalg}$.

Soit $\zeta\colon \Qpt\rightarrow L^{\times}$ un caractère unitaire, que l'on voit comme un caractère de $G$ en précomposant avec le déterminant, et comme un caractère $\sG_{\Qp}$ en précomposant avec morphisme de réciprocité. Alors on définit
\begin{itemize}
\itemb $R^{\ps,\zeta}_{\CB}$  l'anneau de pseudo-déformation universel de $(\tr \bar\rho_{\CB}, \det \bar \rho_{\CB})$ de déterminant $\zeta\omega$,
\itemb $Z_{\CB}^{\zeta}$ le centre de la catégorie abélienne $\Tors_{\CB}^{\zeta}G$.
\end{itemize}
L'un des outils principaux est le théorème de type $R=T$ de Paškūnas et Tung\footnote{Comme on a fait l'hypothèse $p>3$, c'est en fait \cite{pas} que l'on utilise. Comme suggéré dans \cite{codonifac}, il est possible que cette hypothèse soit superflue, même s'il n'est pas tout à fait clair pour nous comment adapter certains arguments.}, (\cf \cite{patu}) :

\medskip
\begin{theo}\label{thm:patu}
Soit $\zeta\colon\Qpt\rightarrow \rmO_L^{\times}$ un caractère lisse. Il existe un unique isomorphisme
$$
\iota_{\CB}^{\zeta}\colon R^{\ps,\zeta}_{\CB}\xrightarrow{\sim}Z_{\CB}^{\zeta}.
$$
\end{theo}
\medskip

\medskip
\begin{rema}
Il n'est pas clair comment déduire, dans le cas spécial, la conjecture de Breuil-Mézard  des résultats de cette section, comme dans \cite{codonikis}. L'une des raisons est que la réduction modulo $\varpi_L$ du réseau que l'on construit au \num \ref{subsubsec:latt} n'est \emph{a priori} pas semi-simple\footnote{Pour un $K$-type de la forme $\tau\otimes W$ avec $\tau$ lisse et $W$ algébrique, sa réduction modulo $\varpi_L$ n'est \emph{a priori} pas semi-simple.}. Dans l'énoncé de \cite[Proposition 5.9]{codonikis} et sa preuve, il est implicitement utilisé que l'on peut choisir un réseau dans un $K$-type cuspidal dont la réduction modulo $\varpi_L$ soit semi-simple, ce qui permet la décomposition en somme directe. 
\end{rema}
\medskip

\subsection{La série spéciale $p$-adique}\label{subsec:seriesperap}
\subsubsection{Définition des représentations spéciales}\label{subsubsec:defrepspe}
Soit $k\in \BZ$ un entier, on définit le $L$-$(\varphi,N)$-module spécial $\Sp_L(k)$ comme le $L$-espace vectoriel $\Sp_L(k)=Le_0\oplus L e_1$ muni des endomorphismes
$$
\begin{cases}
\varphi(e_1)=p^{-\frac{k-1}{2}}e_1\\
\varphi(e_0)=p^{-\frac{k+1}{2}}e_0,
\end{cases}
\quad
\begin{cases}
Ne_1=e_0\\
Ne_0=0,
\end{cases}
$$
qui est de pente $-\tfrac k 2$. Tout $L$-$(\varphi,N,\sG_{\Qp})$-module spécial non exceptionnel est de la forme $(\Sp_L(k)\otimes_{\Qp}\Qp^{\nr})\otimes_L \chi$ pour $\chi \colon \Qpt \rightarrow L^{\times}$ un caractère lisse (\cf la preuve de \cite[Proposition 1.3]{sch}). 

Soit $k\in \BZ$ un entier, on définit le $L$-$(\varphi,N)$-module exceptionnel $\wt{\Sp}_L(k)$ comme le $L$-espace vectoriel $\wt{\Sp}_L(k)=Le_0\oplus L e_1$ muni des endomorphismes
$$
\begin{cases}
\varphi(e_1)=p^{-\frac{k-1}{2}}e_1\\
\varphi(e_0)=p^{-\frac{k+1}{2}}e_0,
\end{cases}
\quad
\begin{cases}
Ne_1=0\\
Ne_0=0,
\end{cases}
$$
C'est-à-dire que $\wt{\Sp}_L(k)$ est le $\varphi$-module $\Sp_L(k)$ avec $N=0$ ; sa pente est également $-\tfrac k 2$. Tout $L$-$(\varphi,N,\sG_{\Qp})$-module exceptionnel est de la forme $(\wt{\Sp}_L(k)\otimes_{\Qp}\Qp^{\nr})\otimes_L \chi$ pour un caractère lisse $\chi \colon \Qpt \rightarrow L^{\times}$.

Soit $\lambda\in P_+$ et $\CL\in L\cup\{\infty\}$.

\begin{itemize}
\itemb Si $\CL\in L$,  on définit $D_{\CL}^{\lambda}$ le $L$-$(\varphi,N)$-module filtré dont le $L$-$(\varphi,N)$-module sous-jacent $\Sp_L(\lvert \lambda \rvert-2)$ et muni de la filtration 
$$
\Fil^iD_{\CL}^{\lambda}\coloneqq 
\begin{cases}
D_{\CL}^{\lambda}&\text{ si } i\leqslant -\lambda_2+1,\\
L(e_1-\CL e_0) &\text{ si } -\lambda_2+2\leqslant i\leqslant -\lambda_1+1,\\
0&\text{ si }  -\lambda_1+2\leqslant i.
\end{cases}
$$
\itemb Si $\CL =\infty$,  on définit $D_{\CL}^{\lambda}=D_{\infty}^{\lambda}$ comme le $L$-$(\varphi,N)$-module filtré de $L$-$(\varphi,N)$-module sous-jacent $\wt{\Sp}_L(\lvert \lambda \rvert )$ et de filtration 
$$
\Fil^iD_{\CL}^{\lambda}\coloneqq 
\begin{cases}
D_{\CL}^{\lambda}&\text{ si } i\leqslant -\lambda_2+1,\\
L(e_1- e_0) &\text{ si } -\lambda_2+2\leqslant i\leqslant -\lambda_1+1,\\
0&\text{ si }  -\lambda_1+2\leqslant i.
\end{cases}
$$
\end{itemize}

On pose $V_{\CL}^{\lambda}\coloneqq \bV_{\st}(D_{\CL}^{\lambda}[1])$ (\cf (\ref{eq:defvst})), qui est une $L$-représentation spéciale de $\sG_{\Qp}$ et $\Pi_{\CL}^{\lambda}\coloneqq \bPi(V_{\CL}^{\lambda})$ qui est une $L$-représentation spéciale de $G$. Puisque $w(\lambda)>1$, ces représentations sont absolument irréductibles ; $\Pi_{\CL}^{\lambda}$ est de caractère central\footnote{Pour $M=\Sp_L(\lvert \lambda \rvert-2)$ (ou $M=\wt{\Sp}_L(\lvert \lambda \rvert-2)$), $\zeta_{\lambda}\coloneqq \zeta^{\lambda}_M$ défini dans l'introduction.} $\zeta_{\lambda}\coloneqq\omega^{\lvert \lambda \rvert-1}$ et $V_{\CL}^{\lambda}$ est de déterminant $\zeta_{\lambda}\omega=\omega^{\lvert \lambda \rvert}$.

\subsubsection{Complété unitaire universel de $\St_{\lambda}^{\lalg}$}
Soit $\sC^{\lambda}\coloneqq \wh{\St_{\lambda}^{\lalg}}$ le complété unitaire universel de $\St_{\lambda}^{\lalg}$. On connaît cet espace explicitement, en termes de fonctions de classe $\sC^{u(\lambda)}$, où l'on définit $u(\lambda)\coloneqq \tfrac{w(\lambda)-1}{2}$.

En effet, on définit $\sC^{\lambda}$ comme le quotient $\sC^{\lambda}\coloneqq \wt{\sC}^{\lambda}/W_{\lambda}^*$ où $\wt{\sC^{\lambda}}\subset \sC^0(\Qp,L)$ est le sous-espace des fonctions  $f\colon \Qp\rightarrow L$  de classe $\sC^{u(\lambda)}$ (\cf \cite[I.5]{colfun} pour la définition) et la fonction $x\mapsto x^{w(\lambda)-1}f(1/x)$ s'étend en une fonction de classe $\sC^{u(\lambda)}$ sur $\Qp$ et $W_{\lambda}^*\subset \wt{\sC}^{\lambda}$ est le sous-espace des fonctions polynomiales de degré au plus $w(\lambda)-1$. L'action de $G$ sur $\wt{\sC}^{\lambda}$ est définie par
$$
g=
\begin{pmatrix}
a & b \\ c & d
\end{pmatrix}
\in G, \quad
(g\star_{\lambda} f)(x) =\lvert \det(g)\rvert^{\tfrac{\lvert \lambda \rvert-1} 2}(cx+d)^{w(\lambda)-1} f\left ( \frac{ax+b}{cx+d} \right )
$$

Pour $\CL\in \BP^1(L)$, $a\in \Qp$ et $n\geqslant 1$ un entier, on définit la fonction $f_{\CL,a}^n\in \sC^0(\Qp,L)$ par
$$
x\mapsto f_{\CL,a}^n(x)\coloneqq (x-a)^{n}\log_{\CL}(x-a)
$$ 
où si $\CL=\infty$ on pose $\log_{\infty}=v_p$.

De plus, pour $\CL\in \BP^1(L)=L\cup \{\infty\}$, on définit $\wt{M}_{\CL}^{\lambda}\subset \sC^0(\Qp,L)$ le sous-espace engendré par les fonctions $h\colon \Qp \rightarrow L$ de la forme
$$
h=\sum_{i\in I}\lambda_if^{n_i}_{\CL,a_i},
$$
où $I$ est un ensemble fini, et pour tout $i\in I$, $n_i\geqslant 1$ est un entier, $\lambda_i\in L$ et $a_i\in \Qp$, tels que 
\begin{itemize}
\itemb pour tout $i\in I$, $u(\lambda)< n_i\leqslant w(\lambda)$,
\itemb $\deg\left ( \sum_{i\in I}\lambda_i(x-a_i)^{n_i}\right)< u(\lambda)$.
\end{itemize}
En d'autres termes, $\wt{M}_{\CL}^{\lambda}\subset \sC^0(\Qp,L)$ est l'espace des fonctions qui sont des combinaisons linéaires de fonctions de la forme $f^{n}_{\CL,a}$ avec $u(\lambda)<n\leqslant w(\lambda)$ et qui ont un pôle d'ordre $<u(\lambda)$ en $\infty$ (\cf \cite[V.1., 3]{colprin}). 
\medskip
\begin{lemm}
On a $\wt{M}_{\CL}^{\lambda}\subset \wt{\sC}_{\lambda}$. On note $M_{\CL}^{\lambda}$ l'image de $\wt{M}_{\CL}^{\lambda}$ dans $\sC^{\lambda}$ et $\wh{M}_{\CL}^{\lambda}$ l'adhérence de $M_{\CL}^{\lambda}$ dans $\sC^{\lambda}$. On a une suite exacte de $L$-représentations de Banach unitaires de $G$
$$
0\rightarrow \wh{M}_{\CL}^{\lambda}\rightarrow \sC^{\lambda}\rightarrow \Pi_{\CL}^{\lambda}\rightarrow 0.
$$
\end{lemm}
\medskip
\begin{proof}
Ces résultats sont dus à Breuil, \cf \cite[Lemme 3.3.2]{breprin} pour le fait que $\wt{M}_{\CL}^{\lambda}\subset \wh{\sC}_{\lambda}$ et \cite[Corollaire 3.3.4]{breprin} pour la suite exacte.
\end{proof}
On finit ce numéro par trois énoncés pratiques sur $\sC^{\lambda}$. On note $\sD^{\lambda}$ le dual (stéréotypique) de $\sC^{\lambda}$. La preuve du lemme suivant est exactement celle de \cite[Lemme 5.3 b)]{codonifac} :
\medskip
\begin{lemm}\label{lem:gbded}
Les vecteurs $G$-bornés de $(\St_{\lambda}^{\lalg})'$ et de $(\St_{\lambda}^{\lan})'$ sont isomorphes à $\sD^{\lambda}$.
\end{lemm}
\medskip
On calcule les vecteurs localement algébriques de $\sC^{\lambda}$ : 
\medskip
\begin{prop}\label{prop:localgvecuni}
$\End_{L[G]}(\sC^{\lambda})=L$
\end{prop}
\begin{proof}
Soit $\pi\coloneqq \St_{\lambda}^{\lalg}$, alors $\pi$ est dense dans $\wh\pi$ et comme $\End_{L[G]}(\pi)=L$ par le lemme de Schur classique, il suffit de prouver que $\wh{\pi}^{\lalg}=\pi$ pour conclure. Mais puisque $\wh{\pi}\cong \sC^{\lambda}$, et les vecteurs localement algébriques de $\sC^{u(\lambda)}(\BP^1(\Qp),L)$ sont $\sC^{\lalg}(\BP^1(\Qp),L)\cong \sC^{\infty}(\BP^1(\Qp),L)\otimes W_{\lambda}^*$, donc on a $(\sC^{\lambda})^{\lalg}\cong\St_{\lambda}^{\lalg}$.
\end{proof}
On caractérise les quotients de $\sC^{\lambda}$ de longueur finie comme complétés de $\sC^{\lambda}$ de longueurs finis.
\medskip
\begin{prop}\label{prop:denstoquo}
Soit $\Pi$ une $L$-représentation unitaire de $G$ de longueur finie. Alors $\Pi$ est un quotient de $\sC^{\lambda}$ si et seulement si on a un morphisme $\St_{\lambda}^{\lalg}\incl \Pi$ d'image dense.
\end{prop}
\medskip
\begin{proof}
Si $\Pi$ est un quotient de $\sC^{\lambda}$, puisque par la propriété universelle on a 
$$
\Hom_G(\sC^{\lambda},\Pi)\cong\Hom_G(\St_{\lambda}^{\lalg},\Pi)
$$
alors on a une application $\St_{\lambda}^{\lalg}\rightarrow \Pi$ nécessairement injective puisque $\St_{\lambda}^{\lalg}$ est irréductible et d'image dense par définition du complété unitaire universel.

Réciproquement, on obtient une application $\sC^{\lambda}\rightarrow \Pi$ dont on doit montrer la surjectivité. Comme pour la preuve de \cite[Proposition 3.2]{codonikis}, c'est une conséquence du fait que $\Pi$ est résiduellement de longueur finie \cfc{codopa}.
\end{proof}

\Subsubsection{Réseau de $\sC^{\lambda}$}\label{subsubsec:latt}
Dans ce numéro, on construit un réseau de $\St_{\lambda}^{\lalg}$ dont le complété $p$-adique est un réseau dans $\sC^{\lambda}$. Posons $K\coloneqq \GL_2(\Zp)\subset G$, $I\subset K$ le \emph{sous-groupe d'Iwahori} (\ie le sous-groupe des matrices triangulaires supérieures modulo $p$) et $Z\subset G$ le centre de $G$. De plus, on note $N\coloneqq N_G(IZ)$ le normalisateur de $IZ$ et $KN$ sont produit avec $K$. On a
$$
w_p\coloneqq \begin{pmatrix} 0 & p \\ 1 & 0\end{pmatrix},\quad N=IZ\sqcup w_p IZ,\quad KN=KZ\sqcup w_p KZ
$$
Pour $\lambda\in P$, on a une $\rmO_L$-représentation localement algébrique de $G$ de dimension finie
$$
W_{\lambda}^{+,*}\coloneqq\Sym_{\rmO_L}^{w(\lambda)-1}\otimes_{\rmO_L} {\det}^{\lambda_1}\otimes_{\rmO_L}\lvert \det\rvert^{\tfrac{\lvert \lambda \rvert -1}{2}}.
$$
C'est un $\rmO_L$-réseau de $W_{\lambda}^*$ stable par $KN$. Posons
$$
\CI_0^{\lambda,(+)}\coloneqq \ind_{KZ}^GW_{\lambda}^{(+),*} ,\quad \wt{\CI}_1^{\lambda,(+)}\coloneqq \ind_{IZ}^GW_{\lambda}^{(+),*}
$$
où les induites sont compactes. Alors $\wt{\CI}_1^{\lambda,+}$ est un $\rmO_L$-réseau $G$-stable de $\wt{\CI}_1^{\lambda}$ et  $\CI_0^{\lambda,+}$ est un $\rmO_L$-réseau $G$-stable de $\CI_0^{\lambda}$. L'action de $w_p\in G$ définit une involution $\Pi\colon \wt{\CI}_1^{\lambda,(+)}\rightarrow \wt{\CI}_1^{\lambda,(+)}$ par $\Pi(f)=f \circ w_p$. Ainsi, on définit $\CI_1^{\lambda,(+)}\coloneqq (\wt{\CI}_1^{\lambda,(+)})^{\Pi=-\id}$, de sorte que l'on a 
$$
\wt{\CI}_1^{\lambda,(+)}=\CI_1^{\lambda,(+)}\oplus (\wt{\CI}_1^{\lambda,(+)})^{\Pi=\id}.
$$
La projection $\wt{\CI}_1^{\lambda,(+)}\rightarrow\CI_1^{\lambda,(+)}$ est donnée par $\Pi_-\coloneqq \id-\Pi$.
\medskip
\begin{prop}\label{prop:latt}
On a une surjection $G$-équivariante $\CI_1^{\lambda}\surj \St_{\lambda}^{\lalg}$ telle que si l'on note $I^{\lambda,+}$ l'image de  $\CI_1^{\lambda,+}$ par cette surjection, qui définit un réseau de $\St_{\lambda}^{\lalg}$, alors le complété $p$-adique de $I^{\lambda,+}$ définit un réseau de $\sC^{\lambda}$, \ie 
$$
\sC^{\lambda}=(\varprojlim_nI^{\lambda,+}/p^n)\otimes_{\rmO_L}L.
$$
\end{prop}
\medskip
\medskip
\begin{rema}
Notons que la surjection $\CI_1^{\lambda}\surj \St_{\lambda}^{\lalg}$ est exactement celle obtenue par Chitrao et Ghate \cite[\S 4]{chgh}. Notre preuve est légèrement différente et repose (de manière caché à travers \cite{AdS}) sur l'isomorphisme $\rmH^1_{\dR}(\BH_{\Qp})\cong (\St_L^{\infty})'$ alors que le calcul de Chitrao et Ghate est explicite.
\end{rema}
\medskip
Commençons par expliciter le réseau $I^{\lambda,+}$. Soit $\st_{(0,1)}^{+}\coloneqq \sC(\BP^1(\BF_p),\rmO_L)/\rmO_L$ la \emph{$\rmO_L$-représentation de Steinberg finie}, que l'on considère comme une $\rmO_L$-représentation de $KZ$ par inflation. Pour $\lambda \in P$, posons
$$
{\st}_{\lambda}^+\coloneqq {\st}_{(0,1)}^{+}\otimes_{\rmO_L}W_{\lambda}^{*,+}.
$$
Comme $KZ/IZ\cong \BP^1(\BF_p)$, on obtient  $\ind_{IZ}^{KZ}W_{\lambda}^{*,+}\cong  \sC(\BP^1(\BF_p),W_{\lambda}^{*,+})$ et donc en appliquant le foncteur exact $\ind_{KZ}^G$ à la suite exacte 
$$
0\rightarrow W_{\lambda}^{*,+}\xrightarrow{\iota} \sC(\BP^1(\BF_p),W_{\lambda}^{*,+})\rightarrow {\st}_{\lambda}^+\rightarrow 0,
$$
on obtient une suite exacte 
$$
0\rightarrow \CI_0^{\lambda,+}\xrightarrow{\ind_{KN}^G\iota}\wt{\CI}_1^{\lambda,+}\rightarrow \ind_{KZ}^G{\st}_{\lambda}^+\rightarrow 0,
$$
et l'on obtient un diagramme commutatif définissant $I^{\lambda,+}$ : 
\begin{center}
\begin{tikzcd}
0\ar[r]&\CI_0^{\lambda,+}\ar[r,"\ind_{KZ}^G\iota"]\ar[d,equal]&\wt{\CI}_1^{\lambda,+}\ar[r]\ar[d,"\Pi_{-}"]& \ind_{KZ}^G{\st}_{\lambda}^+\ar[r]\ar[d]&0\\
0\ar[r]& \CI_0^{\lambda,+}\ar[r,"\partial^+"]& \CI_1^{\lambda,+}\ar[r]&I^{\lambda,+} \ar[r]& 0.
\end{tikzcd}
\end{center}
Plus exactement on définit $\partial^+$ comme la composée $\CI_0^{\lambda,+}\rightarrow \wt{\CI}_0^{\lambda,+}\xrightarrow{\ind_{KZ}^G\iota}\wt{\CI}_1^{\lambda,+}\xrightarrow {\Pi_-}\CI_1^{\lambda,+}$ et $I^{\lambda,+}\coloneqq \CI_1^{\lambda,+}/\partial^+(\CI_0^{\lambda,+})$. On définit finalement $\partial=\partial^+\otimes_{\rmO_L}L\colon \CI_0^{\lambda}\rightarrow \wt{\CI}_0^{\lambda}$. Notons que $(\wt{\CI}_1^{\lambda,(+)})^{\Pi=\id}=\ind_{N}^GW_{\lambda}^{+,*}$ et donc
\begin{equation}\label{eq:sttoi}
\ind_{KZ}^G{\st}_{\lambda}^+=\ind_{N}^GW_{\lambda}^{+,*}\oplus I^{\lambda,+}
\end{equation}
en particulier, on obtient le corollaire suivant :
\medskip
\begin{coro}\label{lem:modpfinlen}
Le $k_L[G]$-module lisse $I^{\lambda,+}/\varpi_L$ est de type fini.
\end{coro}
\medskip
On peut expliciter $\partial^+$ (\resp $\partial$) en interprétant les éléments des induites comme des fonctions. Notons $\wt{\sfs}\colon G/IZ\rightarrow G/KZ$ la projection canonique et $\wt{\sft}=\wt{\sfs}\circ w_p\colon G/IZ\rightarrow G/KZ$. L'application $\ind_{KN}^G\iota$ est explicitement donnée sur $f\in \CI_0^{\lambda,+}$, vu comme fonction sur $G/KZ$, pour $a\in G/IZ$ par
$$
(\ind_{KZ}^G\iota)(f)(a)=f(\wt\sfs(a)).
$$
Ainsi,  pour $a\in G/IZ$, 
\begin{equation}\label{eq:expart}
\partial^{(+)}(f)(a)=f(\wt\sfs(a))-f(\wt\sft(a)).
\end{equation}

La preuve de la proposition \ref{prop:latt} se fait par dualité. On note $\Ind$ l'induite sans condition sur le support. On commence par le lemme classique suivant, dont on esquisse une démonstration pour le confort du lecteur.
\medskip
\begin{lemm}\label{lem:indual}
Soit $H\subset G$ un sous-groupe ouvert et compact modulo le centre et soit $\sigma^+$ une $\rmO_L$-représentation localement algébrique de $H$, libre et de rang finie sur $\rmO_L$. Notons $\sigma=\sigma^+\otimes_{\rmO_L}L$, alors
\begin{enumerate}
\item $(\Ind_H^G\sigma)$ est naturellement un espace de Fréchet et $(\Ind_H^G\sigma)' \cong \ind_H^G\sigma^{*}$,
\item $(\Ind_H^G\sigma^+)' \cong\wh{\ind_H^G\sigma^{*,+}}$, où le complété à droite est $p$-adique,
\item $(\Ind_H^G\sigma)^{G-\rmb}\cong \Ind_H^G\sigma^{+}\otimes_{\rmO_L}L$.
\end{enumerate}
\end{lemm}
\medskip
\begin{proof}
Posons $V=(\Ind_H^G\sigma)$ et $W=\ind_H^G\sigma^{*}$. Comme $G/H$ est dénombrable, $V\cong L^{\BN}$, l'ensemble des suites à valeurs dans $L$, qui est un espace de Fréchet et $W\cong L^{(\BN)}$ l'ensemble des suites à valeurs dans $L$ à support compact, qui est un espace de type $LB$ (limite inductive d'espaces de Banach). Comme $L^{(\BN)}$ est le $L$-dual stéréotypique de $L^{\BN}$ on obtient sans peine que $V'\cong W$.

Pour le second point posons $V^+=\Ind_H^G\sigma^+$ et $W^+=\ind_H^G\sigma^{*,+}$. Alors $V^+\otimes L\subset V\cong L^{\BN}$ s'identife aux suites bornées à valeur dans $L$ qui est un espace de Smith (dual stéréotypique d'un banach). Ainsi, son dual s'identifie au $L$-banach 
$$
\left ( \wh{\bigoplus}_{\BN}\rmO_L\right ) \otimes_{\rmO_L}L=\wh{\bigoplus}_{\BN}L,
$$
où le complété est pour la norme $p$-adique, des suites de limite nulle. On en déduit que $(V^+\otimes_{\rmO_L}L)'\cong \wh{W}^+\otimes_{\rmO_L}L$ et comme cet isomorphisme identifie les boules unités on a $V^+\cong \wh{W}^+$.

Pour le dernier point, il suffit de remarquer que comme l'action de $G$ sur $G/H$ est transitive, $f\in V$, identifié à une fonction sur $G/H$, est $G$-bornée si et seulement si l'ensemble $\{f(s)\}_{s\in G/H}$ est borné, c'est-à-dire si et seulement s'il existe $n\geqslant 0$ tel que $p^nf\in V^+$, \ie $f\in V^+\otimes_{\rmO_L}L$. Ceci démontre le dernier point.

\end{proof}
Posons
$$
\CJ_0^{\lambda,(+)}\coloneqq \Ind_{KZ}^GW_{\lambda}^{(+)},\quad \wt{\CJ}_1^{\lambda,(+)}\coloneqq \Ind_{IZ}^GW_{\lambda}^{(+)},
$$
où les induites sont les induites classiques (sans condition sur le support). Comme précédemment, on définit $\Pi\colon \wt{\CJ}_1^{\lambda,(+)}\rightarrow \wt{\CJ}_1^{\lambda,(+)}$ la précomposition par $w_p$ et on pose $\CJ_1^{\lambda,(+)}\coloneqq(\wt{\CJ}_1^{\lambda,(+)})^{\Pi=-\id}$ de sorte que 
$$
\wt{\CJ}_1^{\lambda,(+)}\coloneqq\CJ_1^{\lambda,(+)}\oplus(\wt{\CJ}_1^{\lambda,(+)})^{\Pi=\id}.
$$

On commence par réinterpréter ces espaces en termes de fonctions sur l'arbre de Bruhat-Tits. Soit $\CT_{\bullet}$ l'arbre de Bruhat-Tits, le graphe défini par $\CT_0= G/KZ$ et $\CT_1=G/N$ et muni des deux flèches $\sfs,\sft:\CT_1\rightarrow \CT_0$ source et but que l'on va expliciter. La décomposition de Bruhat donne 
$$
G=P^+KZ\sqcup w_p P^+KZ,\quad P^+\coloneqq \begin{pmatrix}\Zp\setminus\{0\}&\Zp\\ 0 & 1\end{pmatrix}.
$$
Ainsi, comme l'application $G/IZ\rightarrow G/N$ est d'ordre $2$, tout élément $a\in G/N$ a pour antécédents $\{\wt{a},w_p\wt{a}\}\subset G/IZ$ choisis tels que $\wt{a}\in P^+KZ$. On pose alors 
$$
\sfs(a)\coloneqq \wt{\sfs}(\wt{a}),\ \sft(a)\coloneqq \wt{\sft}(\wt{a})=\wt{\sfs}(w_p\wt{a}).
$$
L'idée étant que l'on peut considérer $G/IZ$ comme l'ensemble des arêtes dédoublées et $\wt{\sfs},\wt{\sft}$ définissant une structure de graphe non orienté\footnote{Rappelons que si $\sfs,\sft\colon G_1\rightarrow G_0$ est un graphe (orienté) alors on lui associe un graphe non orienté en posant $\wt{G}_1\coloneqq G_1\sqcup G_1'$ où $G_1'=G_1$ et $\wt{\sfs},\wt{\sft}\colon \wt{G}_1\rightarrow G_0$ sont définis par $\restr{\wt{\sfs}}{G_1}=\sfs$, $\restr{\wt{\sfs}}{G_1'}=\sft$, $\restr{\wt{\sft}}{G_1}=\sft$, $\restr{\wt{\sft}}{G_1'}=\sfs$ } sur l'arbre de Bruhat-Tits. De plus, on obtient
$$
\CJ_i^{\lambda}\cong \sC(\CT_i,W_{\lambda}),
$$
où pour $i=1$ cet isomorphisme identifie les fonctions impaires sur les arêtes non orientées $G/IZ$ (au sens où changer le sens de l'arête change le signe de la fonction sur cette arête) aux fonctions sur les arêtes orientées. \emph{Le laplacien}  $\Delta \colon \CJ_1^{\lambda}\rightarrow \CJ_0^{\lambda}$ est défini sur $f\in\sC(\CT_1,W_{\lambda})$ pour $s\in \CT_0$ par
$$
(\Delta f)(s)\coloneqq \sum_{\substack{a\in \CT_1\\\sfs(a)=s}}f(a)-\sum_{\substack{a\in \CT_1\\\sft(a)=s}}f(a).
$$
On note $\Delta^+\colon \CJ_1^{\lambda,+}\rightarrow \CJ_0^{\lambda,+}$ sa restriction à $\CJ_1^{\lambda,+}$. Le lemme suivant découle immédiatement du lemme \ref{lem:indual} et de l'expression (\ref{eq:expart}) de $\partial$ :
\medskip
\begin{lemm}\label{lem:firdu}
Pour $i=0,1$, on a des isomorphismes $L$-linéaires topologiques $G$-équivariants
$$
(\CJ_i^{\lambda})'\cong\CI_i^{\lambda}.
$$
De plus, sous ces isomorphismes, $\Delta'=\partial$.
\end{lemm}
\medskip
\medskip
\begin{prop}\label{prop:lattdu}
Le laplacien 
$$
\Delta^{(+)}\colon \CJ_1^{\lambda,(+)}\xrightarrow{\Delta^{(+)}}\CJ_0^{\lambda,(+)},
$$
est un morphisme strict surjectif, et de plus :
\begin{itemize}
\itemb $\ker \Delta \cong (\St_{\lambda}^{\lalg})'$,
\itemb $(\ker\Delta^+)\otimes_{\rmO_L}L\cong \sD^{\lambda}$, \ie $\ker\Delta^+$ définit un réseau de $\sD^{\lambda}$.
\end{itemize}
\end{prop}
\medskip
\begin{proof}
Pour la première partie de la proposition, on peut supposer $\lambda=(0,1)$, quitte à utiliser une base de $W_{\lambda}$ et raisonner suivant les coordonnées.

On commence par montrer que $\Delta$ est strict pour la topologie faible et donc strict pour la topologie de Fréchet. Soit $s\in \CT_0$ et $p_s$ la semi-norme associée, pour $a\in \CT_1$, on note de même $p_a$ la semi-norme associée. Alors
$$
p_s(\Delta(f))\leqslant \sum_{\substack{a\in \CT_1\\\sft(a)=s}}\lvert f(a)\rvert - \sum_{\substack{a\in \CT_1\\\sfs(a)=s}}\lvert f(a)\rvert\leqslant 2(p+1)\sup_{\substack{a\\ s\in a}}p_a(f)
$$
où le $\sup$ se prend sur l'ensemble fini des arêtes voisines de $s$. Ainsi $\Delta$ est strict et en particulier d'image fermée (de même pour $\Delta^+$).

Comme $\Delta'=\partial$ par le lemme \ref{lem:firdu}, la surjectivité de $\Delta$ est une conséquence du lemme \ref{lem:indual} et de l'injectivité de $\partial$. On donne une preuve différente, en particulier pour avoir la surjectivité de $\Delta^+$. Comme on vient de montrer que l'image est fermée il suffit de montrer que son image contient une partie dense. Soit $s\in \CT_0$, on choisit une suite d'arêtes $\{a_n(s)\}_{n\in \BN}\subset \CT_1$ telles que $\sfs(a_0)=s$, et pour tout $n\geqslant 0$, $\sft(a_n)=\sfs(a_{n+1})$ (\ie les $a_n$ forment une chaine partant de $s$). Posons alors 
$$
f_s\coloneqq \sum_{n\geqslant 0}\delta_{a_n}.
$$
où pour une arête $a\in \CT_1$ (\resp un sommet $s\in \CT_0$), on note $\delta_a$ le dirac en $a$ ($\delta_s$ le dirac en $s$). On a alors $\Delta f_s=\delta_s$. Par linéarité, les sommes finies de diracs sont dans l'image et comme ce sous-espace est dense, on a montré que $\Delta$ est surjectif (de même\footnote{Remarquons que l'image de $\Delta^+$ est fermée parce que c'est l'image d'un compact par une application continue.} pour $\Delta^+$).

Montrons que $\ker \Delta \cong (\St_{\lambda}^{\lalg})'$ ce qu'il suffit de faire pour $\lambda=(0,1)$ quitte à tensoriser par $W_{\lambda}$. C'est alors un résultat classique sur les cochaînes harmoniques pour lequel on renvoie à \cite[Theorem 1.1]{AdS} (\cf aussi \cite[Lemme 12.4]{van}).

Comme on sait que $\sD^{\lambda}$ s'identifie aux vecteurs $G$-bornés de $(\St_{\lambda}^{\lalg})'$ (\cf lemme \ref{lem:gbded}), le dernier point est une conséquence du point $3$ du lemme \ref{lem:indual}.

\end{proof}

On termine maintenant la preuve de la proposition \ref{prop:latt}.
\begin{proof}
On commence par montrer que le conoyau de $\CI_0^{\lambda}\xrightarrow{\partial} \CI_1^{\lambda}$ est isomorphe à $\St_{\lambda}^{\lalg}$. D'après la proposition \ref{prop:lattdu} on a une suite exacte
$$
0\rightarrow (\St_{\lambda}^{\lalg})'\rightarrow \CJ_1^{\lambda}\xrightarrow{\Delta} \CJ_0^{\lambda}\rightarrow 0
$$
à laquelle on applique le $L$-dual stéréotypique pour obtenir par le lemme \ref{lem:indual}
$$
0\rightarrow \CI_0^{\lambda}\xrightarrow{\partial} \CI_1^{\lambda}\rightarrow \St_{\lambda}^{\lalg}\rightarrow 0.
$$
De même, on considère 
$$
0\rightarrow\sD^{\lambda}\rightarrow \CJ_1^{\lambda,+}\otimes_{\rmO_L}L\xrightarrow{\Delta} \CJ_0^{\lambda,+}\otimes_{\rmO_L}L\rightarrow 0,
$$
dont le dual par le lemme \ref{lem:indual} donne
$$
0\rightarrow \wh{\CI}_0^{\lambda,+}\otimes_{\rmO_L}L\xrightarrow{\Delta} \wh{\CI}_1^{\lambda,+}\otimes_{\rmO_L}L\rightarrow \sC^{\lambda}\rightarrow 0.
$$
Ainsi $I^{\lambda,+}\coloneqq \Coker \partial^+$ définit un $\rmO_L$-réseau stable par $G$ de $\St_{\lambda}^{\lalg}$ dont le complété $p$-adique est un $\rmO_L$-réseau $G$-stable de $\sC^{\lambda}$.
\end{proof}

\subsection{Déformations infinitésimales de $V_{\CL}^{\lambda}$}\label{subsec:definfgal}
Soit $\lambda\in P$ tel que $w(\lambda)>1$, on pose $D=\Sp(\lvert \lambda \rvert-2 )$ ; et, pour $n\geqslant 0$, $L_n\coloneqq L[T]/T^{n+1}$. Soit $\CL\in \BP^1(L)=L\cup\{\infty\}$.
\subsubsection{$\CL\neq \infty$}
Le cas $\CL\neq \infty$ est très similaire au cas cuspidal (\cf \cite{codonikis}). On pose 
$$
D_{\dR}\coloneqq (\Qpbar \otimes_{\Qp} M)^{\sG_{\Qp}}=Le_1\oplus Le_0.
$$
Soit $n\geqslant 0$ un entier, $D_{\CL,n}^{\lambda}$ désigne le $L_n$-$(\varphi,N,\sG_{\Qp})$-module filtré défini par $D_{\CL,n}^{\lambda}=L_n\otimes_L \Sp_L(\vert \lambda \rvert -2)$ en tant que $L_n$-$(\varphi,N,\sG_{\Qp})$-module et tel que $(\Qpbar \otimes_{\Qp}D_{\CL,n}^{\lambda})^{\sG_{\Qp}}=L_n\otimes_LD_{\dR}$ est muni de la filtration $\Fil_{\CL,n}^{\bullet}$ définie par
$$
\Fil_{\CL,n}^{i}\coloneqq
\begin{cases}
 L_n\otimes_L D_{\dR} & \text{ si }i\leqslant -\lambda_2+1\\
 L_n\otimes_L (e_1-(\CL+T)e_0) & \text{ si } -\lambda_2+2\leqslant i \leqslant -\lambda_1+1\\
 0 & \text{ si }i\geqslant -\lambda_1+2
\end{cases}
$$
Remarquons que pour $k\in \BN$ un entier tel que $0\leqslant k\leqslant n$, $T^kD_{\CL,n}^{\lambda}$ est isomorphe à $D_{\CL,n-k}^{\lambda}$ comme $L$-$(\varphi,N)$-module filtré.
\medskip
\begin{lemm}\label{lem:wafa}
Soit $n\geqslant 0$.
\begin{enumerate}
\item Le $L_n$-$(\varphi,N)$-module filtré $D_{\CL,n}^{\lambda}$ est faiblement admissible.
\item Les seuls sous $L$-$(\varphi,N)$-modules filtrés faiblement admissibles de $D_{\CL,n}^{\lambda}$ sont les $T^kD_{\CL,n}^{\lambda}$ pour $0\leqslant k \leqslant n$. 

\end{enumerate}
\end{lemm}
\medskip
\begin{proof}
La preuve est très similaire à celle de \cite[Lemme 4.3]{codonikis} avec une petite différence puisque $D_{\CL}^{\lambda}$ n'est pas irréductible comme $L$-$(\varphi,N)$-module. Posons $w\coloneqq w(\lambda)$. Soit $S\subset D_{\CL,n}^{\lambda}$ un sous-$L$-$(\varphi,N)$-module. Alors $S$ est de la forme $S=\Lambda_0e_0\oplus \Lambda_1e_1$ où $\Lambda_1$ et $\Lambda_0$ sont des $L$-sous-espaces vectoriels de $L_n$ tels que $\Lambda_1\subseteq \Lambda_0$, puisque $S$ est stable par $N$. Notons $d_i\coloneqq \dim_L\Lambda_i$ pour $i=0,1$, on a donc $d_1\leqslant d_0$. Comme la filtration n'a que deux crans, en simplifiant l'inégalité des pentes, le premier point revient à montrer
$$
w \dim_L(S\cap \Fil_{\CL,n}^{-\lambda_1})\leqslant d_0  \frac{w-1}{2}+d_1\frac{w+1}{2},
$$
et le second point à montrer que, si l'inégalité est une égalité, alors $\Lambda_1=\Lambda_0\cong T^kL_n$ pour $k\geqslant 0$. Or, $S\cap \Fil_{\CL,n}^{-\lambda_1}$ est constitué d'éléments de la forme $Pe_1+(\CL+T)Pe_0$ où $P\in \Lambda_1$ et $(\CL+T)P\in \Lambda_0$. En particulier :
$$
w \dim_L(S\cap \Fil_{\CL}^{-\lambda_1})\leqslant d_1 w=d_1\frac{w-1}{2}+d_1  \frac{w+1}{2}\leqslant d_0 \frac{w-1}{2}+d_1\frac{w+1}{2}.
$$
Ceci finit de prouver le premier point. Pour le second point, il est clair que si l'on a égalité, alors $d_1=d_0$, c'est-à-dire que $\Lambda\coloneqq \Lambda_0=\Lambda_1$ et donc $S=\Lambda \otimes_L \Sp(\lvert \lambda \rvert-2)$. Dans ce cas, pour un élément $Pe_1+(\CL+T)Pe_0\in S\cap \Fil_{\CL}^{-\lambda_1}$ on a $(\CL+T)P\in \Lambda$, donc $TP\in \Lambda$, c'est-à-dire que $\Lambda\subset L_n$ est un idéal car il est stable par $T$, \ie il existe $k\geqslant 0$ tel que $\Lambda = T^kL_n$.
\end{proof}
\subsubsection{$\CL=\infty$}
Le cas $\CL=\infty$ est différent du précédent : plutôt que de déformer la filtration, on déforme l'opérateur de monodromie \emph{et} le frobenius. Posons 
$$
\begin{gathered}
\BL_{\infty}\coloneqq  L\llbracket T_1,T_2\rrbracket /\langle T_1T_2\rangle,\\ \BL_{\bn}\coloneqq \BL/\langle T_1^{n_1+1}, T_2^{n_2+1}\rangle,\ \bn\coloneqq(n_1,n_2)\in \BN^2
\end{gathered}
$$

On pose comme précédemment
$$
D_{\dR}\coloneqq (\Qpbar \otimes_{\Qp} M)^{\sG_{\Qp}}=Le_1\oplus Le_0.
$$
Soit $\bn=(n_1,n_2)\in \BN^2$ une paire d'entiers positifs et $k\in \BZ$, on définit $\wt{\Sp}_{\bn}(k)=\BL_{\bn}e_0\oplus \BL_{\bn}e_1$ le $\BL_{\bn}$-$(\varphi,N)$-module par les formules
$$
\begin{cases}
\varphi(e_1)=p^{-\frac{k-1}{2}}(1-T_1)e_1\\
\varphi(e_0)=p^{-\frac{k+1}{2}}(1-T_1)^{-1}e_0,
\end{cases}
\quad
\begin{cases}
Ne_1=T_2e_0\\
Ne_0=0.
\end{cases}
$$
On doit vérifier la relation $p\varphi N =N\varphi$ qui, appliquée à $e_0$, est évidente ; appliquée à $e_1$, donne
$$
(p\varphi N -N\varphi)e_1 = p^{-\frac{k+1}{2}}T_2[(T_1+1)^{-1}-(1+T_1)]e_1=p^{-\frac{k+1}{2}}uT_2T_1e_1=0,
$$
où $u=2+\sum_{i=1}^{n_1}T_1^i\in \BL_{\bn}^{\times}$. Ceci justifie la relation $T_1T_2=0$ dans $\BL_{\infty}$.

Pour $\lambda\in P$, on définit le $\BL_{\bn}$-$(\varphi,N)$-module filtré $\BD_{\infty,\bn}^{\lambda}$ de $\BL_{\bn}$-$(\varphi,N)$-module sous-jacent $\wt{\Sp}_{\bn}(\lvert \lambda \rvert-2)$ et dont la filtration est définie par $\Fil_{\infty,n}^{\bullet}\coloneqq \BL_{\bn}\otimes_L \Fil_{\infty}^{\bullet}$ sur $(\Qpbar \otimes_{\Qp}\BD_{\infty,n}^{\lambda})^{\sG_{\Qp}}=L_{\bn}\otimes_LD_{\dR}$ 
 
Comme précédemment, pour $k_i\in \BN$ un entier tel que $0\leqslant k_i\leqslant n_i$ pour $i=1,2$, $T_1^{k_1}T_2^{k_2}D_{\infty,(n_1,n_2)}^{\lambda}$ est isomorphe à $D_{\infty,(n_1-k_1,n_2-k_2)}^{\lambda}$ comme $L$-$(\varphi,N)$-module filtré et on a le lemme suivant :
\medskip
\begin{lemm}
Soit $\bn=(n_1,n_2)\in \BN^2$.
\begin{enumerate}
\item Le $\BL_{\bn}$-$(\varphi,N)$-module filtré $\BD_{\infty,\bn}^{\lambda}$ est faiblement admissible.
\item Les seuls sous $L$-$(\varphi,N)$-modules filtrés faiblement admissibles de $\BD_{\infty,\bn}^{\lambda}$ sont les $T_1^{k_1}T_2^{k_2}\BD_{\infty,\bn}^{\lambda}$ pour $0\leqslant k_1 \leqslant n_1,\ 0\leqslant k_2 \leqslant n_2$. 

\end{enumerate}
\end{lemm}
\medskip
\begin{proof}
La preuve est très proche de celle du lemme \ref{lem:wafa}. Soit $S\subset \BD_{\infty,\bn}^{\lambda}$ un sous-$L$-$(\varphi,N)$-module. Alors $S=\Lambda_0e_0\oplus \Lambda_1e_1$ avec $\Lambda_0,\Lambda_1$ des $L$-sous-espaces vectoriels de $\BL_{\bn}$. Comme $S$ est stable par $N$, $T_2\Lambda_1\subset \Lambda_0$ et comme $S$ est stable par $\varphi$, $(1-T_1)\Lambda_1\subset \Lambda_1$ et $(1-T_1)^{-1}\Lambda_0\subset \Lambda_0$ ce qui implique que $\Lambda_0$ et $\Lambda_1$ sont stables par $T_1$. Soit $d_i\coloneqq \dim_L \Lambda_i$ pour $i=0,1$ et $d\coloneqq \dim_L(\Lambda_0\cap \Lambda_1)$, notons que $d\leqslant d_0,d_1$. 

Comme la filtration n'a que deux crans, on peut simplifier l'inégalité des pentes comme dans la preuve du lemme \ref{lem:wafa} ; le premier point revient donc à montrer
$$
w \dim_L(S\cap \Fil_{\infty,\bn}^{-\lambda_1})\leqslant d_0  \frac{w-1}{2}+d_1\frac{w+1}{2},
$$
et le second point à montrer qu'en cas d'égalité, $\Lambda_0=\Lambda_1=T_1^{k_1}T_2^{k_2}\BL_{\bn}$. Or, $S\cap \Fil_{\infty,\bn}^{-\lambda_1}$ est constitué d'éléments de la forme $Pe_1+Pe_0$ où $P\in \Lambda_1\cap \Lambda_0$ et donc
\begin{equation}\label{eq:exineq}
w\dim_L(S\cap \Fil_{\infty,\bn}^{-\lambda_1})=w d= d\frac{w-1}{2}+d\frac{w+1}{2}\leqslant d_0\frac{w-1}{2}+d_1\frac{w+1}{2}.
\end{equation}
Ceci fini de prouver le premier point. Pour le second point, il est clair que si l'on a égalité dans (\ref{eq:exineq}), alors $d=d_0=d_1$ et donc $\Lambda\coloneqq \Lambda_0=\Lambda_1$. Dans ce cas, comme $T_2\Lambda\subset \Lambda$, $\Lambda$ est stable par $T_2$ et on a déjà montré que $\Lambda$ est stable par $T_1$ donc c'est un idéal de $\BL_{\bn}$, \ie il existe $k_1,k_2\geqslant 0$ tels que $\Lambda=T_1^{k_1}T_2^{k_2}\BL_{\bn}$.
\end{proof}
Pour $n\geqslant 0$ un entier, on définit $D_{\infty,n}^{\lambda}\coloneqq \BD_{\infty,n}^{\lambda}/T_1\BD_{\infty,n}^{\lambda}$ comme $L_n$-$(\varphi,N,\sG_{\Qp})$-module filtré par l'isomorphisme $L_n\cong L\llbracket T_2 \rrbracket = \BL_{\infty}/T_1$ ; on définit de même $D_{\infty,n}^{\lambda,\cris}\coloneqq \BD_{\infty,n}^{\lambda}/T_2\BD_{\infty,n}^{\lambda}$ comme $L_n$-$(\varphi,N,\sG_{\Qp})$-module filtré par l'isomorphisme $L_n\cong L\llbracket T_1 \rrbracket = \BL_{\infty}/T_2$. Pour $\bn=(n_1,n_2)\in\BN$, on a un diagramme cartésien de $\BL_{\bn}$-$(\varphi,N,\sG_{\Qp})$-modules filtrés
\begin{equation}\label{eq:cartfont}
\begin{tikzcd}
\BD_{\infty,\bn}^{\lambda}\ar[r,"\mod T_2"]\ar[d,"\mod T_1" swap]& D_{\infty,n_1}^{\lambda,\cris}\ar[d]\\
D_{\infty,n_2}^{\lambda}\ar[r]& D_{\infty}^{\lambda}
\end{tikzcd}
\end{equation}
\subsubsection{Déformations}
On rappelle le lemme suivant (\cf \cite[Proposition 7.17]{dos}) : 
\medskip
\begin{lemm}\label{lem:extexp}
Soit $\lambda\in P$ et $\CL\in \BP^1(L)= L\cup\{\infty\}$.
\begin{enumerate}
\item Pour $\CL\neq \infty$, le sous-$L$-espace vectoriel de $\Ext^1_{L[\sG_{\Qp}]}(V_{\CL}^{\lambda},V_{\CL}^{\lambda})$ des extensions de Rham de déterminant $\zeta_{\lambda}\omega$ est de dimension $1$. 
\item Le sous-$L$-espace vectoriel $\Ext^1_{L[\sG_{\Qp}]}(V_{\infty}^{\lambda},V_{\infty}^{\lambda})$ des extensions de Rham de déterminant $\zeta_{\lambda}\omega$ est de dimension $2$.
\end{enumerate}
\end{lemm}
\medskip
\medskip
\begin{rema}
Dans la preuve de la proposition \cite[Proposition 7.17]{dos} il apparaît en outre que, dans le cas exceptionnel, le sous-espace des extensions cristallines de déterminant $\zeta_{\lambda}\omega$, est de dimension $1$ (ou de manière équivalente, le sous-espace des extensions semi-stables de déterminant $\zeta_{\lambda}\omega$ des extensions $W$ telles que $N\neq 0$ sur $\bD_{\pst}(W)$, est de dimension $1$). On le déduit aussi de l'existence de la déformation $\BV_{\infty,(1,1)}^{\lambda}$.
\end{rema}
\medskip
Pour $n\geqslant 0$ et $\bn=(n_1,n_2)\in \BN^2$, on pose pour $\CL\in L\cup \{\infty\}$ 
$$
V_{\CL,n}^{\lambda}\coloneqq \bV_{st}(D_{\CL,n}^{\lambda}[1]),\quad 
\BV_{\infty,\bn}^{\lambda}\coloneqq \bV_{st}(\BD_{\CL,\bn}^{\lambda}[1]), 
\quad V_{\infty,n}^{\lambda,\cris}\coloneqq \bV_{st}(D_{\infty,n}^{\lambda,\cris}[1]).
$$
Comme le foncteur $\bV_{\st}$ est exact, il préserve les limites finies et le diagramme cartésien \ref{eq:cartfont} donne un diagramme cartésien de $\BL_{\bn}$-représentations de $\sG_{\Qp}$ :
\begin{equation}\label{eq:cartgal}
\begin{tikzcd}
\BV_{\infty,\bn}^{\lambda}\ar[r,"\mod T_2"]\ar[d,"\mod T_1" swap]& V_{\infty,n_1}^{\lambda,\cris}\ar[d]\\
V_{\infty,n_2}^{\lambda}\ar[r]& V_{\infty}^{\lambda}
\end{tikzcd}
\end{equation}
Soit $\CB$ le bloc correspondant à la réduction\footnote{Rappelons que la \emph{réduction} d'une $L$-représentation $V$ de $\sG_{\Qp}$ est la semi-simplification de $V^+/\varpi_L$ où $V^+\subset V$ est un $\rmO_L$-réseau stable par $\sG_{\Qp}$. Cette définition est indépendante du réseau $V^+$, grâce à la semi-simplification.} de $V=V_{\CL}^{\lambda}$ et posons $R\coloneqq R_{\CB}^{\ps,\zeta_{\lambda}}[1/p]$. Alors $V$ définit un idéal maximal $\fkm_{\CL}^{\lambda}\subset R$. Soit $\wh{R}_{\CL}$ le complété $\fkm_{\CL}^{\lambda}$-adique de $R$. Cet anneau est muni d'une représentation universelle $\wh{\rho}_{\CL}\colon \sG_{\Qp}\rightarrow \GL_2(\wh{R}_{\CL})$.

Soit maintenant $\wh{R}_{\CL,\dR}$ le quotient de $\wh{R}_{\CL}$ classifiant les représentations de Rham, \ie $\wh{R}_{\CL,\dR}\coloneqq \wh{R}_{\CL}/I$ où $I=\cap \fka$ et $\fka$ parcourt les idéaux de $\wh{R}_{\CL}$ tels que $\wh{R}_{\CL}/\fka$ soit de dimension finie sur $L$ et $(\wh{R}_{\CL}/\fka)\otimes \wh \rho_{\CL}$ soit de Rham.

\medskip
\begin{prop}\label{prop:rdralll}
Pour $\CL\neq \infty$, on a $\wh{R}_{\CL,\dR}\cong L_{\infty}$ et $\wh{R}_{\infty,\dR}/T_2\cong\wh{R}_{\infty,\dR}/T_1\cong L_{\infty}$
\end{prop}
\medskip
\begin{proof}
Supposons d'abord que $\CL\neq \infty$, la preuve est alors la même que celle de \cite[Proposition 4.5]{codonikis}. En effet, d'après le premier point du lemme \ref{lem:extexp}, $V_{\CL,1}^{\lambda}$ est l'unique extension non triviale de $V_{\CL}^{\lambda}$ par $V_{\CL}^{\lambda}$ qui est de déterminant $\zeta_{\lambda}\omega$ et qui soit de de Rham. On en déduit que $\wh{R}_{\CL,\dR}$ est un anneau local régulier de dimension $1$. 

Si $\CL=\infty$, le second point du lemme \ref{lem:extexp} donne que $\wh{R}_{\infty,\dR}$ est un quotient de $L\llbracket T_1,T_2\rrbracket$. Alors $\wh{R}_{\sL,\dR}/T_1$ est un quotient de $L_n$ et l'existence de $V_{\infty,n}^{\lambda}$ pour tout $n\geqslant 0$ fournit une suite compatible de morphismes $\wh{R}_{\infty,\dR}/T_1\rightarrow L_n$ ; ils sont surjectifs puisque $V_{\infty,1}^{\lambda}$ n'est pas scindé. Le même raisonnement pour $\wh{R}_{\infty,\dR}/T_2$ donne le résultat.
\end{proof}
\medskip
\begin{prop}\label{prop:rinf}
On a $\wh{R}_{\infty,\dR}\cong\BL_{\infty}=L\llbracket T_1,T_2\rrbracket/\langle T_1T_2\rangle$.
\end{prop}
\medskip
\begin{proof}
On note $R=\wh{R}_{\infty,\dR}$, $D\coloneqq D_{\st}(\wh{\rho}_{\infty,\dR})$ le $R$-$(\varphi,N)$-module filtré de la déformation infinitésimale universelle $\wh{\rho}_{\infty,\dR}\colon \sG_{\Qp}\rightarrow \GL_2(R)$ et $\varphi_D$ le frobenius sur $D$. Notons que $D$ est un $R$-module libre de rang $2$. On veut montrer qu'il existe un morphisme $\BL_{\infty}\rightarrow R$ tel que $\BD_{\infty}^{\lambda}\otimes_{\BL_{\infty}}R\cong D$. Rappelons que $R$ est un quotient de $\BL\llbracket T_1,T_2\rrbracket$ et notons 
$$
\BD_1\coloneqq \varprojlim_{\bn} \BD_{\infty,\bn}^{\lambda}/T_2,\quad \BD_2\coloneqq \varprojlim_{\bn} \BD_{\infty,\bn}^{\lambda}/T_1;
$$
ils définissent des $L_{\infty}$-$(\varphi,N)$-modules. D'après la proposition \ref{prop:rdralll}, on a des diagrammes commutatifs
$$
\begin{array}{lr}
\begin{tikzcd}
R\ar[r,"\mod T_2"]\ar[d,"\mod T_1" swap]& R_1\cong L\llbracket T_1\rrbracket\ar[d]\\
R_2\cong L\llbracket T_2\rrbracket\ar[r]& L
\end{tikzcd}
&
\begin{tikzcd}
D\ar[r,"\mod T_2"]\ar[d,"\mod T_1" swap]& D_1\cong \BD_1\ar[d]\\
D_2\cong \BD_2\ar[r]& D_{\infty}^{\lambda}
\end{tikzcd}
\end{array}
$$
Soient $\tilde{e}_0,\tilde{e}_1$ des relèvements de $e_0,e_1\in \BD_1$ ; d'après le lemme de Nakayama $$
D=R\tilde{e}_0\oplus R\tilde{e}_1.
$$
Comme $D$ est de rang $2$ sur $R$ et que $N\neq 0$ sur son quotient $D_2$, on sait que $N\neq0$ et $N^2=0$ sur $D$, donc, d'après le lemme de Nakayama, $\ker N = R\tilde{e}_0$. Comme $\ker N$ est stable par $\varphi_D$, on obtient $u\in L\llbracket T_1, T_2 \rrbracket^{\times}$ tel que $\varphi_D(\tilde{e}_0)=u \tilde{e}_0$. De plus, il existe $P\in L\llbracket T_1,T_2\rrbracket$ tel que $N\tilde{e}_1=T_2P\tilde{e}_0$ puisque $N\equiv 0 \mod T_2$. Alors, $\varphi_D(\tilde{e}_1)=\beta \tilde{e}_1+\gamma\tilde{e}_0$ et quitte à faire le changement de base $\tilde{e}_1\leftrightarrow \tilde{e}_1-\gamma u^{-1}\tilde{e}_0$, on peut supposer que $\gamma=0$ et $\beta = pu^{-1}$.

Pour résumer, on a justifié que $D=R\tilde{e}_0\oplus R\tilde{e}_1$,  $u\in L\llbracket T_1, T_2 \rrbracket^{\times}$ et $P\in L\llbracket T_1, T_2 \rrbracket$ tels que 
$$
\begin{array}{cc}
 
\begin{cases}
 \varphi_D(\tilde{e}_1)=pu^{-1}\tilde{e}_1\\
 \varphi_D(\tilde{e}_0)=u\tilde{e}_0
\end{cases}
&
\begin{cases}
 N \tilde{e}_1=PT_2\tilde{e}_0\\
 N\tilde{e}_0=0
\end{cases}

\end{array}
$$
Montrons que $v=(u-u^{-1})\in T_1R^{\times}$. Pour ce faire, il suffit de montrer que $v\equiv 2T_1 \mod \langle T_1^2,T_2\rangle$ ; or, on a $v\equiv (1-T_1)-(1-T_1)^{-1}\equiv 2T_1 \mod \langle T_1^2,T_2\rangle$. Ainsi, on obtient
$$
(p\varphi_D N-N\varphi_D)\tilde{e}_1\in T_1T_2PR^{\times},
$$
et donc $R$ est un quotient de $R'\coloneqq L\llbracket T_1,T_2\rrbracket/\langle T_1T_2P\rangle$. Enfin, le morphisme $\BL_{\infty}\rightarrow R'$ défini par $T_1\mapsto T_1$ et $T_2\mapsto PT_2$ est tel que $\BD_{\infty}^{\lambda}\otimes_{\BL_{\infty}}R\cong D$.
\end{proof}

\subsection{Déformations infinitésimales de $\Pi_{\CL}^{\lambda}$}\label{subsec:definfaut}
Pour $n\geqslant 0$ et $\bn=(n_1,n_2)\in \BN^2$, on pose pour $\CL\in L\cup \{\infty\}$ 
$$
\Pi_{\CL,n}^{\lambda}\coloneqq \bPi(V_{\CL,n}^{\lambda}),\quad \bPi_{\infty,\bn}^{\lambda}\coloneqq \bPi(\bV_{\infty,\bn}^{\lambda}),\quad \Pi_{\infty,n}^{\lambda,\cris}\coloneqq \bPi(V_{\infty,n}^{\lambda,\cris}).
$$
À partir du diagramme cartésien (\ref{eq:cartfont}), un diagramme de $\BL_{\bn}$-représentations de $G$
\begin{equation}\label{eq:cartdiaggl}
\begin{tikzcd}
\bPi_{\infty,\bn}^{\lambda}\ar[r,"\mod T_2"]\ar[d,"\mod T_1" swap]&\Pi_{\infty,n_1}^{\lambda,\cris}\ar[d]\\
\Pi_{\infty,n_2}^{\lambda}\ar[r]& \Pi_{\infty}^{\lambda}
\end{tikzcd}
\end{equation}
\medskip
\begin{rema}
Il n'est pas immédiat que ce diagramme est cartésien puisque $\bPi$ n'est pas exact au sens strict du terme. Il est sûrement possible de montrer qu'il l'est à partir de \ref{eq:cartgal} en appliquant $\bV$ au diagramme \ref{eq:cartdiaggl} et en se ramenant à des représentations $\Pi$ qui satisfont $\bPi\circ \bV(\Pi)=\Pi$.
\end{rema}
\medskip

\subsubsection{Vecteurs localement algébriques}
Dans ce numéro on calcule les vecteurs localement algébriques de $\Pi_{\CL,n}^{\lambda}$ et $\bPi_{\infty,\bn}^{\lambda}$. On commence par rappeler la proposition suivante (\cf \cite[Théorème 4.1]{codonikis}) :
\medskip
\begin{prop}\label{prop:drlalg}
Soit $E$ une algèbre locale de corps résiduel $L$ et $V$ une $E$-représentation de $\sG_{\Qp}$ de rang $2$ sur $E$ qui est absolument irréductible. Alors $\bPi(V)^{\lalg}$ est dense dans $\bPi(V)$ si et seulement si $V$ est de Rham à poids de Hodge-Tate distincts.
\end{prop}
\medskip
Rappelons que\footnote{En raison des nombreuses décorations autour de $\Pi$, on s'autorise à placer l'exposant \og $\lalg$\fg en pré-exposant.} 
$$
\prescript{\lalg}{}{\Pi}_{\CL}^{\lambda}=
\begin{cases}
\St_{\lambda}^{\lalg} & \CL\neq \infty,\\
\wt{\St}_{\lambda}^{\lalg} & \CL=\infty,
\end{cases}
$$
où $\wt{\St}_{\lambda}^{\lalg}$ est l'unique extension de $W^{*}_{\lambda}$ par $\St_{\lambda}^{\lalg}$. Notons que $\St_{\lambda}^{\lalg}$ est dense dans $\prescript{\lalg}{}{\Pi}_{\infty}^{\lambda}$ puisque cette représentation est un quotient de $\sC^{\lambda}$ d'après la proposition \ref{prop:denstoquo}.
\medskip
\begin{lemm}\label{lem:localg}
\begin{enumerate}
\item Si $n\geqslant 0$ est un entier et $\CL\neq \infty$ alors
$$
\prescript{\lalg}{}{\Pi}_{\CL,n}^{\lambda}\cong\St_{\lambda}^{\lalg}\otimes_L L_n.
$$
\item Si $\bn=(n_1,n_2)\in \BN^2$ alors on a un diagramme commutatif de $\BL_{\bn}$-représentations de $G$
\begin{equation}\label{eq:cartlalg}
\begin{tikzcd}
\prescript{\lalg}{}{\bPi}_{\infty,\bn}^{\lambda}\ar[r,"\mod T_2"]\ar[d,"\mod T_1" swap]&\prescript{\lalg}{}{\Pi}_{\infty,n_1}^{\lambda,\cris}\ar[d]\\
\prescript{\lalg}{}{\Pi}_{\infty,n_2}^{\lambda}\ar[r]& \wt{\St}_{\lambda}^{\lalg}
\end{tikzcd}
\end{equation}
où $\prescript{\lalg}{}{\Pi}_{\infty,n_2}^{\lambda}\cong\wt{\St}_{\lambda}^{\lalg}\otimes_L L_n$.
\end{enumerate}
\end{lemm}
\medskip
\begin{proof}
Pour le premier point, la preuve se fait comme \cite[Remarque 4.7 (i)]{codonikis} que l'on rappelle brièvement. D'après \cite[Remarque 4.2 (ii)]{codonikis}, si $\Pi$ est une $L$-représentation unitaire de $G$ dont toutes les composantes de Jordan-Hölder sont isomorphes à $\Pi_0$ et $\Pi_0^{\lalg}$ est irréductible et dense dans $\Pi_0$ alors $\Pi^{\lalg}$ est dense dans $\Pi$ si et seulement si $\lg(\Pi)=\lg(\Pi^{\lalg})$, ce qui permet d'obtenir le résultat grâce à \ref{prop:drlalg}.

Pour le second point, on commence par calculer $\prescript{\lalg}{}{\Pi}_{\infty,n_2}^{\lambda}$, ce qui se fait comme précédemment en adaptant \cite[Remarque 4.2 (ii)]{codonikis} de la façon suivante : si $\Pi$ est une $L$-représentation unitaire de $G$ dont toutes les composantes de Jordan-Hölder sont isomorphes à $\Pi_0$ et $\Pi_0^{\lalg}$ est de longueur $2$ et dense dans $\Pi_0$ alors $\Pi^{\lalg}$ est dense dans $\Pi$ si et seulement si $2\lg(\Pi)=\lg(\Pi^{\lalg})$. On adapte alors la fin de l'argument de \cite[Remarque 4.7 (i)]{codonikis} pour conclure la preuve.

\end{proof} 

\medskip
\begin{rema}
\ 
\begin{itemize}
 \itemb Notons que le diagramme (\ref{eq:cartlalg}) est cartésien si (\ref{eq:cartdiaggl}) est cartésien puisque prendre les vecteurs localement algébriques (et les vecteurs localement analytiques) est un foncteur exact à gauche donc commute aux limites finies.
\itemb Si l'on disposait de théorèmes précis pour la correspondance de Langlands $p$-adique en famille (ce qui ne saurait tarder, \cf Colmez-Rodrigues, Dotto-Emerton-Gee), on pourrait décrire explicitement $\prescript{\lalg}{}{\Pi}_{\infty,n_1}^{\lambda,\cris}$ comme une induite localement algébrique explicite, puisque c'est une représentation ordinaire. Plus précisément, on aurait
$$
\prescript{\lalg}{}{\Pi}_{\infty,n_1}^{\lambda,\cris}\cong \Ind_B^G(x^{\lambda_1}(1-T_1)^{-v_p(x)}\otimes x^{\lambda_2-1}(1-T_1)^{v_p(x)})^{\lalg}\otimes_L\lvert \det \rvert^{\tfrac{\lvert \lambda \rvert-1}{2}}.
$$
La démonstration du lemme suivant serait alors plus immédiate.
\end{itemize}

\end{rema}
\medskip
\begin{lemm}\label{lem:stfac}
Soit $n\in \BN$, tout morphisme $\St_{\lambda}^{\lalg} \rightarrow \Pi_{\infty,n}^{\lambda,\cris}$ se factorise par
$$
\St_{\lambda}^{\lalg} \rightarrow T_1^{n}\Pi_{\infty,n}^{\lambda,\cris}\incl \Pi_{\infty,n}^{\lambda,\cris}.
$$
\end{lemm}
\begin{proof}
Rappelons que $\St_{\lambda}^{\lalg}$ est un quotient de 
$$
B_{\lambda}^{\lalg}\coloneqq \Ind_B^G(x^{\lambda_1}\otimes x^{\lambda_2-1})^{\lalg}\otimes_L\lvert \det \rvert^{\tfrac{\lvert \lambda \rvert-1}{2}}.
$$
Notons $\delta$ le caractère de $T$ défini par $\delta(\diag(x_1,x_2))\coloneqq x_1^{\lambda_1} x_2^{\lambda_2-1}\lvert x_1x_2 \rvert^{\tfrac{\lvert \lambda \rvert-1}{2}}$. Par la réciprocité de Frobenius, on se ramène à montrer que toute application $\delta\rightarrow \Pi_{\infty,n}^{\lambda,\cris}$ $B$-équivariante se factorise par $T_1^{n}\Pi_{\infty,n}^{\lambda,\cris}$. Pour cela, il suffit de démontrer que 
\begin{equation}\label{eq:cristfor}
g_p\coloneqq 
\begin{pmatrix}
p & 0\\ 0 & 1
\end{pmatrix},\ v\in \prescript{\lalg}{}{\Pi}_{\infty,n}^{\lambda,\cris}\text{, alors } g_p\cdot v\in (1-T_1)^{\pm1}Lv
\end{equation}
En effet si c'est le cas, l'image de $\delta\rightarrow \prescript{\lalg}{}{\Pi}_{\infty,n}^{\lambda,\cris}$ est un élément $v\in \Pi_{\infty,n}^{\lambda,\cris}$ tel que $g_p\cdot v \in Lv$, ce qui donnerait $Lv=(1-T_1)^{\pm 1}Lv$ et donc $v\in T_1^{n}\Pi_{\infty,n}^{\lambda,\cris}$.

Pour démontrer (\ref{eq:cristfor}) on utilise le fait que $\prescript{\lalg}{}{\Pi}_{\infty,n}^{\lambda,\cris}\cong (\Delta_{\dif}^-)^{\sU(\fkg)-\finie}$ en tant que $L_n$-représentation de $B$ comme dans la preuve de \cite[Théorème 4.2]{codonikis}. L'action de $g_p$ est alors donnée par $\varphi$, d'après \cite[Proposition I.3.6]{colp}. On en déduit que si $v\in \prescript{\lalg}{}{\Pi}_{\infty,n}^{\lambda,\cris}$, alors $g_p\cdot v\in (1-T_1)^{\pm1}Lv$.
\end{proof}
\medskip
\begin{coro}\label{coro:speg}
\begin{enumerate}
\item Pour $n\geqslant 0$ un entier, la représentation $\Pi_{\CL,n}^{\lambda}$ est un quotient de $\sC^{\lambda}$.
\item Pour $\bn=(n_1,n_2)\in \BN^2$, si $\St_{\lambda}^{\lalg} \rightarrow \bPi_{\infty,\bn}^{\lambda}$ est un morphisme d'image dense, alors $n_1=0$.
\end{enumerate}
\end{coro}
\medskip
\begin{proof}
Pour le premier point, par la proposition \ref{prop:denstoquo}, il suffit de construire une inclusion $\St_{\lambda}^{\lalg}\incl \Pi_{\CL,n}^{\lambda}$ d'image dense. Comme $V_{\CL,n}^{\lambda}$ est de Rham, d'après la proposition \ref{prop:drlalg}, $\prescript{\lalg}{}{\Pi}_{\CL,n}^{\lambda}\incl \Pi_{\CL,n}^{\lambda}$ est d'image dense. Par le lemme \ref{lem:localg} et l'irréductibilité de $\St_{\lambda}^{\lalg}$, on obtient
$$
\Hom_G(\St_{\lambda}^{\lalg},\prescript{\lalg}{}{\Pi}_{\CL,n}^{\lambda})\cong L_n
$$
On choisit $f\colon \St_{\lambda}^{\lalg}\rightarrow \Pi_{\CL,n}^{\lambda}$ correspondant à $\lambda\in L_n^{\times}$ dont la projection sur $(\Pi_{\CL,n}^{\lambda}/T)^{\lalg}$ est non nulle ; elle est donc injective et dense dans $\Pi_{\CL,n}^{\lambda}/T$. On en déduit que l'image de $f$ est dense dans $\Pi_{\CL,n}^{\lambda}$ car ses seules sous-$L$-représentations (topologiques) de $G$ sont les $T^k\Pi_{\CL,n}^{\lambda}$ pour $0\leqslant k \leqslant n$.

Pour le second point, supposons que l'on ait un morphisme $\St_{\lambda}^{\lalg} \rightarrow \bPi_{\infty,\bn}^{\lambda}$ d'image dense et montrons que $n_1=0$. D'après la proposition \ref{prop:denstoquo}, on a une suite de surjections $\sC^{\lambda}\twoheadrightarrow  \bPi_{\infty,\bn}^{\lambda} \twoheadrightarrow \Pi_{\infty,n_1}^{\lambda,\cris}$. Alors, d'après la proposition \ref{prop:denstoquo} on a un morphisme injectif $\St_{\lambda}^{\lalg}\incl \Pi_{\infty,n_1}^{\lambda,\cris}$ d'image dense. Ainsi, d'après le lemme \ref{lem:stfac}, $n_1=0$.
\end{proof}

\Subsubsection{Complétés de longueur finie de $\St_{\lambda}^{\lalg}$ }
\medskip
\begin{prop}\label{prop:complst}
Soit $\lambda\in P$ tel que $w(\lambda)>1$. Soit $\Pi$ une $L$-représentation unitaire de $G$ de caractère central $\zeta_{\lambda}$ muni d'un morphisme $\St_{\lambda}^{\lalg}\rightarrow \Pi$ d'image dense telle que les composantes de Jordan-Hölder\footnote{En fait, il suffit de supposer qu'une seule des composantes de Jordan-Hölder soit $\Pi_{\CL}^{\lambda}$.} de $\Pi$ soient tous $\Pi_{\CL}^{\lambda}$ avec multiplicité $n$, pour $\CL\in L\cup\{\infty\}$. Alors $\Pi\cong \Pi^{\lambda}_{\CL,n}$.
\end{prop}
\medskip
\begin{proof}
La preuve est très similaire à \cite[Proposition 4.8]{codonikis}. On rappelle brièvement l'argument.

Soit $J_{\CL}$ l'enveloppe injective de $\Pi_{\CL}^{\lambda}$ dans la catégorie des limites inductives des $L$-représentations unitaires de $G$. Alors $\End J_{\CL}=\wh{R}_{\CL}$ d'après \cite[Corollary 6.26]{patu} et comme le bloc de $\Pi_{\CL}^{\lambda}$ dans la catégorie des représentations unitaires de $G$ n'a qu'un objet (puisque $V_{\CL}^{\lambda}$ est absolument irréductible, \cf \cite[Corollary 6.11]{patu}) on a 
$$
\Pi\cong (\Hom_G(J_{\CL}',\Pi')\otimes_{\wh{R}_{\CL}}J_{\CL}')'.
$$
Comme $\St_{\lambda}^{\lalg}$ est dense dans $\Pi$, $\Pi^{\lalg}$ est dense dans $\Pi$ et donc $\Hom_G(J_{\CL}',\Pi')$ est en fait un $\wh{R}_{\CL,\dR}$-module.
\begin{itemize}
\itemb Si $\CL\neq \infty$ la proposition \ref{prop:rdralll} implique que
$$
\Hom_G(J_{\CL}',\Pi')\cong \bigoplus_{i\in I}L_{\infty}/T^{n_i},\quad \Pi\cong\bigoplus_{i\in I}\Pi_{\CL,n_i}^{\lambda},
$$
où $I$ est un ensemble fini et $n_i\in \BN$ pour tout $i\in I$. Comme dans la preuve de \cite[Proposition 4.8]{codonikis}, en utilisant la densité de $\St_{\lambda}^{\lalg}$, on conclut que $\lvert I \rvert =1$.
\itemb Si $\CL=\infty$ la proposition \ref{prop:rinf} implique que
$$
\Hom_G(J_{\CL}',\Pi')\cong \bigoplus_{i\in I}\BL_{\infty}/\langle T_1^{n_i}, T_2^{m_i}\rangle ,\quad \Pi\cong\bigoplus_{i\in I}\bPi_{\infty,(n_i,m_i)}^{\lambda},
$$
où $I$ est un ensemble fini et $(n_i,m_i)\in \BN^2$ pour tout $i\in I$. Comme précédemment on conclut que $\lvert I \rvert =1$ ce qui donne que $\St_{\lambda}^{\lalg}$ est dense dans $\Pi\cong \bPi_{\infty,(n,m)}^{\lambda}$ avec $(n,m)\in \BN^2$ et donc $n=0$ d'après le point $2$ du corollaire \ref{coro:speg}. Donc $\Pi\cong \Pi_{\infty,m}^{\lambda}$.
\end{itemize}

\end{proof}

\subsection{Complétion $\CB$-adique et anneaux de Kisin}\label{subsec:Bcomkis}
\subsubsection{Définitions}\label{subsubsec:defibcomp}
Rappelons que (\cf \cite[\S 2.1]{codonikis}) si $X$ est un $\rmO_L[G]$-module topologique, on définit $X_{\CB}$, \emph{le complété $\CB$-adique de $X$}, comme la limite inverse des quotients de $X$ appartenant à $\Tors_{\CB} G$. Le complété profini $\wh{X}$ de $X$ admet une décomposition
$$
\wh{X}=\prod_{\CB}X_{\CB},
$$
où $\CB$ décrit l'ensemble des blocs des $k_L$-représentations de $G$. On appelle dans ce cas $X_{\CB}$ le \emph{complété $\CB$-adique de $X$}. En particulier, on note $I^{\lambda,+}_{\CB}$ le complété $\CB$-adique de $I^{\lambda,+}$ {\cf proposition \ref{prop:latt}). On a
$$
\wh{I}^{\lambda,+}=\prod_{\CB}I^{\lambda,+}_{\CB},
$$ 
et on pose $\sC_{\CB}^{\lambda}\coloneqq I^{\lambda,+}_{\CB}\otimes_{\rmO_L}L$. 

On écrit $I^{\lambda,+}_{\CB}=\varprojlim_i \Pi^+_i$ où les $\Pi^+_i\in \Tors_{\CB}^{\zeta_{\lambda}}G$ et on pose
$$
\sigma_{\CB}^{\lambda}\coloneqq (\varprojlim_i \bV(\Pi^+_i))\otimes_{\rmO_L}L.
$$
Par l'isomorphisme du théorème \ref{thm:patu} $R_{\CB}^{\lambda,+}$ agit sur $I^{\lambda,+}_{\CB}$, on note $R_{\CB}^{\lambda,+}$ le quotient à travers lequel $R^{\ps,\zeta_{\lambda}}_{\CB}$ agit sur $I^{\lambda,+}_{\CB}$ et on pose $R_{\CB}^{\lambda}\coloneqq R_{\CB}^{\lambda,+}[1/p]$ (ce sont les \emph{anneaux de Kisin} définis dans l'introduction). Ainsi $\sC_{\CB}^{\lambda}$ est un $R_{\CB}^{\lambda}$-module et par la fonctorialité de $\bV$, on obtient une action de $R_{\CB}^{\lambda}$ sur $\sigma_{\CB}^{\lambda}$. Le but de cette sous-section est de montrer que $R_{\CB,M}^{\lambda}=T_{M,\CB}^{\lambda}$.
\subsubsection{Propriétés de finitude}
\medskip
\begin{lemm}\label{lem:rhotf}
Le $R_{\CB}^{\lambda}$-module $\rho_{\CB}^{\lambda}$ est de type fini.
\end{lemm}
\medskip
\begin{proof}
On suit la preuve de \cite[Lemme 5.2]{codonikis}. Il suffit de prouver que $\bV(I^{\lambda,+}_{\CB})$ est un $R^+=R^{\ps,\zeta_{\lambda}}_{\CB}$-module de type fini. Par le lemme de Nakayama (topologique), il suffit de montrer que $\bV(I^{\lambda,+}_{\CB})\otimes_{\rmO_L}k_L$ est de type fini. Or,
$$
\bV(I^{\lambda,+}_{\CB})\otimes_{\rmO_L}k_L=\bV(I^{\lambda,+}_{\CB}\otimes_{\rmO_L}k_L)=\bV(({I^{\lambda,+}\otimes_{\rmO_L}k_L})_{\CB}).
$$
En effet, pour la première égalité, $I^{\lambda,+}_{\CB}=\varprojlim_i \Pi^+_i$ où les morphismes de transition sont surjectifs, alors, en considérant la suite exacte
$$
0\rightarrow \Pi^+_i\xrightarrow{\varpi_L}\Pi^+_i\rightarrow \Pi^+_i\otimes_{\rmO_L}k_L\rightarrow 0,
$$
à laquelle on applique le foncteur exact $\bV$ (\cf \cite[Théorème IV. 2.24]{colp}), puis on passe à la limite projective ce qui donne encore une suite exacte puisque $R^1\lim \bV(\Pi_i^+)=0$ ($\Pi_{i+1}^+\rightarrow \Pi_i^+$ étant surjectif et $\bV$ exacte, $\bV(\Pi_{i+1}^+)\rightarrow \bV(\Pi_{i}^+)$ est surjectif). L'égalité $I^{\lambda,+}_{\CB}\otimes_{\rmO_L}k_L=({I^{\lambda,+}\otimes_{\rmO_L}k_L})_{\CB}$ découle de l'exactitude de la $\CB$-complétion (\cf \cite[Corollaire 2.3]{codonikis}). 

On conclut alors la preuve par le corollaire par \cite[Corollaire 2.8]{codonikis} comme dans la preuve de \cite[Lemme 5.2]{codonikis}.
\end{proof}

\medskip
\begin{lemm}\label{lem:finitcompdef}
Soit $n\geqslant 0$ un entier, $\fkm_x\subset R_{\CB}^{\ps,\zeta}[1/p]$ un idéal maximal dont on note $L_x$ le corps résiduel. Le $L_x[G]$-module $\sC_{\CB}^{\lambda}/\fkm_x^{n+1}$ est de longueur finie.
\end{lemm}
\medskip
\begin{proof}
On note 
$$
R^+\coloneqq R_{\CB}^{\ps,\zeta},\quad R\coloneqq R^+\bigl[\tfrac 1 p \bigr],\quad \Pi_{x,n}\coloneqq\sC_{\CB}^{\lambda}/\fkm_x^{n+1}.
$$
Par \cite[Corollary 1.6]{codopa} il suffit de montrer que $\Pi_{x,n}$ est résiduellement de longueur finie. Soit $\tilde{\fkn}_{n,x}\coloneqq R^+\cap \fkm_{x,n}$ et $\fkn_{x,n}\subset R^+$ l'idéal engendré par $\tilde{\fkn}_{x,n}$ et $\varpi_L$. Alors, il suffit de montrer que $\pi_{x,n}\coloneqq I_{\CB}^{\lambda,+}/\fkn_{x,n}$ est de longueur finie.
D'après \cite[Proposition 2.2]{codonikis} on a
\begin{equation}\label{eq:itop}
 I_{\CB}^{\lambda,+}\cong (P_{\CB}^{\zeta}\otimes_{E_{\CB}^{\zeta}}\Hom_G(P_{\CB}^{\zeta},(I^{\lambda,+})^{\vee}))^{\vee}.
\end{equation}

où $P_{\CB}^{\zeta}\coloneqq \bigoplus_{\pi\in \CB}P_{\pi}^{\zeta}$ et $E_{\CB}^{\zeta}\coloneqq\End(P_{\CB})$ où $P_{\pi}$ est l'enveloppe projective de $\pi^{\vee}$ (définie dans \cite[\S 2]{pas}, \cf aussi \num \ref{subsubsec:blocth} ) et $(\ \cdot\ )^{\vee}\coloneqq \Hom_G(\ \cdot\ ,L/\rmO_L)$ est le dual de Pontryagin. Mais d'après \cite[Corollary 6.7]{patu}$, P_{\CB}^{\zeta}/\fkn_{x,n}\cong \pi^{\vee}$ où $\pi$ est un $\rmO_L[G]$-module de longueur finie, la formule \ref{eq:itop} donne, par compacité et projectivité de $P_{\CB}^{\zeta}$ :
$$
\pi_{x,n}\cong (\pi^{\vee}\otimes_{E_{\CB}^{\zeta}}\Hom_G(I^{\lambda,+}/\varpi_L,\pi))^{\vee}.
$$
Donc il reste à montrer que $\Hom_G(I^{\lambda,+}/\varpi_L,\pi)$ est de $k_L$-dimension finie. L'injection $I^{\lambda,+}\incl \ind_{KZ}^G\st_{\lambda}^+$ (\cf (\ref{eq:sttoi})) combinée à la réciprocité de Frobenius donne une injection
$$
\Hom_G((I^{\lambda,+}/\varpi_L),\pi)\incl \Hom_{KZ}(\st_{\lambda}^+/\varpi_L,\pi)=\Hom_{KZ}(\st_{\lambda}^+/\varpi_L,\pi^{K_1}),
$$
où la dernière égalité est vérifiée car  $K_1\coloneqq \ker (K\surj \GL_2(\BF_p))$ agit trivialement sur $\st_{\lambda}^+/\varpi_L$. Mais comme $\pi$ est admissible, $\pi^{K_1}$ est de dimension finie et comme $\st_{\lambda}^+/\varpi_L$ est aussi de dimension finie, on en déduit que $\Hom_G(I^{\lambda,+}/\varpi_L,\pi)$ est de dimension finie, ce qui termine la preuve. 
\end{proof}

\subsubsection{La famille universelle}\label{subsubsec:lafamun}
Soit $\Ban_{\CB}^{\zeta}G$ la catégorie abélienne des $L$-représentations de Banach unitaires de $G$, de caractère central $\zeta$, contenant un réseau dont la réduction modulo $\varpi_L$ est dans $\CB$. Notons 
$$
X\coloneqq \Spec\, R_{\CB}^{\ps,\zeta_{\lambda}}\big [\tfrac 1 p\big ] ,\quad X_{\Sp}\coloneqq \Spec\, R_{\CB}^{\lambda}.
$$
On pose $U_{\CB}^{\lambda}\coloneqq \{\CL\in L\cup \{\infty\}\mid \Pi_{\CL}^{\lambda}\in \Ban_{\CB}^{\zeta_{\lambda}}G\}$. Alors, de façon équivalente $\CL\in U_{\CB}^{\lambda}$ si et seulement si $V_{\CL}^{\lambda}$ a pour réduction\label{} $\bar \rho_{\CB}$. 

Soit $x\in X_{\Sp}$, on note $\fkm_x\subset R^{\lambda}_{\CB}$ l'idéal maximal associé et $L_x\coloneqq R^{\lambda}_{\CB}/\fkm_x$. 

\medskip
\begin{lemm}\label{lem:specrep}
Pour $n\geqslant 0$ un entier, pour tout $x\in X_{\Sp}$, il existe $\CL(x)\in U_{\CB}^{\lambda}$ tel que
$$
\begin{gathered}
\sC_{\CB}^{\lambda}\otimes_{R^{\lambda}_{\CB}}(R^{\lambda}_{\CB}/\fkm_x^{n+1})\cong \Pi_{\CL(x),n}^{\lambda},\\
\sigma_{\CB}^{\lambda}\otimes_{R^{\lambda}_{\CB}}(R^{\lambda}_{\CB}/\fkm_x^{n+1})\cong V_{\CL(x),n}^{\lambda}, \quad \rho_{\CB}^{\lambda}\otimes_{R^{\lambda}_{\CB}}(R^{\lambda}_{\CB}/\fkm_x^{n+1})\cong V_{\CL(x),n}^{\lambda}.
\end{gathered}
$$
En particulier, l'application $x\mapsto \CL(x)$ définit une bijection $X_{\Sp}(\Qpbar)\cong U_{\CB}^{\lambda}$. 
\end{lemm}
\begin{proof}

Notons 
$$
L_{x,n}\coloneqq R^{\lambda}_{\CB}/\fkm_x^{n+1},\quad \Pi_{x,n}\coloneqq\sC_{\CB}^{\lambda} \otimes_{R^{\lambda}_{\CB}}L_{x,n},\quad V_{x,n}\coloneqq \sigma_{\CB}^{\lambda}\otimes_{R^{\lambda}_{\CB}}L_{x,n}.
$$
Par construction, $\St_{\lambda}^{\lalg}$ est dense dans $\Pi_{x,n}$ et d'après \ref{lem:finitcompdef}, $\Pi_{x,n}$ est de longueur finie sur $L_x$ donc par \cite[Theorem 1.11]{pas}, le bloc de $\Ban_{\CB}^{\zeta}G$ correspondant à $\Pi_{x,n}$ contient une unique représentation de la forme $\Pi_{\CL(x)}^{\lambda}$ pour un certain $\CL(x)\in U_{\CB}^{\lambda}$, puisque cette représentation correspond à $V_{\CL(x)}^{\lambda}$ qui est absolument irréductible. Ainsi, par la proposition \ref{prop:complst}, $\Pi_{x,n}\cong \Pi_{\CL(x),n}^{\lambda}$ (l'entier $n$ provient du fait que $\Pi_{x,n}$ est annulé par $\fkm_x^{n+1})$. Comme $V_{x,n}=\bV(\Pi_{x,n})$ par définition on a $V_{x,n}=V_{\CL(x),n}^{\lambda}$ et d'après \cite[Theorem 1.11]{pas} $\bV(\Pi_{x,0})\cong \rho_{\CB}^{\lambda}\otimes_{R_{\CB,M}^{\lambda}} L_{x}$ dont on déduit le dernier isomorphisme. Finalement l'injectivité de $x\mapsto \CL(x)$ provient de l'injectivité de $\CL\mapsto V_{\CL}^{\lambda}$.

\end{proof}
On obtient alors l'équivalent de \cite[Théorème 5.4]{codonikis}.
\medskip
\begin{prop}\label{prop:dimun}
\begin{enumerate}
\item Le $R_{\CB}^{\lambda}$-module $\rho_{\CB}^{\lambda}$ est localement libre de rang $2$.
\item Le $L$-schéma $X_{\Sp}$ est lisse, réduit et purement de dimension $1$.
\item On a $\det_{R_{\CB}^{\lambda}}\rho_{\CB}^{\lambda}=R_{\CB}^{\lambda}\cdot t^{\lvert \lambda \rvert}$.
\end{enumerate}
\end{prop}
\medskip
\begin{proof}
La preuve est \emph{verbatim} la même que celle de \cite[Théorème 5.4]{codonikis} (sauf le troisième point), dont on rappelle les ingrédients principaux. On note 
$$
R=R_{\CB}^{\lambda},\quad \rho^+=\rho_{\CB}^{\lambda,+},\quad \rho=\rho_{\CB}^{\lambda}
$$
et $\wh{\rho}_x$ le complété du localisé de $\rho$ en $x\in X_{\Sp}$.

D'après le lemme \ref{lem:specrep}, $\wh{\rho}_x$ est libre de rang $2$ sur $L_x\llbracket T_x\rrbracket$ et $\End_{\sG_{\Qp}}\wh{\rho}_x\cong L_x\llbracket T_x\rrbracket$ ce qui donne une injection $R\incl \prod_x L_x\llbracket T_x\rrbracket$. Ceci montre que $R$ est réduit. 

Ensuite, comme $\rho$ est de type fini par le lemme \ref{lem:rhotf}, et $\wh{\rho}_x$ est libre de rang $2$, on obtient que $\rho$ est localement libre de rang $2$ sur $R$. Finalement, ceci implique que l'anneau local complété de $R$ en $x$ est $\wh{R}_x\cong L_x\llbracket T_x \rrbracket$, ce qui prouve que $X_{\Sp}$ est lisse purement de dimension $1$.

Pour le dernier point, notons $N=\det_{R}\rho$ qui est un $R$-module localement libre de rang $1$ et $N^+\coloneqq \det_{R^+}\rho^+$. On sait que $N^+$ est de type fini sur $R^+$ et $N^+/\varpi_L\cong \det_{k_L}\bar\rho_{\CB}$. Par le lemme de Nakayama, $N^+$ est engendré par un élément et donc $N$ est engendré par un élément. Comme $N$ est localement libre de rang $1$ et engendré par un élément il est libre de rang $1$ sur $R$. Or, $N\otimes_RL_x\cong L_x t^{\lvert \lambda \rvert}$ on en déduit que $N=R_{\CB}^{\lambda}\cdot t^{\lvert \lambda \rvert}$.

\end{proof}
\medskip

\subsubsection{L'ouvert analytique $\CU_{\CB}^{\lambda}$}\label{subsubsec:lopen}
On montre maintenant que  $U_{\CB}^{\lambda}$ est l'ensemble des points classiques d'un ouvert analytique de $\BP^{1,\an}_{L}$. 
\medskip
\begin{prop}\label{prop:uouvan}
L'application $x\mapsto \CL(x)$ est la restriction aux points classiques d'une application analytique $f\colon X_{\Sp}^{\an}\rightarrow \BP^{1,\an}_{L}$. De plus, $f$ est une immersion ouverte, c'est-à-dire que $U_{\CB}^{\lambda}$ est l'ensemble des points classiques d'un ouvert analytique $\CU_{\CB}^{\lambda}\subset \BP^{1,\an}_{L}$.
\end{prop}

On commence par un lemme plus général, qui est une conséquence de \cite{beco}.
\medskip
\begin{lemm}\label{lem:bcfam}
Soit $A^+$ une $\rmO_L$-algèbre locale complète noethérienne (topologiquement) de type fini, posons $A \coloneqq A^+[\tfrac 1 p]$ et supposons que $A$ est un anneau de Dedekind. Soit $V$ une $A$-représentation (continue) de $\sG_{\Qp}$, telle que
\begin{itemize}
\itemb $V$ est libre de rang $2$,
\itemb $\det_A V= A\cdot e_{\chi}$ où $\sG_{\Qp}$ agit sur $e_{\chi}$ par un caractère $\chi\colon \sG_{\Qp}\rightarrow L^{\times}$, en particulier, $\det_AV$ est libre de rang $1$ sur $A$,
\itemb pour tout idéal maximal $\fkm_x\subset A$, dont on note $L_x$ le corps résiduel, $V_x\coloneqq V\otimes_AL_x$ est spéciale à poids de Hodge-Tate $\lambda$ \ie $V_x\cong V_{\CL(x)}^{\lambda}$ pour $\CL(x)\in \BP^{1}(L_x)$.
\end{itemize}
Alors :
\begin{enumerate}
\item $\bD_{\dR}(V)\cong A\otimes_LD_{\dR}=Ae_0\oplus Ae_1$, en particulier, $\bD_{\dR}(V)$ est libre de rang $2$ sur $A$,
\item  $\Gr_{-\lambda_2}\bD_{\dR}(V)$ est localement libre de rang $1$ sur $A$,
\item il existe une application analytique $f\colon (\Spec A)^{\an}\rightarrow \BP^{1,\an}_L$ telle que $f(x)=\CL(x)$ sur les points fermés.
\end{enumerate}
\end{lemm}
\medskip
\begin{proof}
Notons que si les deux premiers points sont vérifiés, on définit l'application $f$ du dernier point par $(D_{\dR}\otimes_LA)\cong\bD_{\dR}(V)\rightarrow \Gr_{-\lambda_2}\bD_{\dR}(V)$.

Pour le premier point, on commence par montrer que $\bD_{\dR}(V)$ est localement libre de rang $2$ sur $A$. Pour le second point, il suffit de montrer que $\bD_{\HT}(V)$ est localement libre puisque 
$$
\bD_{\HT}(V)\cong \bigoplus_{i\in\BZ} \Gr_i\bD_{\dR}(V)=\Gr_{-\lambda_2}\bD_{\dR}(V)\oplus \Gr_{-\lambda_1}\bD_{\dR}(V).
$$
On peut écrire $A$ comme limite projective de $L$-algèbres affinoïdes $A\coloneqq \varprojlim_i A_i$. On considère $V_i\coloneqq V\otimes_AA_i$, et pour tout $i\in \BN$ on a $\bD_{\dR}(V)\coloneqq \varprojlim_i \bD_{\dR}(V_i)$ (\resp  $\bD_{\HT}(V)\coloneqq \varprojlim_i \bD_{\HT}(V_i)$) d'après \cite[Lemma 4.4]{liu} (\resp \cite[Lemma 3.5]{liu}) puisque $\bD_{\dR}(V_i)\cong \bD_{\dR}(V)\otimes_{A}A_i$ (\resp $\bD_{\HT}(V_i)\cong \bD_{\HT}(V)\otimes_{A}A_i$).

D'après \cite[Théorème 5.3.2]{beco}, $\bD_{\dR}(V_i)$ est un $A_i$-module localement libre de rang $2$ et d'après \cite[Théorème 5.1.4]{beco},  $\bD_{\HT}(V_i)$ est un $A_i$-module localement libre de rang $2$. Mais comme 
$$
\bD_{\HT}(V_i)=\Gr_{-\lambda_2}\bD_{\dR}(V_i)\oplus \Gr_{-\lambda_1}\bD_{\dR}(V_i),
$$
on en déduit que $\Gr_{-\lambda_2}\bD_{\dR}(V)$ est localement libre de rang $1$, ce qui prouve le second point.

Pour conclure la preuve du premier point, comme $A$ est un anneau de Dedekind, par la classification des modules projectifs sur les anneaux de Dedekind, il suffit de montrer que $\det_A\bD_{\dR}(V)$ est libre de rang $1$. On a
$$
{\det}_A\bD_{\dR}(V)=\bD_{\dR}({\det}_AV)=(\Bdr\wotimes_{\Qp}A\cdot e_{\chi})^{\sG_{\Qp}}\cong A\otimes_{\Qp}\bD_{\dR}(\Qp\cdot e_{\chi}),
$$
où la première égalité provient du fait que $\bD_{\dR}$ est un foncteur monoïdal. Ceci finit de prouver le premier point.

Finalement, pour $x\in \Sp A_i$ un point classique, comme $\bD_{\dR}(V_i)_x\cong \bD_{\dR}(V_x)$ d'après \cite[Théorème 5.3.2]{beco}, on a $f(x)=\CL(x)$ ce qui finit de prouver le dernier point.
\end{proof}
\medskip
\begin{rema}
Le dernier point du lemme \ref{lem:bcfam} peut s'énoncer un peu différemment. On peut construire sur $\BP^{1,\an}_L$ un $\CO_{\BP^{1,\an}_L}[\sG_{\Qp}]$-module universel $\CV^{\lambda}$ tel que $\CV_{x}^{\lambda}=V_{x}^{\lambda}$ où 
$x$ est vu comme un élément de $\BP^1(L_x)$. Alors on a $V\cong f^*\CV^{\lambda}$.
\end{rema}
\medskip
On peut maintenant conclure la preuve de la proposition \ref{prop:uouvan}.
\begin{proof}
Comme $R_{\CB}^{\lambda}=R_{\CB}^{\lambda,+}[\tfrac 1 p]$ où $R_{\CB}^{\lambda,+}$ est un anneau local complet noethérien, les hypothèses du lemme \ref{lem:bcfam} sont satisfaites par la proposition \ref{prop:dimun} donc on obtient une application $f\colon X_{\Sp}^{\an}\rightarrow \BP^1_L$ et il s'agit de montrer que c'est une immersion ouverte. On va montrer que $f$ est étale et injective sur tous les points avant de conclure en utilisant qu'un morphisme étale et injectif sur les points est une immersion ouverte (\cf. \cite[Lemme I.14]{far}).

Commençons par montrer que $f$ est étale. D'après le lemme \ref{lem:specrep}, on sait que $f$ est formellement étale aux points classique, donc $f$ est formellement étale puisque la condition d'être formellement étale est une condition ouverte. De plus, $R_{\CB}^{\lambda}$ et $\BP^{1,\an}_L$ sont de présentation finie donc $f$ est de présentation finie. Ainsi, $f$ est étale.

D'après le lemme \ref{lem:specrep}, $f$ est injective sur les points classiques, il suffit de montrer que $f$ est injective sur les $C$-points pour $C$ une extension de $\Qp$ complète et algébriquement close. Soit $x_1,x_2\coloneqq R_{\CB}^{\lambda}\rightarrow C$. Si $f(x_1)=f(x_2)$ ceci signifie que $x_1$ et $x_2$ définissent la même filtration sur $D_{\dR}\otimes_LC$, \ie la même suite exacte
$$
0\rightarrow \Gr_{-\lambda_1}(D_{\dR}\otimes_LC)\rightarrow (D_{\dR}\otimes_LC)\rightarrow \Gr_{-\lambda_2}(D_{\dR}\otimes_LC)\rightarrow 0.
$$
Par la théorie de Sen, on en déduit des isomorphismes canoniques
$$
\rho_{\CB}^{\lambda}\otimes_{R_{\CB}^{\lambda},x_1}C\cong \Gr_{-\lambda_1}(D_{\dR}\otimes_LC)(\lambda_1)\oplus\Gr_{-\lambda_2}(D_{\dR}\otimes_LC)(\lambda_2)\cong \rho_{\CB}^{\lambda}\otimes_{R_{\CB}^{\lambda},x_2}C.
$$
Comme $\rho_{\CB}^{\lambda}$ est la famille universelle sur $R_{\CB}^{\lambda}$ on en déduit que $x_1=x_2$. Ceci conclut la preuve de la proposition.
\end{proof}
\medskip
\begin{rema}
La preuve de la liberté de $\bD_{\dR}(\rho_{\CB}^{\lambda})$ utilisant le fait que $R_{\CB}^{\lambda}$ soit un anneau de Dedekind semble artificielle. Il est possible que pour $R$ un anneau de Kisin plus général (disons pour les représentations de $\sG_F$ de dimension $n$ de type cuspidal ou spécial à poids de Hodge-Tate fixés), le $\bD_{\dR}$ de la famille universelle soit libre et qu'on puisse en déduire une application de $(\Spec\,R)^{\an}$ vers une variété de drapeau (partielle) qui soit une immersion ouverte, du moins lorsque $M$ est cuspidal (\cf \cite{waer}).
\end{rema}
\medskip
\Subsection{Fin de la preuve et corollaires}\label{subsec:finsecd}
On obtient finalement le théorème suivant :
\medskip
\begin{theo}
\begin{enumerate}
\item On a $R_{\CB}^{\lambda}\cong \CO(\CU_{\CB}^{\lambda})$ qui est donc un produit fini d'anneaux principaux.
\item Le $R_{\CB}^{\lambda}$-module $\rho_{\CB}^{\lambda}$ est libre de rang $2$.
\item On a $R_{\CB}^{\lambda}=T_{\CB}^{\lambda}$, \ie $R_{\CB}^{\lambda}=\End_G\sC_{\CB}^{\lambda}$.
\item On a un isomorphisme de $R_{\CB}^{\lambda}[\sG_{\Qp}]$-modules $\sigma_{\CB}^{\lambda}\cong\rho_{\CB}^{\lambda}$.
\end{enumerate}
\end{theo}
\medskip
\begin{proof}
Dans la proposition \ref{prop:uouvan} on a montré que $X_{\Sp}^{\an}\cong \CU_{\CB}^{\lambda}$ ce qui montre le premier point et c'est donc un produit d'anneaux principaux parce que $R_{\CB}^{\lambda}$ est noethérien. Le second point est alors une conséquence de la proposition \ref{prop:dimun}.

Pour le troisième point, la preuve est exactement la même que celle de \cite[Corollaire 5.6]{codonikis}, que l'on rappelle brièvement. On a une inclusion immédiate $R_{\CB}^{\lambda}\incl \End_G\sC_{\CB}^{\lambda}$, l'autre inclusion revient à considérer $\alpha \in \End_G\sC_{\CB}^{\lambda}$ et montrer que $\bV(\alpha)\in \End_{\sG_{\Qp}}\rho_{\CB}^{\lambda}$ est une homothétie en utilisant que s'en est une modulo $\fkm_x$ pour tout $x\in X_{\Sp}$. Ceci donne $\bV(\alpha)\in R_{\CB}^{\lambda}$ et on conclut en remarquant que $\bV(\alpha-\bV(\alpha))=0$.

Par \cite[Lemme 5.12]{codonifac} le dernier point est une conséquence immédiate du premier point et du lemme \ref{lem:specrep}.
\end{proof}

\subsubsection{Universalité de $\rho_{\CB}^{\lambda}$}
En suivant mot pour mot la preuve de \cite[Théorème 5.8]{codonikis} on obtient le corollaire suivant :
\medskip
\begin{coro}
Soit $E$ une $L$-algèbre commutative de dimension $d$ et $\rho\colon \sG_{\Qp}\rightarrow \GL_2(E)$ a pour réduction $\bar \rho_{\CB}^{\oplus d}$, est de déterminant $\zeta_{\lambda}$, semi-stable à poids de Hodge-Tate $\lambda$ et $\bD_{\st}(\rho) = \Sp(\lvert \lambda \rvert)\otimes_L E$. Alors il existe un morphisme $R_{\CB}^{\lambda}\rightarrow E$ tel que $\rho=E\otimes_{R_{\CB}^{\lambda}}\rho_{\CB}^{\lambda}$.
\end{coro}
\medskip

\medskip
\begin{rema}\label{rem:degenexp}
Notons qu'à cause du cas exceptionnel, la réciproque du corollaire est fausse, \ie l'équivalent de \cite[Corollaire 5.7]{codonikis} est faux dans le cas spécial. En effet, si $\infty\in U_{\CB}^{\lambda}$ alors le quotient de $R_{\CB}^{\lambda}\rightarrow L$ correspondant à $\CL=\infty$ fournit une représentation dont le $\bD_{\st}$ n'est pas $\Sp(\lvert \lambda \rvert)$. Le corollaire reste vrai si l'on suppose  $\infty\not \in U_{\CB}^{\lambda}$, on obtient alors le corollaire suivant :
\end{rema}
\medskip

\medskip
\begin{coro}
Soit $E$ un quotient de $R_{\CB}^{\lambda}$ de degré fini sur $L$, alors $E\otimes_{R_{\CB}^{\lambda}} \rho_{\CB}^{\lambda}$ est une $E$-représentation de $\sG_{\Qp}$ de réduction $\bar \rho_{\CB}^{\oplus[E\colon L]}$, semi-stable à poids de Hodge-Tate $\lambda$. De plus, si $\infty\not \in U_{\CB}^{\lambda}$ alors $\bD_{\st}(E\otimes_{R_{\CB}^{\lambda}} \rho_{\CB}^{\lambda})\cong E\otimes_L\Sp(\lvert \lambda \rvert)$.
\end{coro}
\medskip

\section{Factorisation de la cohomologie étale du système local}
Soit $\lambda\in P_+$, $M\in \phinl$ et $\CB$ un bloc de $\Tors^{\zeta^{\lambda}_M}G$. Si $M$ est spécial supposons que $w(\lambda)>1$, on renvoie à la remarque \ref{rem:princ} pour la convention si $M$ est spécial et $w(\lambda)=1$. On commence par définir le $R_{\CB,M}^{\lambda}[G]$-module $\bPi(\rho^{\lambda}_{\CB,M})'$ comme dans \cite[Définition 5.14]{codonifac}.
\medskip
\begin{defi}\label{def:plangfam}
Soit $\rho_{\CB,M}^{\lambda,\vartriangle}$ un $R_{\CB,M}^{\lambda,+}$-réseau $\sG_{\Qp}$-stable dans le $R_{\CB,M}^{\lambda}$-dual de $\rho_{\CB,M}^{\lambda}$. Le dual de Pontryagin de $\rho_{\CB,M}^{\lambda,\vartriangle}$ est de la forme $\varinjlim_i V_i$ où les $V_i$ sont des représentations de $\sG_{\Qp}$ de rang $2$ sur des quotients de $R_{\CB,M}^{\lambda,+}$. On définit
$$
\bPi(\rho^{\lambda}_{\CB,M})'\coloneqq \varprojlim_i\bPi(V_i)'\otimes_{\rmO_L}L.
$$
\end{defi}
\medskip
Pour tout idéal maximal $\fkm_x\subset R_{\CB,M}^{\lambda}$ on a donc
$$
\Pi_x'\coloneqq\bPi(\rho_x)'\cong L_x\otimes_{R_{\CB,M}^{\lambda}}\bPi(\rho^{\lambda}_{\CB,M})'.
$$
Si $\lambda\in P$ et $M\in \phinl$ on définit
$$
\niv(M)\coloneqq \min\{n\in \BN \mid \czG^1_n\subset \Stab(\JL_M)\},
$$
 où l'on rappelle que $\czG^1_n\subset \czG_n^+$ est le noyau de la norme réduite. L'entier $\niv(M)$ est bien défini puisque $\JL_M$ est une $L$-représentation lisse de dimension finie et $\czG^1$ est un groupe compact. Pour $n\geqslant 0$ un entier on note $\phinl^n\subset \phinl$ le sous-ensemble des $M\in \phinl$ tels que $\niv(M)\leqslant n$ (notons que $\phinl^n$ est un ensemble fini). Dans toute cette section, on fixe $n\geqslant 0$ un entier, et on choisit $L$ une extension finie de $\Qp$ toutes les extensions quadratiques de $\Qp$ et telle que tous les $\JL_M$ pour $M\in \phinl^n$ soient définis sur $L$. Le but consiste à démontrer le résultat de factorisation, \ie le théorème suivant :
\medskip
\begin{theo}\label{thm:repprin}
On a un isomorphisme de $L[\BG]$-modules topologiques :
$$
\rmH^1_{\et}\Harg{\presp{\rmM}_{\Qpbar}^{n}}{\Sym \BV(1)}= \bigoplus_{\lambda \in P_+}\bigoplus_{M\in \phinl^n}\left (\wh{\bigoplus_{\CB}}\ \bPi(\rho^{\lambda}_{\CB,M})'\otimes_{R_{\CB,M}^{\lambda}} \rho_{\CB,M}^{\lambda}
\otimes_{R_{\CB,M}^{\lambda}}\cz{R}_{\CB,M}^{\lambda} \right ) \otimes_L \JL_M^{\lambda}
 $$
\end{theo}
\medskip
\begin{rema}\label{rem:princ}
\ 
\begin{itemize}
\itemb Comme annoncé précédemment, pour les facteurs associés à $\lambda\in P_+$ tels que $w(\lambda)=1$, le théorème est précisément \cite{codonifac} tordu par un caractère. Dans le cas où $M$ est spécial et $w(\lambda)=1$ (\cf quatrième point de la remarque \ref{rem:kisthm}) la convention pour $\CB$ le bloc de la steinberg est $\bPi(\rho_{\CB,M}^{\lambda})\coloneqq \St_L^{0}\otimes (\zeta_M^{\lambda}\circ \det)$ où $\St_L^{0}$ est la steinberg continue à coefficients dans $L$.
\itemb L'une des difficultés principales dans \cite{codonifac} est de démontrer des théorèmes de finitude (\cf \cite[Corollaire 0.12]{codonifac}) pour la cohomologie modulo $p$. Comme le système local modulo $\varpi_{\rmD}$ se trivialise sur le premier revêtement il n'est pas très difficile, en utilisant la suite spectrale d'Hochschild-Serre, de déduire ces théorèmes dans le cas des coefficients non triviaux (\cf \num \ref{subsubsec:finit}).
\end{itemize}
\end{rema}
\medskip
Soit $\lambda\in P_+$ un poids que l'on fixe dans toute la suite. À fin d'alléger les notations, posons 
$$
\rmH^{\lambda}_{\Qpbar}\coloneqq (\varprojlim_m\varinjlim_{\substack{K/\Qp\\ [K:\Qp]<\infty}}\rmH^1_{\et}\Harg{\presp{\rmM}_K^{n}}{\BV^+_{\lambda}(1)/p^m})\otimes_{\rmO_L}L.
$$
Par définition et d'après \cite{van}, on a
$$
\rmH^1_{\et}\Harg{\presp{\rmM}_{\Qpbar}^{n}}{\Sym \BV(1)}=\bigoplus_{\lambda\in P_+}\rmH^{\lambda}_{\Qpbar},
$$
$$
\rmH^1_{\et}\Harg{\presp{\rmM}_{\Qpbar}^{n}}{\Sym \BV(1)}=\bigoplus_{\lambda\in P_+}\bigoplus_{M\in \phinl^n}\rmH^{\lambda}_{\Qpbar}[M]\otimes_{L}\JL_M^{\lambda}.
$$
Il suffit donc de montrer que pour $M\in \phinl$ on a un $L[G\times \sG_{\Qp}]$-isomorphisme topologique
\begin{equation}\label{eq:toprov}
\rmH^{\lambda}_{\Qpbar}[M]\cong \wh{\bigoplus_{\CB}} (\bPi(\rho^{\lambda}_{\CB,M})'\otimes_{R_{\CB,M}^{\lambda}} \rho_{\CB,M}^{\lambda}
\otimes_{R_{\CB,M}^{\lambda}}\cz{R}_{\CB,M}^{\lambda} )
\end{equation}
En particulier, pour $M\cong \Sp(\lvert \lambda \rvert-2)$, on obtient :
\medskip
\begin{coro}
Soit $\lambda\in P_+$ tel que $w(\lambda)>1$. Alors
$$
 \rmH^1_{\et}\Harg{\BH_{\Qpbar}}{\BV_{\lambda}(1)}=\wh{\bigoplus_{\CB}}\ \bPi(\rho^{\lambda}_{\CB,\Sp})\otimes_{R_{\CB,\Sp}^{\lambda}} \rho_{\CB,\Sp}^{\lambda}\otimes_{R_{\CB,\Sp}^{\lambda}}\cz{R}_{\CB,\Sp}^{\lambda}.
$$
\end{coro}
\medskip
Comme (\ref{eq:toprov}) est démontré si $w(\lambda)=1$, quitte à tordre par un caractère, on suppose que $w(\lambda)>1$ dans toute la suite.

\subsection{Rappels et compléments de \cite{van}}
On définit 
$$
\rmH^{\lambda,+}_{C}\coloneqq\rmH^1_{\et}\Harg{\presp{\rmM}_C^{n}}{\BV_{\lambda}^+(1)}/{\rm torsion}, \quad \rmH^{\lambda}_{C}\coloneqq \rmH^{\lambda,+}_{C}\otimes_{\rmO_L}L,
$$
où l'on sait que la torsion est finie (\cf \cite[Proposition 10.9]{van}).
\subsubsection{Entrelacement}
On a les théorèmes suivants démontrés dans \cite{van} :
\medskip
\begin{theo}\label{thm:oldvan}
Soit $\Pi$ une $L$-représentation unitaire et absolument irréductible de $G$ sur un espace de Banach, 
$$
\Hom_{G}\intnn{\Pi'}{\rmH^{\lambda}_{C}} \cong 
\begin{cases}
V_{M,\sL}^{\lambda}\otimes_L \JL_M^{\lambda} & \text{ si } \Pi=\Pi_{M,\sL}^{\lambda},\quad M\in \phinl^n\\
0 & \text{ si $\Pi$ n'est pas du type ci-dessus}
\end{cases}
$$
\end{theo}
\medskip

\Subsubsection{Une inclusion}
Pour ce numéro, on doit rentrer un peux dans le détail de \cite{van}. Notons simplement $X\coloneqq \presp{\rmM}_{\Qp}^n$ et posons $\CO^n\coloneqq \CO_X(X)\otimes_{\Qp}L$ laquelle est une $L$-représentation de $G\times \czG$. Pour $k,m\in \BZ$ on définit $\CO^n\{k,m\}$ par $\CO^n\{k,m\}=\CO^n$, en tant que $L[\czG]$-module topologique, et l'action de $g\in G$  est définie pour $f\in\CO^n$ par :
$$
g = 
\begin{pmatrix}
 a & b \\ c & d
\end{pmatrix}\in G,\quad 
(g\star f)(z) = j(z,g)^{-k}\det(g)^{m}f(g\cdot z),
$$
où l'on note $\pi \colon \presp{\rmM_{\Qp}^n}\rightarrow \presp{\rmM_{\Qp}^0}$ la projection naturelle et $j(z,g) = (a-c\pi(z))\in\CO$ le facteur d'automorphie. L'action considérée a un sens puisque $\CO^n$ est naturellement un $\CO^0$-module. Pour $\lambda\in P_+$, on pose
$$
\CO_{\lambda}^n\coloneqq \CO^n\{1-w(\lambda),1-\lambda_2\}\otimes_L\lvert \det \rvert_p^{\tfrac{1-\lvert \lambda \rvert}{2}}.
$$
\medskip
\begin{lemm}\label{lem:nogbdedfun}
Soit $\lambda\in P_+$ tel que $w(\lambda)>1$, le $L[G]$-module $\CO_{\lambda}^n$ n'a pas de vecteurs $G$-bornés, \ie $(\CO_{\lambda}^n)^{G-\rmb}=0$.
\end{lemm}
\medskip
\begin{proof}
Par \cite[Remarque A.2]{codoni}, pour $w(\lambda)=1$, on a $(\CO_{\lambda}^n)^{G-\rmb}=L$. On peut supposer que $\lambda=(0,k+1)$ avec $k\geqslant 1$, et quitte à tordre par un caractère. Soit $f\in \CO_{\lambda}^n$ une fonction $G$-bornée, en agissant par le sous-groupe unipotent de $G$, on en déduit en particulier que 
$$
\{f(z+b)\}_{b\in \Qp}\subset \CO_{\lambda}^n
$$
est une partie bornée. Mais ceci implique que $f$ est une fonction bornée sur $X$ et donc, par \cite[Remarque A.2]{codoni}, $f\in L$ est une constante. Mais dans ce cas, par l'action de $G$
$$
\{\lvert ad\rvert^{\tfrac{1-k}{2}}a^{-k}f\}_{a,d\in \Qpt}\subset L,
$$
est une partie bornée. Donc $f=0$, ce qui conclut la preuve.
\end{proof}
Soit $M\in \phinl$, on introduit quelques notations : 
\begin{itemize}
\itemb $\rmH^{\lambda}_{\pet,C}\coloneqq \rmH^1_{\pet}\Harg{\presp{\rmM}_C^n}{\BV_{\lambda}(1)}$, le premier groupe de cohomologie proétale de $\BV_{\lambda}(1)$,
\itemb pour $Z$ un $L[\czG]$-module, $Z[M]\coloneqq \Hom_{\czG}(\JL_M\otimes_L\lvert \nrd\rvert_p^{\tfrac{1-w}{2}}, Z)$,
\itemb $X_{\st}^1(M)\coloneqq (\Bstp\otimes_{\Qp^{\nr}}M)^{N=0,\varphi=p}$.
\end{itemize}
\medskip
\begin{prop}\label{prop:injco}
Soit $\lambda\in P_+$ tel que $w(\lambda)>1$ et $M\in\phinl^n$. On a une inclusion de $L[\sG_{\Qp}\times G]$-modules topologiques
$$
\rmH^{\lambda}_{C}[M]\incl t^{\lambda_1}X_{\st}^1(M[\lambda_1])\wotimes_L(\wh{\LL}_M^{\lambda})'
$$
\end{prop}
\medskip
\begin{proof}
On peut supposer $\lambda_1=0$ quitte à tordre par un caractère et poser $w=\lambda_2$. On commence par supposer $M$ cuspidal. Alors, d'après \cite[Corollaire 11.11]{van} on a une suite exacte
$$
0\rightarrow \rmB_{\lambda}\wotimes_{\Qp}\CO^n_{\lambda}[M]\rightarrow \rmH_{\pet,C}^{\lambda}[M]\rightarrow X_{\st}^1(M)\wotimes_L(\LL_M^{\lambda})'\rightarrow 0,
$$
où $\rmB_{\lambda}\coloneqq \Bdrp/t^{w}$.
On passe aux vecteurs $G$-bornés dans la suite exacte et comme
\begin{itemize}
\itemb $(\CO_{\lambda}^n[M])^{G-\rmb}=0$ par le lemme \ref{lem:nogbdedfun},
\itemb $ (\rmH_{\pet,C}^{\lambda}[M])^{G-\rmb}=\rmH_{C}^{\lambda}[M]$ par \cite[Théorème 10.1]{van},
\itemb $((\LL_M^{\lambda})')^{G-\rmb}\cong (\wh{\LL}_M^{\lambda})'$ par \cite[Lemme 5.3 b)]{codonifac},
\end{itemize}
on obtient l'inclusion voulue. 

On suppose maintenant que $M$ est spécial et, quitte à tordre par un caractère, on peut supposer $M=\Sp(\lvert \lambda \rvert -2)$. Alors, d'après \cite[Corollaire 12.15]{van} on a une suite exacte 
$$
0\rightarrow Q_{\lambda}^0\otimes_L W_{\lambda}\rightarrow  \rmH^{\lambda}_{\pet}[M]\rightarrow \rmB_{\lambda}\wotimes_{\Qp}\intn{{\St}^{\lan}_{\lambda}}'\rightarrow Q_{\lambda}^1\wotimes_L \intn{{\St}^{\lalg}_{\lambda}}'\rightarrow 0.
$$
où
$$
\rmU_{\lambda}^0(L)\coloneqq L\otimes_{\Qpd}\intn{\Bcrisp}^{\varphi^2=p^{w+1}},\quad \rmU_{\lambda}^1\coloneqq L\otimes_{\Qpd}\intn{\Bcrisp}^{\varphi^2=p^{w-1}},\quad Q_{\lambda}^i\coloneqq \frac{\rmB_{\lambda}\otimes_{\Qp}L}{\rmU_{\lambda}^i(L)}
$$
Notons $K\coloneqq \ker (\rmB_{\lambda}\wotimes_{\Qp}\intn{{\St}^{\lan}_{\lambda}}'\rightarrow Q_{\lambda}^1\wotimes_L \intn{{\St}^{\lalg}_{\lambda}}')$. Alors on obtient la suite exacte courte
$$
0\rightarrow Q_{\lambda}^0\otimes_L W_{\lambda}\rightarrow  \rmH^{\lambda}_{\pet}[M]\rightarrow K \rightarrow 0,
$$
où l'on passe aux vecteurs $G$-bornés. Or comme
\begin{itemize}
\itemb $W_{\lambda}^{G-\rmb}=0$,
\itemb $ (\rmH_{\pet,C}^{\lambda}[M])^{G-\rmb}=\rmH_{C}^{\lambda}[M]$ par \cite[Théorème 10.1]{van} comme précédemment,
\itemb $K^{G-\rmb}\incl \rmU_{\lambda}^1(L)\wotimes_L \sD^{\lambda}$ par le lemme \ref{lem:gbded}, en prenant les vecteurs $G$-bornés dans la définition de $K$,
\end{itemize}
on obtient une injection $\rmH^{\lambda}_{\pet}[M]\incl \rmU_{\lambda}^1(L)\wotimes_L \sD^{\lambda}$ et comme $\rmU_{\lambda}^1(L)\subset X_{\st}^1(M)$ on obtient l'inclusion voulue.
\end{proof}
\subsubsection{Résultats de finitude}\label{subsubsec:finit}
D'après \cite[Proposition 10.12]{van}, on a la proposition suivante : 
\medskip
\begin{prop}\label{prop:finit}
Soit $K$ une extension finie de $\Qp$, alors $\rmH_{\et}^1\Harg{\presp{\rmM}_K^n}{\BV_{\lambda}^+/\varpi_D}$ est le dual d'un $\rmO_L[G]$-module lisse de longueur finie.
\end{prop}
\medskip
Une preuve similaire par descente galoisienne donne la proposition suivante : 
\medskip
\begin{prop}\label{prop:finidim}
Soit $\pi$ une $k_L$-représentation lisse et admissible de $G$, alors
$$
\Hom_{k_L[G]}\Harg{\pi'}{\rmH^1_{\et}\Harg{\presp{\rmM}_C^n}{\BV_{\lambda}^+/\varpi_{\rmD}}}
$$
est un $k_L$-espace vectoriel de dimension finie.
\end{prop}
\medskip
\begin{proof}
Le corollaire \cite[Corollaire 4.21]{codonifac} implique que pour tout $k_L$-système local (de rang fini) trivial $T$,  $\Hom_{k_L[G]}\Harg{\pi'}{\rmH^1_{\et}\Harg{\presp{\rmM}_C^n}{T}}$ est de dimension finie et comme $\BV_{\lambda}^+/\varpi_{\rmD}$ est trivial sur $\presp{\rmM}_C^n$ dès que $n\geqslant 1$, il suffit de montrer le résultat pour $n=0$ ce que l'on fait par la suite spectrale de Hochschild-Serre appliquée au cas $n=1$. Fixons $\pi$ un $k_L[G]$-module lisse et pour $X$ un $k_L[G]$-module, notons provisoirement $X[\pi']\coloneqq \Hom_{k_L[G]}\Harg{\pi'}{X}$. Notons
$$
\rmH\coloneqq \rmH^1_{\et}\Harg{\presp{\rmM}_C^0}{\BV_{\lambda}^+(1)/\varpi_{\rmD}},\quad \wt{\rmH}\coloneqq \rmH^1_{\et}\Harg{\presp{\rmM}_C^1}{\BV_{\lambda}^+(1)/\varpi_{\rmD}}.
$$
Comme on l'a rappelé, $\wt{\rmH}[\pi']$ est de dimension finie. Comme $\presp{\rmM}^1_C\rightarrow \presp{\rmM}^0_C$ est un recouvrement étale de groupe de Galois $\BF\coloneqq \BF_{q^2}^{\times}$, la suite spectrale de Hochschild-Serre donne une suite exacte
$$
0\rightarrow \rmH^1(\BF,\wt{\rmH})\rightarrow \rmH\rightarrow \wt{\rmH}^\BF
$$
Premièrement, justifions que $\rmH^1(\BF,\wt{\rmH})[\pi']=0$. Comme l'action de $G$ et $\czG$ commutent, on a $\rmH^1(\BF,\wt{\rmH})[\pi']=\rmH^1(\BF,\wt{\rmH}[\pi'])$. Mais comme $\wt{\rmH}[\pi']$ est un $k_L$-espace vectoriel de dimension finie et que $\BF$ est un groupe fini d'ordre premier à $p$,  $\rmH^1(\BF,\wt{\rmH}[\pi'])=0$. On obtient une inclusion
$$
\rmH[\pi']\incl \wt{\rmH}^\BF[\pi']
$$
Or, $\wt{\rmH}^\BF$ est un sous-$k_L[G]$-module de $\wt{\rmH}$ et donc on obtient une application injective 
$$
\rmH[\pi']\incl\wt{\rmH}[\pi'].
$$
Comme $\wt{\rmH}[\pi']$ est de dimension finie, ceci termine la preuve.
\end{proof}
\subsubsection{Vecteurs presque lisses}\label{subsubsec:vpl}
Soit $V^+$ une $\rmO_L$-représentation de $\sG_{\Qp}$. Rappelons qu'on dit que $x\in V^+$ est \emph{presque lisse} si pour tout entier $k\geqslant 1$, la classe de $x$ modulo $p^k$ est fixe par un sous-groupe ouvert de $\sG_{\Qp}$. Les vecteurs presque lisses de $V^+$ forment une sous-$\rmO_L$-représentation de $V^+$.
\medskip
\begin{prop}\label{prop:gqpplis}
L'application naturelle
$$
\rmH^{\lambda,+}_{\Qpbar}\rightarrow \rmH^{\lambda,+}_{C}
$$
est injective et identifie $\rmH^{\lambda,+}_{\Qpbar}$ aux vecteurs presque lisses de $\rmH^{\lambda,+}_{C}$ sous l'action de $\sG_{\Qp}$. 
\end{prop}
\medskip
\begin{proof}
Soit $k\geqslant 1$ un entier, et $K/\Qp$ une extension finie. Comme on a supposé que $w(\lambda)>1$ on a $\rmH_{\et}^0(\presp{\rmM}_C^n,\BV_{\lambda}^+/\varpi_D^k)=0$ dès que $k>2n$. Donc la suite spectrale d'Hochschild-Serre donne
$$
\rmH_{\et}^1(\presp{\rmM}_K^n,\BV_{\lambda}^+/\varpi_D^k)\cong\rmH_{\et}^1(\presp{\rmM}_C^n,\BV_{\lambda}^+/\varpi_D^k)^{\sG_{K}}
$$
En passant à la limite projective sur $k$, on obtient le résultat.
\end{proof}
\medskip
En particulier, comme les $V_{M,\sL}^{\lambda}$ ont des réseau $\sG_{\Qp}$-stable dont tous les vecteurs sont presque lisses on déduit du théorème \ref{thm:oldvan} le corollaire suivant :
\medskip
\begin{coro}\label{coro:entr}
Soit $\Pi$ une $L$-représentation unitaire et irréductible de $G$ sur un espace de Banach, alors
$$
\Hom_{G}\intnn{\Pi'}{\rmH^{\lambda}_{{\Qpbar}}} \cong 
\begin{cases}
V_{M,\sL}^{\lambda}\otimes_L \JL_M^{\lambda} & \text{ si } \Pi=\Pi_{M,\sL}^{\lambda} \text{ et }\niv(M)\leqslant n\\
0 & \text{ si $\Pi$ n'est pas du type ci-dessus}
\end{cases}
$$
\end{coro}
\medskip
\begin{coro}\label{cor:propveclis}
\ 
\begin{enumerate}
\item $\rmH^{\lambda,+}_{\Qpbar}$ est sans $p$-torsion, $p$-adiquement complet et $p$-saturé dans $\rmH^{\lambda,+}_{C}$.
\item Pour tout entier $k\geqslant 1$, $\rmH^{\lambda,+}_{\Qpbar} /p^k$ est une $\rmO_L$-représentation lisse de $\sG_{\Qp}$.
\item Pour tout entier $k\geqslant 1$ et toute extension finie $K/\Qp$, le $\rmO_L[G]$-module $(\rmH^{\lambda,+}_{\Qpbar} /p^k)^{\sG_K}$ est lisse de longueur finie.
\end{enumerate}
\end{coro}
\medskip
\begin{proof}
Les deux premiers points sont des conséquences de \ref{prop:gqpplis} et du fait que $\rmH^{\lambda,+}_{C}$ est $p$-adiquement complet et sans $p$-torsion.

Pour le dernier point, on se ramène au cas $k=1$ par un dévissage immédiat. La preuve se fait comme \cite[Proposition 5.19]{codonifac}, on la rappelle brièvement. Notons
$$
Y\coloneqq (\rmH^{\lambda,+}_{\Qpbar} /p)^{\sG_K},\ Z\coloneqq\rmH_{\et}^1\Harg{\presp{\rmM}_C^n}{\BV_{\lambda}^+(1)/\varpi_D}^{\sG_K}
$$
On a une injection continue $\iota\colon Y\incl Z$, $G$-équivariante, définie comme composée des morphismes$Y\rightarrow (\rmH^{\lambda,+}_{C} /p)^{\sG_K}\rightarrow Z$ donnée par les points précédents. Par la proposition \ref{prop:finit}, il suffit de montrer que $\iota$ est un homéomorphisme sur son image. Or, $Z$ étant profini, $Y$ est séparé. Soit $H\subset G$ un sous-groupe ouvert compact. Alors $Y$ et $Z$ sont naturellement des $\rmO_L\llbracket H\rrbracket$-modules et $\iota$ est $\rmO_L\llbracket H\rrbracket$-linéaire. De plus, $Z$ est un $\rmO_L\llbracket H\rrbracket$-module de type fini et donc $Y$ aussi (comme $\rmO_L\llbracket H\rrbracket$ est noethérien). Or, $Y$ est séparé, sa topologie est définie par sa topologie naturelle de $\rmO_L\llbracket H\rrbracket$-module de type fini. En particulier, $Y$ est profini et $\iota$ est un isomorphisme sur son image.
\end{proof}

\Subsection{Décomposition en blocs}
\subsubsection{Rappels sur la théorie des blocs}\label{subsubsec:blocth}
Soit $\zeta\colon \Qpt\rightarrow \rmO_L^{\times}$ un caractère unitaire. On rappelle que l'on note $\Tors\, G$ (\resp $\Tors^{\zeta}\,G$) la catégorie abélienne des $\rmO_L[G]$-modules lisses de longueur finie (et de caractère central $\zeta$). Rappelons que pour $\pi$ un $\rmO_L[G]$-module lisse, on définit son dual de Pontryagin par
$$
\pi^{\vee}\coloneqq \Hom_{\rmO_L}(\pi,L/\rmO_L),
$$
où l'on considère les applications continues et $\pi^{\vee}$ est muni de l'action contragrédiente de $G$.

Pour $\pi\in \Tors^{(\zeta)}\,G$ on note $P_{\pi}^{(\zeta)}\surj \pi^{\vee}$ l'enveloppe projective de $\pi^{\vee}$ (\cf \cite[\S II.5]{gab} et \cite[\S 2]{pas}) ; c'est un $\rmO_L[G]$-module compact. Pour $\CB$ un bloc de $\Tors^{(\zeta)}\,G$ posons
$$
P_{\CB}^{(\zeta)}\coloneqq \bigoplus_{\pi\in \CB}P_{\pi}^{(\zeta)},\quad E_{\CB}^{(\zeta)}\coloneqq \End_G(P_{\CB}^{(\zeta)}),\quad Z_{\CB}^{(\zeta)}\coloneqq Z(E_{\CB}^{(\zeta)}),
$$
où $Z(E_{\CB}^{(\zeta)})$ désigne le centre de $E_{\CB}^{(\zeta)}$.

\medskip
\begin{prop}\label{prop:blocs}
On a une décomposition 
$$
\Tors^{(\zeta)}\, G=\prod_{\CB}\Tors_\CB^{(\zeta)}\, G.
$$
De plus, pour $\Pi\in \Tors^{(\zeta)}\, G$, on a 
$$
\Pi\cong \bigoplus_{\CB}(P_{\CB}^{(\zeta)}\otimes_{E_{\CB}^{(\zeta)}}\bm_{\CB}(\Pi^{\vee}))^{\vee},
$$
où $\bm_{\CB}^{(\zeta)}(\Pi^{\vee})\coloneqq \Hom_G(P_{\CB}^{(\zeta)},\Pi^{\vee})$.
\end{prop}
\medskip
\medskip
\begin{lemm}\label{lem:quotp}
Soit $\fkm_x\subset R_{\CB,M}^{\lambda}[1/p]$ un idéal maximal ; notons $L_x$ son corps résiduel, \hbox{$\zeta\coloneqq \zeta_M^{\lambda}$} et rappelons que $\Pi_x\coloneqq \bPi(\rho_x)$ est la $L_x$-représentation de Banach unitaire associée à la spécialisation de la représentation universelle au point $x$. Alors,
$$
\Pi_x'\cong P_{\CB}\otimes_{E_{\CB}}\bm_{\CB}(\Pi_x')
$$
où l'on rappelle que $\bm_{\CB}(\Pi_x')\cong \Hom_G(P_{\CB},\Pi_x')$.
\end{lemm}
\medskip
\begin{proof}
Notons 
$$
R^+\coloneqq R_{\CB}^{\ps,\zeta}\subset R\coloneqq R_{\CB,M}^{\lambda},\quad L\coloneqq L_x,\quad E\coloneqq E_{\CB}^{\zeta},\quad E_x\coloneqq E\otimes_{R^+}L
$$

Soit $\Pi_x^+\subset \Pi_x$ un $\rmO_{L_x}$-réseau stable et pour $k\geqslant 0$ un entier, $\pi_k\coloneqq \Pi_x^+/\varpi_L^k$. Alors $\pi_k\in \Tors_\CB^{\zeta}\, G$ et d'après la proposition \ref{prop:blocs}, on a 
$$
\pi_k\cong (P_{\CB}\otimes_{E}\bm_{\CB}(\pi_k^{\vee}))^{\vee}
$$
En passant à la limite projective sur $k$, comme 
$$
\varprojlim_k\pi_k^{\vee}\cong\varprojlim_k \Hom_{\rmO_L}(\pi,\varpi_L^{-k}\rmO_L/\rmO_L)\cong\Hom_{\rmO_L}(\varprojlim_k\pi_k,\rmO_L) , 
$$
ce qui vaut pour toute famille projective $\pi_k$ de $\rmO_L[G]$-modules tels que $\pi_{k+1}/\varpi_L^{k}=\pi_{k}$, on obtient, en prenant le $\rmO_L$-dual continu
$$
(\Pi_x^+)'\cong P_{\CB}\otimes_E\bm_{\CB}((\Pi_x^+)')
$$
on inverse $p$ pour obtenir le résultat.
\end{proof}

\subsubsection{Décomposition au niveau entier}\label{subsubsec:decentier}
On commence par démontrer une décomposition au niveau entier, posons $\bm_{\CB,\et}^+\coloneqq \Hom_G(P_{\CB},\rmH^{\lambda,+}_{\Qpbar})$. 

\medskip
\begin{prop}
On a un isomorphisme naturel de $L[\BG]$-modules topologiques
$$
\rmH^{\lambda,+}_{\Qpbar}\cong\wh{\bigoplus_{\CB}}\left [ P_{\CB}\otimes_{E_{\CB}}\bm_{\CB,\et}^+\right ].
$$
où le complété à droite est le complété $p$-adique.
\end{prop}
\medskip
\begin{proof}
La preuve est strictement la même que celle de \cite[Proposition 5.21]{codonifac} ; rappelons-en les points principaux. Notons $X\coloneqq \rmH^{\lambda,+}_{\Qpbar}$ et pour $K/\Qp$ une extension finie  $X_{k,K}\coloneqq (X/p^k)^{\sG_{K}}$. D'après la proposition \ref{prop:finit}, $X_{k,K}$ est le dual d'un $\rmO_L[G]$-module lisse de longueur finie. On peut donc appliquer la décomposition en blocs de la proposition \ref{prop:blocs} au dual de $X_{k,K}$ \ie
$$
X_{k,K}=\bigoplus_{\CB}(P_{\CB}\otimes_{E_{\CB}}\bm_{\CB}(X_{k,K}).
$$
En passant à la limite inductive sur les extensions finies $K/\Qp$, puisque $\bm_{\CB}(X/p^k)=\varinjlim_K\bm_{\CB}(X_{k,K})$ car $P_{\CB}$ est compact et les $X_{k,K}$ sont profinis, on obtient, par le corollaire \ref{cor:propveclis}
$$
X/p^k=\bigoplus_{\CB}(P_{\CB}\otimes_{E_{\CB}}\bm_{\CB}(X/p^k)).
$$
Comme $X=\varprojlim_kX/p^k$ d'après le corollaire \ref{cor:propveclis} il reste à montrer $\bm_{\CB}(X/p^k)\cong \bm_{\CB}(X)/p^k$. La suite de la preuve de \cite[Proposition 5.21]{codonifac} démontre essentiellement le lemme suivant
\medskip
\begin{lemm}
Soit $X$ un $\rmO_L[G]$-module topologique tel que
\begin{itemize}
\itemb $X$ est sans $p$-torsion et $p$-adiquement complet,
\itemb $X/p$ est la réunion de $\rmO_L$-modules de longueur finie.
\end{itemize}
Alors, pour tout $k\geqslant 1$ on a un isomorphisme $\bm_{\CB}(X/p^k)\cong \bm_{\CB}(X)/p^k$.
\end{lemm}
\medskip
Dans notre cas ce lemme s'applique puisque le premier point est vérifié par le corollaire \ref{cor:propveclis} et le second point par la proposition \ref{prop:finit}.
\end{proof}

Posons 
$$
\check{\bm}^+_{\CB,\et}\coloneqq \Hom_{\rmO_L}(\bm_{\CB,\et}^+,\rmO_L),
$$
le dual continu topologique de $\bm_{\CB,\et}^+$. C'est un $\rmO_L$-module compact sans $p$-torsion muni d'une action de $Z_{\CB}$.
\medskip
\begin{prop}\label{prop:checktypfin}
 Le $Z_{\CB}$-module $\check{\bm}_{\CB,\et}^+$ est de type fini.
\end{prop}
\medskip
\begin{proof}
La preuve se fait exactement comme pour le \cite[Théorème 5.22]{codonifac} en utilisant la proposition \ref{prop:finidim}, rappelons la brièvement. Soit $\fkm\subset Z_{\CB}$ l'idéal maximal. Par le lemme de Nakayama (topologique) et la dualité, il suffit de montrer que $(\bm_{\CB,\et}^+/\varpi_L)[\fkm]$ est de dimension finie sur $k_L$. On a une injection
$$
\bm_{\CB,\et}^+/\varpi_L[\fkm]\incl \Hom_G(P_{\CB}\otimes_{Z_{\CB}}k_L,\rmH^1_{\et}\Harg{\presp{\rmM}_C^n}{\BV_{\lambda}^+(1)/\varpi_{\rmD}}).
$$
Mais comme, d'après \cite[corollaire 6.7]{patu}, $P_{\CB}\otimes_{Z_{\CB}}k_L$ est le dual d'une $L$-représentation lisse et admissible de $G$, on conclut la preuve de la proposition par la proposition \ref{prop:finidim}.
\end{proof}

\subsubsection{$E_{\CB,M}^{\lambda}$ est une algèbre d'Azumaya}\label{subsubsec:azu}
Posons 
$$
E_{\CB,M}^{\lambda}\coloneqq E_{\CB}^{\zeta_M^{\lambda}}\otimes_{R_{\CB}^{\ps,\zeta_M^{\lambda}}}R_{\CB,M}^{\lambda}.
$$
Dans ce numéro, on démontre que $E_{\CB,M}^{\lambda}$ est une $R_{\CB,M}^{\lambda}$-algèbre d'Azumaya. Notons que $E_{\CB}^{\zeta_M^{\lambda}}$ n'est pas une algèbre d'Azumaya et que le résultat tient car les représentations apparaissant dans l'anneau de Kisin sont absolument irréductibles. Pour démontrer ce résultat on doit raisonner bloc par bloc. Rappelons que (sous l'hypothèse que $p>3$), les blocs absolus $\CB$ de $\Tors\, G$ sont de la forme
\begin{enumerate}
\item[{\it(i)}]  $\CB=\{\pi\}$ où $\pi$ est une $k_L$-représentation supersingulière de $G$,
\item[{\it(ii)}] $\CB=\{\Ind_B^G\chi_1\otimes \chi_2\omega^{-1},\Ind_B^G\chi_2\otimes \chi_1\omega^{-1}\}$ où $\chi_1,\chi_2\colon \Qpt\rightarrow k_L^{\times}$ sont des caractères lisses tels que $\chi_1\chi_2^{-1}\neq 1,\omega^{\pm 1}$,
\item[{\it(iii)}] $\CB=\{\Ind_B^G\chi\otimes \chi\omega^{-1}\}$ où $\chi\colon \Qpt\rightarrow k_L^{\times}$ est un caractère lisse,
\item[{\it(iv)}] $\CB=\{\eta\circ \det, \Sp\otimes (\eta\circ \det), (\Ind_B^G\omega\otimes \omega^{-1})\otimes(\eta\circ \det)\}$ où $\eta\colon \Qpt \rightarrow k_L^{\times}$ est un caractère lisse et $\Sp\coloneqq \St_{\rmO_L}\otimes_{\rmO_L }k_L$ est la Steinberg lisse à coefficients dans $k_L$.
\end{enumerate}
\medskip
\begin{prop}
Soit $\fkm_x\subset R_{\CB,M}^{\lambda}$ un idéal maximal dont on note $L_x$ le corps résiduel. Posons $E_x\coloneqq E_{\CB,M}^{\lambda}\otimes_{R_{\CB,M}^{\lambda}}L_x$, alors si $\CB$ est de type
\begin{enumerate}
\item[{\it(i)}]  $E_x\cong L_x$,
\item[{\it(ii)}] $E_x\cong M_2(L_x)$,
\item[{\it(iii)}] $E_x\cong M_2(L_x)$,
\item[{\it(iv)}] $E_x\cong M_3(L_x)$.
\end{enumerate}
En particulier, $E_{\CB,M}^{\lambda}$ est une $R_{\CB,M}^{\lambda}$-algèbre d'Azumaya.
\end{prop}
\medskip
\begin{proof}
Comme annoncé on raisonne bloc par bloc. Soit
$$
E\coloneqq E_{\CB}^{\zeta^{\lambda}_M},\quad R^+\coloneqq R_{\CB}^{\zeta_M^{\lambda}},\quad R_M\coloneqq R_{\CB,M}^{\lambda}
$$
On note $\fkr\coloneqq R^+\cap \bigcap_{x}\fkm_x$ \emph{l'idéal de réductibilité} où l'intersection est prise sur tous les idéaux maximaux de $R^+[1/p]$ tels que la représentation $\rho_x$ soit réductible. Comme pour $\fkm_x\subset R_M$, $\rho_x$ est irréductible, l'image de tout élément non nul de $\fkr$ dans $L_x$ n'est pas nulle.

Si $\CB$ est de type {\it (i)} alors par \cite[Theorem 1.5]{pas}, $E\cong R^+$ donc $E_x\cong L_x$.

Si $\CB$ est de type {\it (ii)} alors par \cite[Corollary B.20]{pas}, $\fkr=(c)$ est principal et par \cite[Corollary B.27]{pas}, $E_{c}\cong M_2(R^+_{c})$ et donc $E_x\cong M_2(L_x)$.

Si $\CB$ est de type {\it (iii)} c'est exactement l'énoncé de \cite[Corollary 9.28]{pas} (combiné a \cite[Lemme 9.21]{pas} donnant la condition sur \og $(uv-vu)(uv-vu)^*$\fg puisque $\rho_x$ est absolument irréductible.

Si $\CB$ est de type {\it (iv)}, l'énoncé n'est pas n'apparaît pas tel quel dans \cite[\S 10]{pas} mais on peut le déduire des résultats de cette section. 

Quitte à tordre par un caractère on peut supposer que $\eta=1$, posons
$$
\pi_1\coloneqq \Ind_B^G\omega\otimes \omega^{-1},\quad \pi_2\coloneqq \Sp,\quad \pi_3\coloneqq \boldsymbol{1}_G.
$$
Pour $i=1,2,3$ on pose $P^i$ l'enveloppe projective de $\pi_i^{\vee}$, $P_x^i\coloneqq P^i\otimes_{R^+}L_x$ et pour $j=1,2,3$, $E^{i,j}\coloneqq \Hom_G(P^i,P^j)$, $E_x^{i,j}\coloneqq E^{i,j}\otimes_{R^+}L_x$(on change légèrement les notations par rapport à \cite[\S 10]{pas}). Comme les $P^i$ sont compacts et projectifs on a $E_x^{i,j}=\Hom_G(P^i_x,P^j_x)$ et pour $\varphi\in E^{i,j}$ on note $\varphi_x\in E^{i,j}_x$ l'élément associé. Alors, $E=(E^{i,j})_{1\leqslant i,j\leqslant 3}$ et $E_x =(E^{i,j}_x)_{1\leqslant i,j\leqslant 3}$. On veut donc montrer que pour tous les entiers $i,j=1,2,3$, $E^{i,j}_x\cong L_x$. 

Pour $i=j$ le \cite[Corollary 10.78]{pas} donne $E^{ii}\cong R^+$ donc $E^{ii}_x\cong L_x$ pour $i=1,2,3$. On a des idéaux $\fka^{ii}\subset E^{ii}$ (\cf \cite[Corollary 10.79]{pas} et avant le Lemme 10.74 pour la définition) tels que $\fka^{ii}\cong \fkr$ d'après \cite[Lemma 10.80]{pas} par l'isomorphisme $E^{ii}\cong R^+$ ci-dessus. Ainsi $\fka^{ii}_x\coloneqq \fka^{ii}\otimes_{R^+}L_x\cong L_x$ . Le \cite[Lemma 10.74]{pas} donne $\varphi^{31}\in E^{13}$ et $\varphi^{12}\in E^{21}$ tels que :
\begin{itemize}
\itemb $(237)\implies \varphi^{31}\circ\colon E^{11}\xrightarrow{\sim} E^{13}\implies E^{13}_x\cong L_x$,
\itemb $(238)\implies \varphi^{12}\circ \colon E^{22}\xrightarrow{\sim} E^{21}\implies E^{21}_x\cong L_x$,
\itemb $(240)\implies E^{12}\cong \fka^{11}\implies E^{12}_x\cong L_x$,
\itemb $(241)\implies E^{31}\cong \fka^{33}\implies E^{31}_x\cong L_x$,
\itemb $(242)\implies E^{32}\cong \fka^{33}\implies E^{32}_x\cong L_x$.
\end{itemize}
Il reste à montrer que $E^{23}_x\cong L_x$. Or, $E^{23}$ est engendré sur $R^+$ par deux éléments (\cf \cite[]{pas}) $\varphi^{32}\coloneqq \varphi^{31}\circ\varphi^{12}$ et un autre élément $\beta$ (\cf \cite[Lemma 10.75]{pas}). On sait d'après les deux premiers points ci-dessus (correspondant à $(237)$ et $(238)$) que $\varphi^{32}$ est un isomorphisme et donc il reste à montrer qu'il existe une constante $u\in L_x$ telle que $\beta_x=u\varphi^{32}_x$.

D'après \cite[Lemma 10.92]{pas} il existe $c,d\in R^+$ tels que $c\beta =d\varphi^{32}$. Mais d'après \cite[Corollary 10.94]{pas}, $c,d\in \fkr$. Ainsi $c$ et $d$ sont des constantes non nulles modulo $\fkm_x$ et $\beta_x=u\varphi^{32}_x$ où $u\in L_x$ est défini par $uc\equiv d \mod \fkm_x$.

\end{proof}

Le corollaire suivant est alors immédiat :
\medskip
\begin{coro}\label{cor:cdfss}
Soit $F$ le corps des fractions d'un facteur de $R_{\CB,M}^{\lambda}$. Alors $E_{\CB,M}^{\lambda}\otimes_{R_{\CB,M}^{\lambda}}F$ est une algèbre semi-simple sur $F$.
\end{coro}
\medskip

\medskip
\begin{coro}\label{cor:pbatpt}
Soit $\fkm_x\subset R_{\CB,M}^{\lambda}$ un idéal maximal dont on note $L_x$ le corps résiduel, notons $R^+\coloneqq R_{\CB}^{\zeta_M^{\lambda}}$, alors 
$$
P_{\CB}\otimes_{R^+}L_x\cong \Pi_x'\otimes_{L_x}\bm_{\CB}(\Pi_x')^{*}.
$$
De plus, $\bm_{\CB}(\Pi_x')$ est un $L_x$-espace vectoriel de dimension fini.
\end{coro}
\medskip
\begin{proof}
Notons 
$$
E\coloneqq E_{\CB}^{\zeta^{\lambda}_M},\quad R^+\coloneqq R_{\CB}^{\zeta_M^{\lambda}},\quad E_x\coloneqq E\otimes_{R^+}L_x,\quad P_x\coloneqq P_{\CB}\otimes_{R^+}L_x.
$$
Par le lemme \ref{lem:quotp} on a 
\begin{equation}\label{eq:pitop}
\Pi_x'\cong P_{\CB}\otimes_{E}\bm_{\CB}(\Pi_x')\cong P_x\otimes_{E_x}\bm_{\CB}(\Pi_x').
\end{equation}
Notons que comme cet isomorphisme est un isomorphisme de $E_x$-modules $\Pi_x'\cong \Pi_x'\otimes_{L_x}E_x$. De plus, $\bm_{\CB}(\Pi_x')$ est un $L_x$-espace vectoriel de dimension finie par \cite[Corollary 6.7]{patu} car c'est un $R^+$-module de type fini. Comme $\bm_{\CB}(\Pi_x')$ est un $E_x$-module absolument irréductible par \cite[Corollary 1.19]{pas}, on a 
$$
\bm_{\CB}(\Pi_x')\otimes_{E_x}\bm_{\CB}(\Pi_x')^*\cong E_x.
$$
Il suffit alors d'appliquer $\otimes_{E_x}\bm_{\CB}(\Pi_x)^*$ pour obtenir le résultat, soit
$$
\Pi_x'\otimes_{L_x}\bm_{\CB}(\Pi_x)^*\cong  \Pi_x'\otimes_{E_x}\bm_{\CB}(\Pi_x')\cong P_x\otimes_{E_x}\bm_{\CB}(\Pi_x')\otimes_{E_x}\bm_{\CB}(\Pi_x')\cong P_x\otimes_{E_x}E_x\cong P_x.
$$
\end{proof}

\subsubsection{Structure de $\bm_{\CB,\et}$ comme $R_{\CB}^{\ps}$-module}\label{subsubsec:final}
Dans ce dernier numéro on finit la preuve du théorème \ref{thm:repprin}. On commence par un lemme :
\medskip
\begin{lemm}\label{lemm:priptto}
Soit $V,W$ deux $E_{\CB,M}^{\lambda}[\sG_{\Qp}]$-modules localement libres et de type fini sur $R_{\CB,M}^{\lambda}$. Supposons que pour tout idéal maximal $\fkm_x\subset R_{\CB,M}^{\lambda}$, dont on note $L_x$ le corps résiduel, on a un isomorphisme $\sG_{\Qp}$-équivariant
$$
V\otimes_{R_{\CB,M}^{\lambda}}L_x\cong W\otimes_{R_{\CB,M}^{\lambda}}L_x,
$$
et $V\otimes_{R_{\CB,M}^{\lambda}}L_x$ est irréductible, alors $V\cong W$.
\end{lemm}
\medskip
\begin{proof}
Comme $R_{\CB,M}^{\lambda}$ est un produit d'anneaux principaux, il suffit de se ramener à l'un de ses facteurs, que l'on note $R$ par un projecteur. Soit $F$ le corps des fractions de $R$, notons 
$$
E_R\coloneqq E_{\CB,M}^{\lambda}\otimes_{R_{\CB,M}^{\lambda}}R,\quad E_F\coloneqq E_R\otimes_{R}F,\quad V_F\coloneqq V\otimes_{R_{\CB,M}^{\lambda}}F,\quad W_F\coloneqq W\otimes_{R_{\CB,M}^{\lambda}}F.
$$
Les traces de $V$ et $W$ coïncident car leurs spécialisations en tout idéal maximal de $R$ coïncident (comme il y a une infinité d'idéaux premiers). De plus, $V_F$ et $W_F$ sont des $E_{F}[\sG_{\Qp}]$-modules irréductibles, sinon il existerait un ouvert $U\subset \Spec\, R$ tel que la spécialisation de $V$ et $W$ en les idéaux maximaux de $U$ seraient réductibles. Comme $E_F$ est une algèbre semi-simple d'après le corollaire \ref{cor:cdfss} et par le théorème de Brauer-Nesbitt, on obtient un isomorphisme $V_F\cong W_F$ de $E_F[\sG_{\Qp}]$-modules. En recopiant la preuve de \cite[Lemme 5.12]{codonifac} on en déduit un isomorphisme $V\cong W$ de $R[\sG_{\Qp}]$-modules, mais par construction, cet isomorphisme étendu à $F$ redonne l'isomorphisme précédent $V_F\cong W_F$ qui commute à l'action de $E_F$, on en déduit que l'isomorphisme $V\cong W$ commute à l'action de $E_R$ et donc que c'est un isomorphisme de $E_R[\sG_{\Qp}]$-modules.
\end{proof}
Notons
$$
\bm_{\CB,\et}^{\lambda}\coloneqq \Hom_G(P_{\CB},\rmH^{\lambda,+}_{\Qpbar})\otimes_{\rmO_L}L
$$
Remarquons que puisque $P_{\CB}$ est compact on a $\bm_{\CB,\et}^{\lambda}=\bm_{\CB,\et}^{\lambda,+}$. On a de plus
$$
\bm_{\CB,\et}^{\lambda}=\bigoplus_{M\in \phinl^n}\bm_{\CB,M}^{\lambda}\otimes_L \JL_M,\quad \bm_{\CB,M}^{\lambda}\coloneqq\Hom_{G\times \czG}(P_{\CB}\otimes_{\rmO_L}\JL_M,\rmH^{\lambda}_{\Qpbar})
$$
On achève maintenant la démonstration du théorème \ref{thm:repprin}. Notons que le centre de $\czG$ agit par le caractère $\zeta_M^{\lambda}$ sur $\rmH^1_{\Qpbar}[M]$, le centre de $G$ agit par le caractère inverse $\zeta\coloneqq(\zeta_M^{\lambda})^{-1}$ (\cf \cite[Remarque 5.23, (ii)]{codonifac}) ce qui fait de $\bm_{\CB,M}^{\lambda}$ un $R_{\CB}^{\ps,\zeta}[1/p]$-module par le théorème \ref{thm:patu}.
\medskip
\begin{theo}
Soit $\lambda \in P_+$ et $M\in \phinl^n$.
\begin{enumerate}
\item L'anneau $R_{\CB}^{\ps,\zeta_M^{\lambda}}$ agit sur $\bm_{\CB,M}^{\lambda}$ au travers de son quotient $R_{\CB,M}^{\lambda}$.
\item On a un isomorphisme de $E_{\CB}^{\zeta_M^{\lambda}}[\sG_{\Qp}]$-modules
$$
\bm_{\CB,M}^{\lambda}\cong \bm_{\CB}^{\zeta}(\bPi(\rho_{\CB,M}^{\lambda})')\otimes_{R_{\CB,M}^{\lambda}}\rho_{\CB,M}^{\lambda}\otimes_{R_{\CB,M}^{\lambda}}\check{R}_{\CB,M}^{\lambda}
$$
où $\check{R}_{\CB,M}^{\lambda}$ est le $L$-dual continu de $R_{\CB,M}^{\lambda}$.
\end{enumerate}
\end{theo}
\medskip
\begin{proof}
La preuve se fait exactement comme pour la preuve de \cite[Théorème 5.24]{codonifac}, donc on rappelle les éléments principaux.

Pour le premier point, notons $I\coloneqq \ker(R_{\CB}^{\ps,\zeta_M^{\lambda}}[1/p]\rightarrow R_{\CB,M}^{\lambda})$, on veut montrer que $I$ annule $\bm_{\CB,M}^{\lambda}$. Comme $p>3$, on peut utiliser le théorème \ref{thm:patu} et considérer $J\coloneqq I\cap Z_{\CB}^{\zeta_M^{\lambda}}$. En vertu de la proposition \ref{prop:injco}, il suffit de montrer que $J$ tue $\Hom_G(P_{\CB},(\wh{\LL}^{\lambda}_M)')$. Soit $Y$ le complété $p$-adique d'un $\rmO_L$-réseau $G$-stable de $\LL^{\lambda}_M$. Par compacité de $P_{\CB}$ il suffit alors de montrer que pour tout $n\geqslant 1$, $J$ annule
$$
\Hom_G(P_{\CB},(Y/p^nY)')
$$
On a un morphisme naturel d'image dense $Y\rightarrow Y_{\CB}$ où $Y_{\CB}$ est le complété $\CB$-adique du réseau fixé de $\LL^{\lambda}_M$ ce qui induit une injection $\Hom_G(P_{\CB},(Y_{\CB}/p^nY_{\CB})')\incl \Hom_G(P_{\CB},(Y/p^nY)')$. Il reste ensuite à montrer que cette injection est un isomorphisme et que $J$ annule $\Hom_G(P_{\CB},(Y_{\CB}/p^nY_{\CB})')$ ; la preuve est \emph{verbatim} la même que la preuve du premier point de \cite[Théorème 5.24]{codonifac}.

Pour le second point on pose $X\coloneqq \Spm(R_{\CB}^{\ps,\zeta_M^{\lambda}}[1/p])$ et $X_M^{\lambda}\coloneqq \Spm(R^{\lambda}_{\CB,M})$ et soit $\fkm_x\subset R_{\CB}^{\ps,\zeta_M^{\lambda}}[1/p]$ l'idéal maximal d'un point $x\in X$ dont on note $L_x$ le corps résiduel. On commence par calculer $(\bm_{\CB,M}^{\lambda})[\fkm_x]$ :
\begin{itemize}
\itemb si $x\notin X_M^{\lambda}$, alors les images des éléments non nuls de $\fkm_x$ par $R_{\CB}^{\ps,\zeta_M^{\lambda}}[1/p]\rightarrow R_{\CB,M}^{\lambda}$ sont inversibles et comme $(\bm_{\CB,M}^{\lambda})$ est un $R_{\CB,M}^{\lambda}$-module,  $(\bm_{\CB,M}^{\lambda})[\fkm_x]=0$,
\itemb d'après le corollaire \ref{cor:pbatpt} on a 
$$
(\bm_{\CB,M}^{\lambda})[\fkm_x]=\Hom_G(P_{\CB}\otimes_{R^+}L_x,\rmH^{\lambda}_{\Qpbar})\cong\Hom_G(\Pi_x'\otimes_{L_x} \bm_{\CB}(\Pi_x'),\rmH^{\lambda}_{\Qpbar})\cong \bm_{\CB}(\Pi_x') \otimes_L\rho_x.
$$
\end{itemize}
En résumé on a
$$
(\bm_{\CB,\et}^{\lambda}[M])[\fkm_x]\cong
\begin{cases}
\bm_{\CB}(\Pi_x') \otimes_{L_x}\rho_x& \text{ si }x\in X_M^{\lambda},\\
0 & \text{ si }x\notin X_M^{\lambda}.
\end{cases}
$$
Notons que $(\bm_{\CB,\et}^{\lambda}[M])[\fkm_x]$ est le $L_x$-dual continu de $\check{\bm}_{\CB,M}^{\lambda}/\fkm_x$. On calcule le $L$-dual continu de $\bn_{\CB,M}^{\lambda}\coloneqq \bm_{\CB}^{\zeta_M^{\lambda}}(\bPi(\rho_{\CB,M}^{\lambda})')\otimes_R\rho_{\CB,M}^{\lambda}\otimes_R\check{R}_{\CB,M}^{\lambda}$ pour obtenir 
$$
\check{\bn}_{\CB,M}\cong \Hom_{R_{\CB,M}^{\lambda}}(\bm_{\CB}^{\zeta_M^{\lambda}}(\bPi(\rho_{\CB,M}^{\lambda})')\otimes\rho_{\CB,M}^{\lambda},R_{\CB,M}^{\lambda}).
$$
Ainsi $\check{\bn}_{\CB,M}/\fkm_x\cong \Hom_{L_x}(\bm_{\CB}(\Pi'_x)\otimes_{L_x}\rho_x,L_x)\cong \check{\bm}_{\CB,M}^{\lambda}/\fkm_x$. Or par la proposition \ref{prop:checktypfin}, comme $\check{\bm}_{\CB,M}^{\lambda}$ est un $R_{\CB}^{\ps,\zeta_M^{\lambda}}[1/p]$-module par l'isomorphisme du théorème \ref{thm:patu} on en déduit que c'est un $R_{\CB,M}^{\lambda}$-module de type fini. On est donc dans la situation du lemme \ref{lemm:priptto} qui permet de conclure la preuve.
\end{proof}

\end{document}